\documentclass[review,12pt]{elsarticle}

\usepackage{setspace}
\usepackage{subcaption}
\usepackage{multirow}
\usepackage{multicol}
\usepackage{graphicx}
\usepackage{amsmath,amsfonts,amssymb,amscd,amsthm,xspace}
\usepackage{bm}
\usepackage{xfrac}
\usepackage[most]{tcolorbox}
\usepackage{float}

\usepackage{amsmath,esint}

\usepackage{lineno,hyperref}
\journal{International Journal of Heat and Mass Transfer}

\biboptions{sort&compress} %USE THIS
\usepackage{cleveref}
\begin{document}

\begin{frontmatter}

\title{Uncertainty Quantification in Three Dimensional Natural Convection using Polynomial Chaos Expansion and Deep Neural Networks}
	%Elsevier \LaTeX\ template\tnoteref{mytitlenote}}
%\tnotetext[mytitlenote]{Fully documented templates are available in the elsarticle package on \href{http://www.ctan.org/tex-archive/macros/latex/contrib/elsarticle}{CTAN}.}

%% Group authors per affiliation:
\author{Shantanu Shahane\fnref{Corresponding Author}}
\author{Narayana R. Aluru} \author{Surya Pratap Vanka}
\address{Department of Mechanical Science and Engineering\\
	University of Illinois at Urbana-Champaign \\
	Urbana, Illinois 61801}
\fntext[Corresponding Author]{Corresponding Author}
%%\author{Author Names}%\fnref{myfootnote}}
%%\address{Radarweg 29, Amsterdam}
%%\fntext[myfootnote]{Since 1880.}
%
%%% or include affiliations in footnotes:
%%\author[mymainaddress,mysecondaryaddress]{Elsevier Inc}
%%\ead[url]{www.elsevier.com}
%
%%\author[mysecondaryaddress]{Global Customer Service\corref{mycorrespondingauthor}}
%%\cortext[mycorrespondingauthor]{Corresponding author}
%%\ead{support@elsevier.com}
%%
%%\address[mymainaddress]{1600 John F Kennedy Boulevard, Philadelphia}
%%\address[mysecondaryaddress]{360 Park Avenue South, New York}
%
\begin{abstract}
This paper analyzes the effects of input uncertainties on the outputs of a three dimensional natural convection problem in a differentially heated cubical enclosure. Two different cases are considered for parameter uncertainty propagation and global sensitivity analysis. In case A, stochastic variation is introduced in the two non-dimensional parameters (Rayleigh and Prandtl numbers) with an assumption that the boundary temperature is uniform. Being a two dimensional stochastic problem, the polynomial chaos expansion (PCE) method is used as a surrogate model. Case B deals with non-uniform stochasticity in the boundary temperature. Instead of the traditional Gaussian process model with the Karhunen-Lo$\grave{e}$ve expansion, a novel approach is successfully implemented to model uncertainty in the boundary condition. The boundary is divided into multiple domains and the temperature imposed on each domain is assumed to be an independent and identically distributed (i.i.d) random variable. Deep neural networks are trained with the boundary temperatures as inputs and Nusselt number, internal temperature or velocities as outputs. The number of domains which is essentially the stochastic dimension is 4, 8, 16 or 32. Rigorous training and testing process shows that the neural network is able to approximate the outputs to a reasonable accuracy. For a high stochastic dimension such as 32, it is computationally expensive to fit the PCE. This paper demonstrates a novel way of using the deep neural network as a surrogate modeling method for uncertainty quantification with the number of simulations much fewer than that required for fitting the PCE, thus, saving the computational cost.
\end{abstract}

\end{frontmatter}
\section{Introduction}
Flow due to natural convection has been studied extensively in the literature \cite{de1983natural,fusegi1991numerical,le1991accurate,hu2017flow,li2016three,shu2003local,yigit2015effects} since it is practically useful in cooling or heating systems for applications like electronics, nuclear reactors, computing servers etc \cite{icoz2004design,sharma2007conjugate,herrlin2006gravity}. \citet{de1983natural}, \citet{shu2003local} and \citet{le1991accurate} simulated natural convection in a two dimensional differentially heated square cavity with gravity in a direction orthogonal to the applied temperature difference. \citet{de1983natural} used the stream function-vorticity formulation for laminar flow. Accurate benchmark solution was obtained using mesh refinement and extrapolation. \citet{le1991accurate} used a pseudo-spectral algorithm combining spatial expressions of Chebyshev polynomial series with a finite difference time marching scheme. The results for Rayleigh number upto $10^8$ are presented. \citet{shu2003local} solved the same problem with local radial basis function based differential quadrature (RDF-DQ) method. This method is a mesh-free approach with the RBFs as test functions to estimate the derivatives at any node as a weighted sum of values at the neighboring nodes. They discussed the effects of the RBF shape parameter and its fine tuning to get accurate solutions. \citet{fusegi1991numerical} presented the results for three dimensional differentially heated cubical enclosure. They used the finite difference discretization with SIMPLE algorithm \cite{patankar1980numerical} for laminar flow at Rayleigh numbers in the range of $10^3$ to $10^6$. Rayleigh-B{\'e}nard is another class of natural convection problems in which the temperature difference is applied parallel to the direction of gravity with the lower wall heated and upper wall cooled. \citet{hu2017flow}, \citet{li2016three} and \citet{yigit2015effects} studied the Rayleigh-B{\'e}nard convection numerically in two and three dimensional cavities with cubical and cylindrical shapes having various aspect ratios.
\par The numerical and experimental results of natural convection are utilized for verification and validation of the numerical software packages. Due to inaccuracies in the measurement and control, the experimental results are prone to errors. These errors can be estimated by introducing stochastic variations in the inputs and propagating these to the outputs. It is possible that the stochastic mean of an output parameter is different from its value at the input mean. Deterministic simulations alone cannot estimate the shift of mean. Thus in the recent years, there has been a growing interest in analysis of the effects of stochastic variations in the inputs on the outputs. There are multiple examples in the literature in which the uncertainty propagation techniques are combined with the deterministic numerical simulations \cite{xiu2002modeling,xiu2003modeling,knio2001stochastic,le2004natural,venturi2012supercritical,ganapathysubramanian2007sparse,carnevale2013uncertainty,marepalli2014quantifying,fajraoui2017analyzing,fezi2017uncertainty,hosder2006non,kumar2016efficient}. It is popular to use the polynomial chaos expansion (PCE) for uncertainty propagation in which the output is approximated as a summation of polynomial basis which are functions of the stochastic inputs. The two main classes of methods to estimate the coefficients of the PCE are stochastic Galerkin projection \cite{xiu2002modeling,xiu2003modeling,knio2001stochastic,le2004natural} and collocation \cite{carnevale2013uncertainty,marepalli2014quantifying,ganapathysubramanian2007sparse}. Stochastic Galerkin method requires modification of the underlying deterministic code since it requires solution of a new set of equations and thus, is called as intrusive method. This becomes a significant additional effort of software development and it is difficult to couple with the legacy codes. Hence, recently non-intrusive stochastic collocation methods have gained popularity. The basic idea is to have multiple evaluations of the deterministic simulation at predefined collocation points which are samples from the underlying probability distribution function of the input parameters. The PCE coefficients are then estimated from the output values obtained from the deterministic solution at these input samples. The coefficients can be used for post-processing operations like output statistics estimation, response surface plotting and sensitivity analysis.
\par Uncertainty quantification (UQ) for various types of natural convection problems has been studied in the literature \cite{le2004natural,venturi2012supercritical,ganapathysubramanian2007sparse,fajraoui2017analyzing}. \citet{le2004natural} used the zero-Mach-number model to simulate natural convection in a two dimensional differentially heated square cavity with uncertainty in the cold wall temperature. The random component of the cold wall temperature is modeled using the Gaussian process with an auto-correlation function which is approximated by the truncated Karhunen-Lo$\grave{e}$ve (KL) expansion. PCE coefficients are estimated by the stochastic Galerkin projection. Output statistics for various values of the non-Boussinesq parameter $\epsilon$ are presented. \citet{ganapathysubramanian2007sparse} presented an adaptive refinement based approach for reducing the number of deterministic simulations in high dimensional stochastic simulations. The adaptive sampling method is applied to the two dimensional natural convection problem with random boundary condition which is modeled by the KL expansion method. \citet{venturi2012supercritical} studied the stability of the two dimensional Rayleigh-B{\'e}nard convection subject to stochastic boundary temperatures. The random boundary condition is assumed to be a non-uniform Gaussian random processes approximated by the KL expansion. It is found that the stochastic wall temperatures can extend the stability range of quasi-conduction states beyond the classical bifurcation point. \citet{fajraoui2017analyzing} analyzed the natural convection of porous media in a two dimensional differentially heated square cavity with uncertainty in the Rayleigh number, permeability anisotropy ratio, dispersion coefficients and heterogeneity variation. PCE method is used to estimate the statistics and sensitivity of the output parameters such as temperature and Nusselt number distributions. 
\par The PCE method is extremely useful for low dimensional uncertainty quantification. But at higher dimensions, it faces the problem known as `curse of dimensionality' i.e., for a linear increase in the stochastic dimensions, the number of samples grows exponentially. The Smolyak algorithm \cite{smolyak1963quadrature} addresses this problem to some extent by reducing the number of samples in high dimensions without compromising the interpolation accuracy. Even with the use of the Smolyak algorithm,  number of samples of the order of $10^3-10^4$ are required for five or more input dimensions. For instance, an eight and sixteen dimensional problem needs 3905 and 51073 samples respectively, for the accuracy level of five \cite{heiss2008likelihood}. Practically, it is computationally expensive to simulate the deterministic solution thousands of times. Thus, an alternate method is required for uncertainty propagation. The Monte Carlo method is a simple approach which approximates the statistics of the output by running the deterministic simulations at pseudo random samples of the inputs \cite{caflisch1998monte}. Since the error using the Monte Carlo method is $\mathcal{O}(\sfrac{1}{\sqrt{n}})$, the number of samples is practically too high which makes using the Monte Carlo method directly with the deterministic simulation difficult. Thus, it is popular to use a surrogate model which is trained and tested using deterministic simulations. A good surrogate model can be trained with a small number of deterministic simulations and its evaluation is cheap. A well tested surrogate model is further used to estimate the outputs at multiple sample inputs. Since the surrogate model evaluation is cheap, there is practically no limit on the number of input samples for the Monte Carlo method. Note that the PCE is also a surrogate model which is ideal for low stochastic dimensions with a possibility of direct estimation of the output statistics without the use of the Monte Carlo method.
\par In order to simulate a high dimensional stochastic problem, a neural network (NN) is used as a surrogate model. \citet{hornik1989multilayer} showed that multilayer feed forward networks are universal approximators i.e., with mild assumptions on the underlying function to be approximated, the network can achieve any desired degree of accuracy by choosing suitable number of neurons. The NNs can handle the `curse of dimensionality' by multiple nonlinear activation functions. In recent years, the NNs have been extremely popular in many fields of work as discussed by the review paper by \citet{schmidhuber2015deep}. Here, only the applications related to surrogate modeling for numerical simulations are discussed \cite{gholami2015simulation,sablani2005non,santra2009prediction,czel2014simultaneous,zhang2015machine,tripathy2018deep,khoo2017solving,zhang2018quantifying}. \citet{sablani2005non} used a NN as a surrogate model for inverse heat conduction problem of estimation of heat transfer coefficient from the temperature-time history at different locations. Since the NN is trained for the inverse problem directly using forward deterministic simulations, the estimation can be done non-iteratively. \citet{czel2014simultaneous} similarly used NN for non-iterative estimation of heat capacity and temperature dependent thermal conductivity using the experimental transient temperature histories. A radial basis function type NN is trained using the numerical solution of the direct heat conduction problem. Both the above publications show that computational time is saved by non-iterative estimation due to the NN surrogate model coupled with the forward numerical simulations. \citet{gholami2015simulation} trained a NN for a three dimensional two fluid flow in a $90^o$ curved channel and compared both the numerical simulations and NN predictions with experimental data. It is reported that the NN model is reasonably accurate and significantly faster compared to the full numerical simulation. \citet{tripathy2018deep} trained a NN to solve a steady state two dimensional diffusion process with spatially varying uncertainty in the diffusion coefficient. This uncertainty is modeled as a log normal random field with mean and covariance functions of the Gaussian random field which is approximated by the Karhunen-Lo$\grave{e}$ve (KL) expansion. Using the trained NN as a surrogate with diffusion coefficient as the input, statistics of the output parameter are estimated. \citet{zhang2018quantifying} have shown the utility of the physics informed neural networks (PINNs) for uncertainty quantification in direct and inverse stochastic problems. The basic idea of a PINN is to minimize the residual when the NN is substituted in the model differential equation together with the standard loss function of the NN. Automatic differentiation is used to estimate the residual. The paper claims that minimizing the residual along with loss function enhances the accuracy of the prediction.
\par This paper presents results of input uncertainty propagation for a three dimensional natural convection problem in a differentially heated cubical enclosure. Two different cases of input uncertainties are considered. Case A assumes that the boundary conditions are uniform. Thus, uncertainty is introduced in the two non-dimensional parameters (Rayleigh and Prandtl numbers). For this case, the polynomial chaos expansion (PCE) method is used as a surrogate model with stochastic collocation to estimate the PCE coefficients. Since this is a two dimensional stochastic problem, the number of samples required is small enough and thus, the estimation of the PCE coefficients is feasible. Case B deals with non-uniform stochastic boundary condition with deterministic material properties of the fluid. Since the temperature difference between the opposite walls drives the natural convection flow, the cold wall is held at a constant temperature and uncertainty is introduced in the hot wall temperature. Although the conventional method to deal with boundary condition uncertainties is to use the Gaussian process model with the Karhunen-Lo$\grave{e}$ve expansion \cite{ganapathysubramanian2007sparse,venturi2012supercritical,le2004natural,tripathy2018deep}, in this work, a novel approach is successfully implemented. The basic idea is to divide the hot wall into multiple domains and impose a temperature boundary condition on each domain. It is assumed that each domain temperature is an independent and identically distributed (i.i.d) random variable. Sets of forward deterministic simulations are used to train deep neural networks (DNNs) with the boundary temperatures as inputs and Nusselt number, internal temperature or velocities as outputs. The number of domains which is equal to the stochastic dimension ranges from 4 to 32. The DNN is successfully trained and tested with less number of samples compared to those required for PCE coefficient estimation for the high dimensional problem. This new approach to deal with stochastic boundary conditions with DNN as a surrogate model is found to be much better than those methods presented in the literature so far.
\section{Deterministic Problem Description}
%\subsection{Geometry and Boundary Conditions}
\subsection{Governing Equations}\label{Sec:Governing Equation}
In this paper, we consider the three dimensional natural convection in a differentially heated cube. A temperature gradient is applied on two opposite faces ($X=0$ and $X=1$) of a cube with sides of length $L$. The remaining four faces are thermally insulated. Gravity is imposed in the $Y$ direction which is orthogonal to the direction of the temperature gradient. Because of the thermal expansion of the fluid inside the cube due to the temperature variation, a buoyant force causes the lighter fluid to move upwards thus creating currents.
\par The flow field can be described by three dimensional incompressible Navier-Stokes and energy equations. For moderate density variations, natural convection can be modeled using the Boussinesq approximation. The system of equations is written in terms of non-dimensional variables as follows \cite{le1991accurate}:
\begin{equation} 
\nabla \cdot \textbf{u} =0
\label{Eq:Continuity}
\end{equation}
\begin{equation} 
\frac{\partial \textbf{u}}{\partial t} + (\textbf{u} \cdot \nabla)\textbf{u} = \frac{Pr}{Ra^{0.5}} \nabla^2 \textbf{u} - \nabla P - \hat{\textbf{g}} Pr \Theta
\label{Eq:Momentum}
\end{equation}
\begin{equation} 
\frac{\partial \Theta}{\partial t} + \nabla \cdot \textbf{u} \Theta = \frac{1}{Ra^{0.5}} \nabla^2 \Theta
\label{Eq:Energy}
\end{equation}
where, $\textbf{u}$ is the velocity vector, $\Theta$ is the temperature, $t$ is time, $\hat{\textbf{g}}$ is the unit vector in the direction of gravity, $P$ is the pressure, $Pr=\sfrac{\nu}{\alpha}$ is the Prandtl number and $Ra=\sfrac{g \beta \Delta T L^3}{\nu \alpha}$ is the Rayleigh number. Characteristic values for non-dimensionalization are as follows: velocity $u_{c}=\left(\sfrac{\alpha}{L}\right) Ra^{0.5}$, time $t_c=\left(\sfrac{L^2}{\alpha}\right) Ra^{-0.5}$, pressure $P_c=\rho u_c^2$ and $L$ is the cavity length. Non-dimensional temperature is defined as $\Theta=\sfrac{(T-T_m)}{(T_h-T_c)}$ where, $T_h$ and $T_c$ are hot and cold wall temperatures respectively and $T_m=\sfrac{(T_h+T_c)}{2}$ is the mean temperature.
\subsection{Solution Algorithm}
The governing equations (\ref{Eq:Continuity}-\ref{Eq:Energy}) are solved using the software OpenCast \cite{shahane2018finite} with finite volume method on a collocated grid. The fractional step method \cite{harlow1965numerical} is used to integrate the equations. An intermediate velocity field ($\textbf{u*}$) is first estimated by solving the modified momentum equation (\ref{Eq:Momentum_frac_step_u*}) without the pressure gradient. The diffusion term is discretized implicitly using second order Crank-Nicolson method whereas an explicit second order Adams-Bashforth is used for the convection term.
\begin{equation} 
\frac{\textbf{u*} - \textbf{u}^n}{\Delta t} = -Conv(\textbf{u}^n,\textbf{u}^{n-1}) + Diff(\textbf{u*},\textbf{u}^n) + Buoy(T^n)
\label{Eq:Momentum_frac_step_u*}
\end{equation}
The velocity correction equation (\ref{Eq:Velocity_correction}) is obtained by subtracting the modified momentum equation from the original equation.
\begin{equation} 
\textbf{u}^{n+1}=\textbf{u*}- (\nabla \Phi)^{n+1} \frac{\Delta t}{\rho}
\label{Eq:Velocity_correction}
\end{equation}
Imposing divergence free condition on the $(n+1)$ velocity field gives the Poisson equation for $\Phi$ (\ref{Eq:Pressure_poisson}).
\begin{equation} 
\nabla \cdot \left(\frac{\nabla \Phi }{\rho}\right)^{n+1} = \frac{\nabla \cdot \textbf{u*}}{\Delta t} 
\label{Eq:Pressure_poisson}
\end{equation}
The overall algorithm for marching from time step $n$ to $n+1$ can be summarized as follows:
 \begin{enumerate}
 	\item Solve for $\textbf{u*}$ using \cref{Eq:Momentum_frac_step_u*}
 	\item Solve the Poisson equation for $\Phi$ (\ref{Eq:Pressure_poisson}) iteratively to estimate $\Phi^{n+1}$
 	\item Correct the velocities ($\textbf{u}^{n+1}$) using \cref{Eq:Velocity_correction}
 	\item Solve the energy equation~(\ref{Eq:Energy})
 \end{enumerate}
Note that the pressure $P$ can be estimated from $\Phi$ if required:
\begin{equation} 
P = \Phi - \mu \frac{\Delta t}{\rho} \nabla^2 \Phi
\label{Eq:Pressure_Phi_relation}
\end{equation}
%\section{Uncertainty Quantification}
\section{Grid Independence Study and Verification}
\begin{figure}[H]
	\centering
	\begin{subfigure}[t]{0.49\textwidth}
		\includegraphics[width=\textwidth]{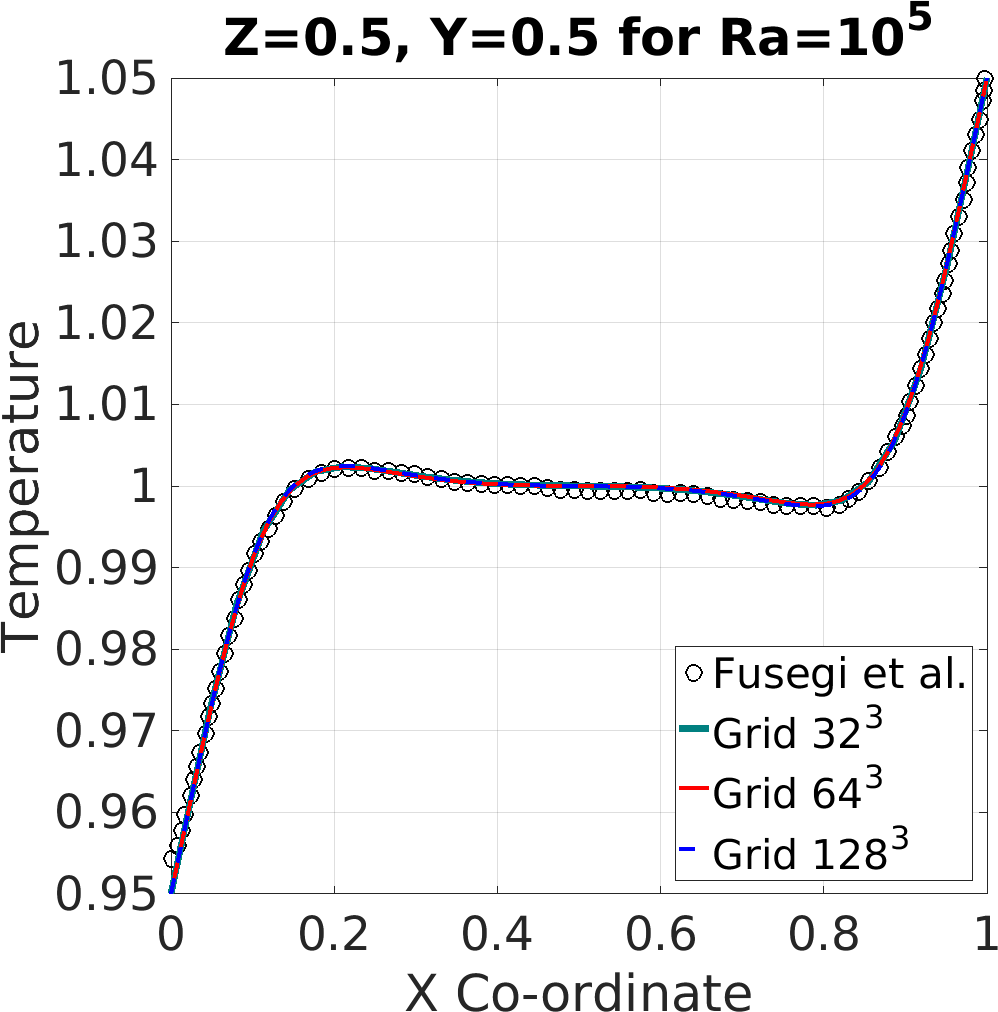}
		\caption{Temperature vs X: $Ra=10^5$}
		\label{Fig:Ra_1E5_Verification_Temp}
	\end{subfigure}
	\begin{subfigure}[t]{0.49\textwidth}
		\includegraphics[width=\textwidth]{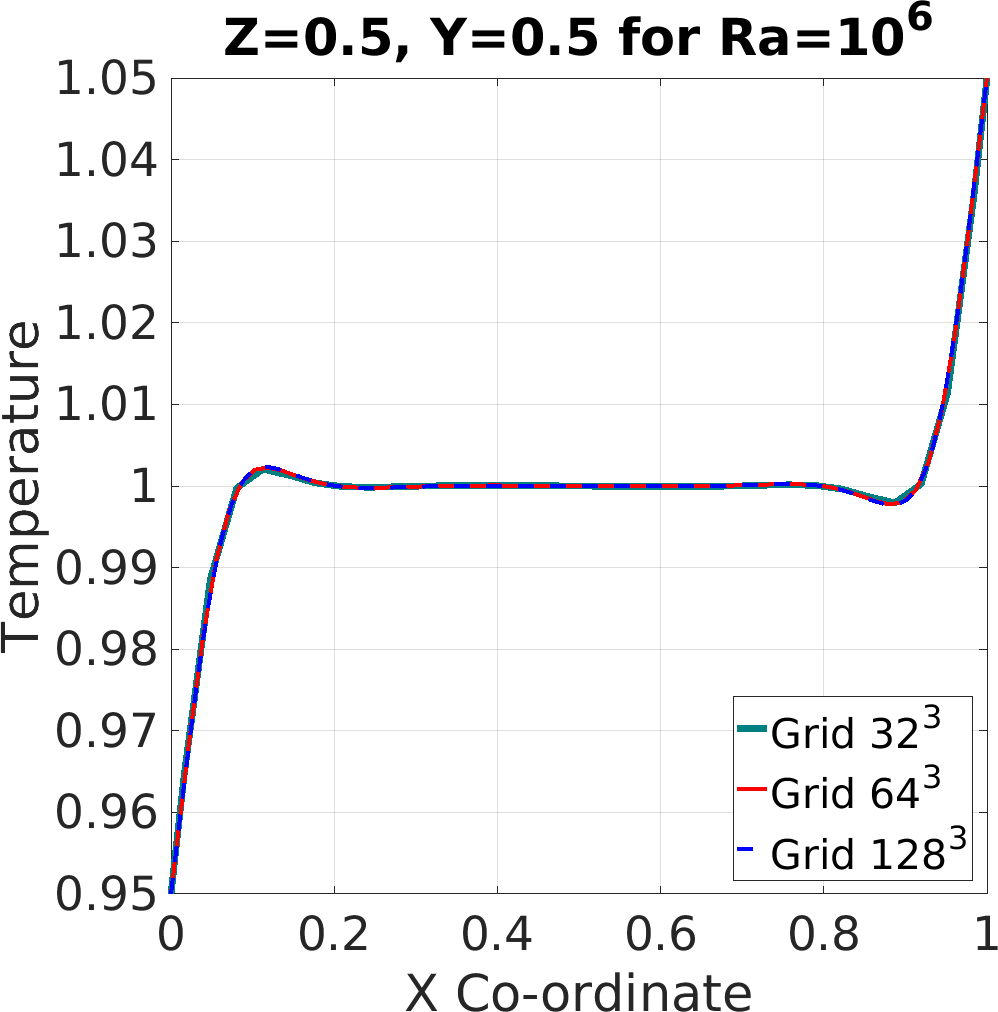}
		\caption{Temperature vs X: $Ra=10^6$}
		\label{Fig:Ra_1E6_Verification_Temp}
	\end{subfigure}
\end{figure}
\begin{figure}[H]\ContinuedFloat
	\centering
	\begin{subfigure}[t]{0.49\textwidth}
		\includegraphics[width=\textwidth]{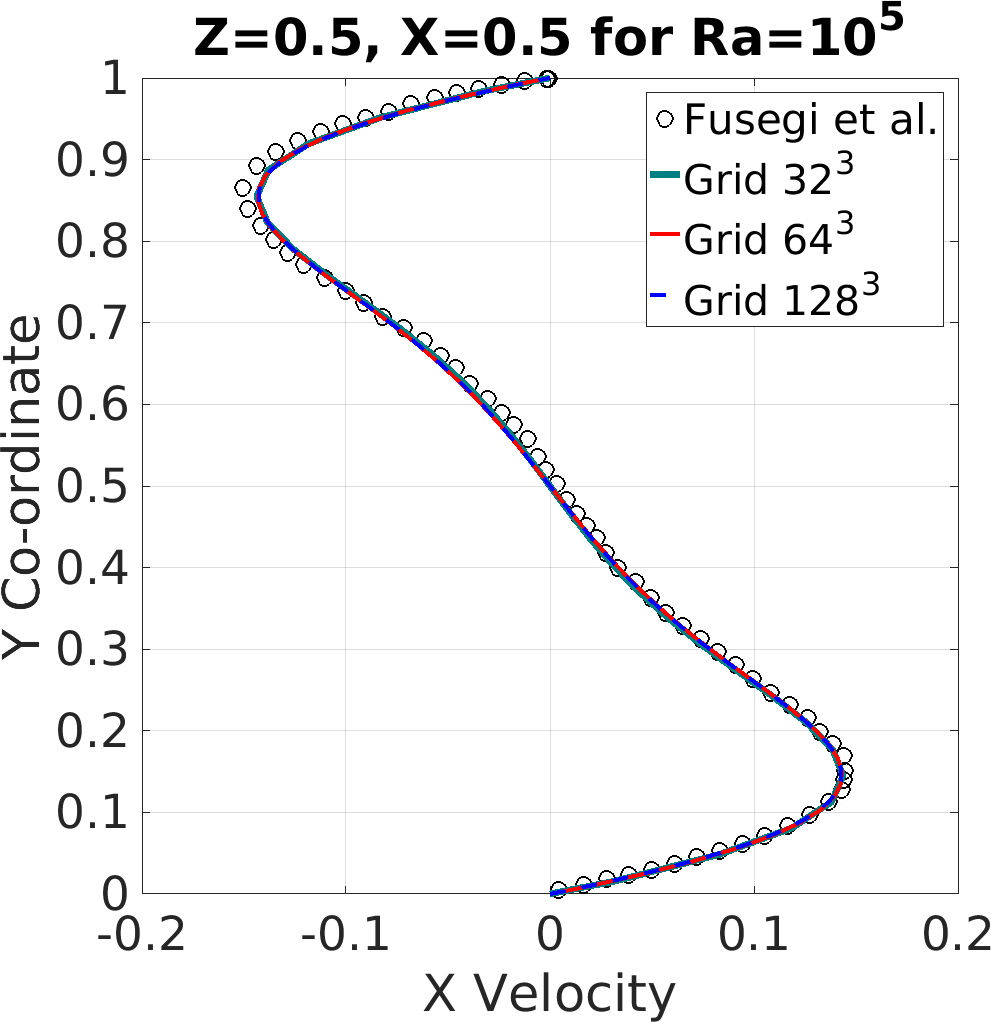}
		\caption{Y vs X Velocity: $Ra=10^5$}
		\label{Fig:Ra_1E5_Verification_Xvel}
	\end{subfigure}	
	\begin{subfigure}[t]{0.49\textwidth}
		\includegraphics[width=\textwidth]{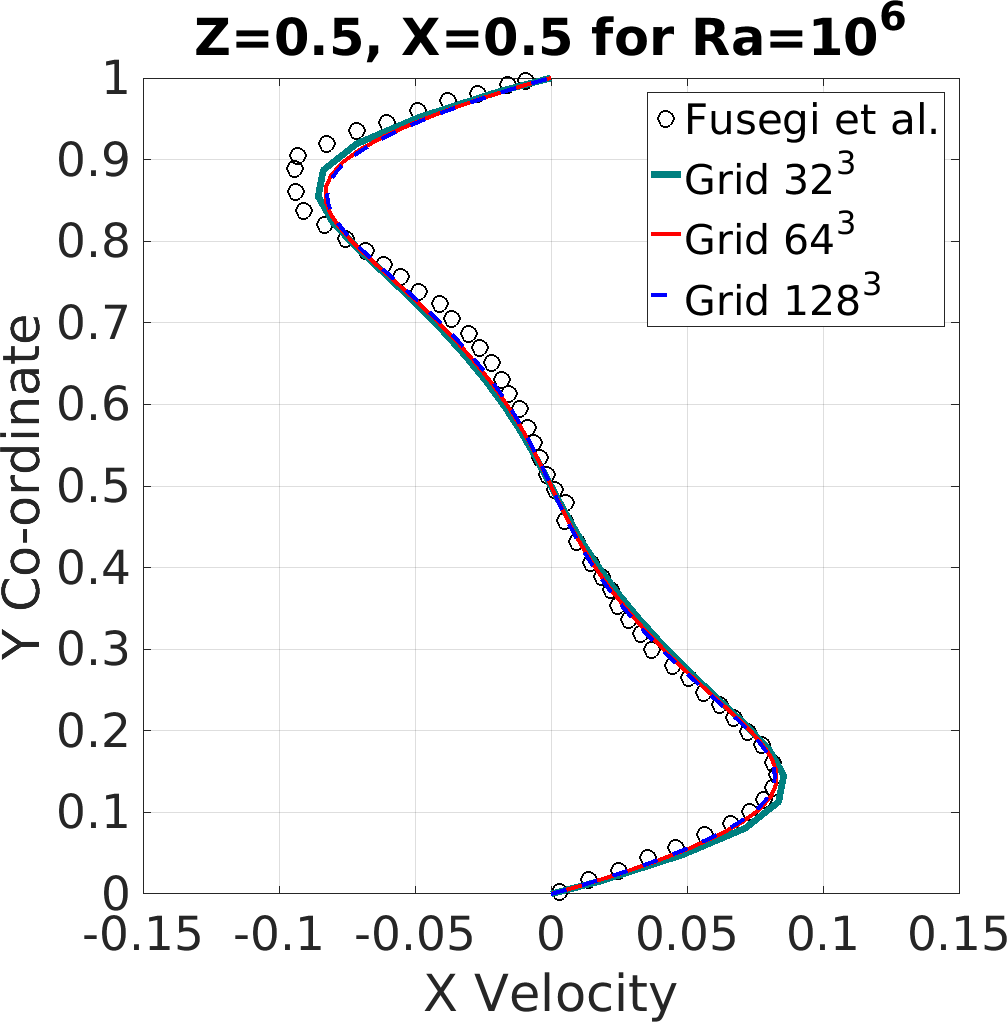}
		\caption{Y vs X Velocity: $Ra=10^6$}
		\label{Fig:Ra_1E6_Verification_Xvel}
	\end{subfigure}
\end{figure}
\begin{figure}[H]\ContinuedFloat
	\centering
	\begin{subfigure}[t]{0.49\textwidth}
		\includegraphics[width=\textwidth]{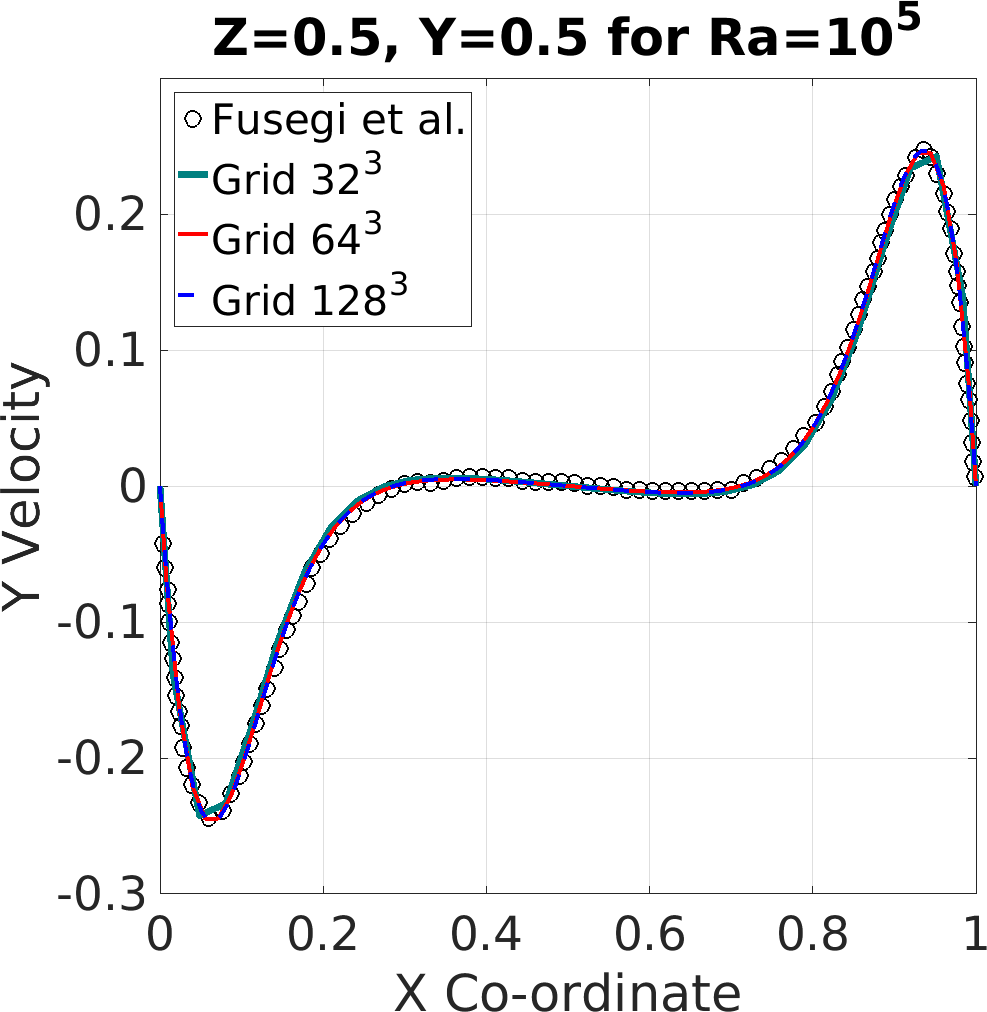}
		\caption{Y Velocity vs X: $Ra=10^5$}
		\label{Fig:Ra_1E5_Verification_Yvel}
	\end{subfigure}
	\begin{subfigure}[t]{0.49\textwidth}
		\includegraphics[width=\textwidth]{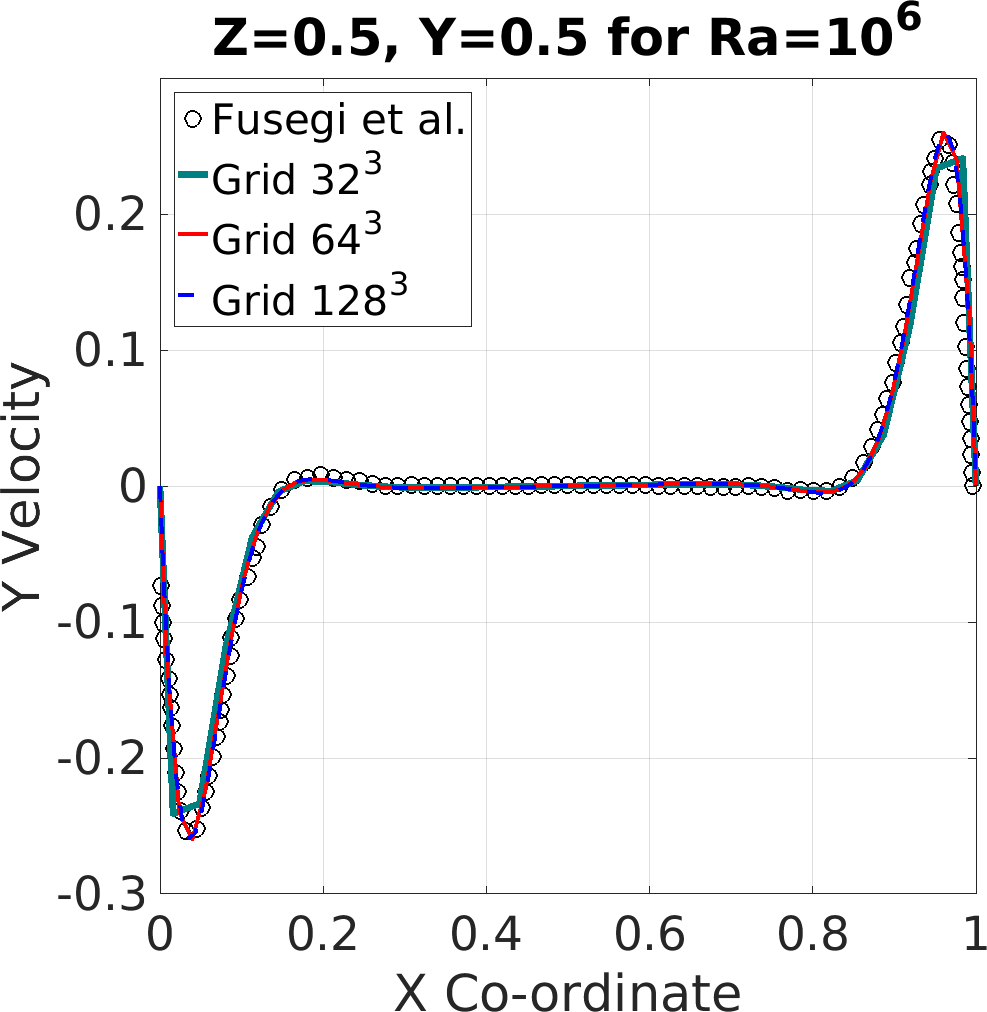}
		\caption{Y Velocity vs X: $Ra=10^6$}
		\label{Fig:Ra_1E6_Verification_Yvel}
	\end{subfigure}
	\caption{Temperature and Velocity Contours}
	\label{Fig:Nat_Conv_Verification}
\end{figure}
%\begin{figure}[h]
%	\centering
%	\begin{subfigure}[t]{0.32\textwidth}
%		\includegraphics[width=\textwidth]{Figures/Verification/Ra_1E5_T_X.png}
%		\caption{Temperature vs X}
%	\end{subfigure}
%	\begin{subfigure}[t]{0.32\textwidth}
%		\includegraphics[width=\textwidth]{Figures/Verification/Ra_1E5_U_Y.png}
%		\caption{Y vs X Velocity}
%	\end{subfigure}	
%	\begin{subfigure}[t]{0.32\textwidth}
%		\includegraphics[width=\textwidth]{Figures/Verification/Ra_1E5_V_X.png}
%		\caption{Y Velocity vs X}
%	\end{subfigure}
%	\caption{Temperature and Velocity Contours for $Ra=10^5$}
%	\label{Fig:Ra_1E5_Verification}
%\end{figure}
%\begin{figure}[h]
%	\centering
%	\begin{subfigure}[t]{0.32\textwidth}
%		\includegraphics[width=\textwidth]{Figures/Verification/Ra_1E6_T_X.png}
%		\caption{Temperature vs X}
%	\end{subfigure}
%	\begin{subfigure}[t]{0.32\textwidth}
%		\includegraphics[width=\textwidth]{Figures/Verification/Ra_1E6_U_Y.png}
%		\caption{Y vs X Velocity}
%	\end{subfigure}	
%	\begin{subfigure}[t]{0.32\textwidth}
%		\includegraphics[width=\textwidth]{Figures/Verification/Ra_1E6_V_X.png}
%		\caption{Y Velocity vs X}
%	\end{subfigure}
%	\caption{Temperature and Velocity Contours for $Ra=10^6$}
%	\label{Fig:Ra_1E6_Verification}
%\end{figure}
The numerical results from the software are verified using the three dimensional natural convection simulations of \citet{fusegi1991numerical}. For Rayleigh number up to $10^6$, the flow remains laminar \cite{fusegi1991numerical}. Hence Rayleigh numbers of $10^5$ and $10^6$ are used for verification. The cube is meshed with three different grid sizes with $32^3$, $64^3$ and $128^3$ structured hexahedrons. Steady state solution is estimated by time marching. For all the scalar fields $\phi$, when the non-dimensional steady state error computed over the entire domain defined as $\sfrac{max(|\phi^{new}-\phi^{old}|)}{max(|\phi^{new}|)}$ is less than $10^{-4}$, it is assumed that the steady state is reached. The temperatures at boundary faces $X=0$ and $X=1$ are set to 0.95 and 1.05 respectively.
\par The temperature and velocities are plotted along a centerline for both the Rayleigh numbers (\cref{Fig:Nat_Conv_Verification}). Note that the characteristic velocity used by \citet{fusegi1991numerical} is different compared to the discussion in section~\ref{Sec:Governing Equation}. Hence only for verification, the velocities are scaled by $u_{c}=\sqrt{g \beta L (T_h - T_c)}$.  In each figure, estimates from the three grid levels computed by OpenCast \cite{shahane2018finite} are superimposed with the results from \citet{fusegi1991numerical} whenever available. The numerical estimates from OpenCast match well with the published results thus, verifying the code. For the $Ra=10^5$ case (\cref{Fig:Ra_1E5_Verification_Temp,Fig:Ra_1E5_Verification_Xvel,Fig:Ra_1E5_Verification_Yvel}), it can be seen that all the three grid results from OpenCast overlap with each other. For $Ra=10^6$ (\cref{Fig:Ra_1E6_Verification_Temp,Fig:Ra_1E6_Verification_Xvel,Fig:Ra_1E6_Verification_Yvel}), the coarsest grid ($32^3$) plot is slightly off but the remaining two finer grid plots overlap. This shows that grid independence is achieved and for all further computations, a grid of size $101^3$ is used.
\section{Parameter Uncertainty Quantification} \label{Sec:Parameter Uncertainty Quantification}
Engineering problems have uncertainties in input parameters due to difficulty in precise measurement and control. It is necessary to propagate the input stochasticity to the output. Let $w(\bm{x},\bm{\xi})$ be a function which maps inputs to an output. Here, the vector $\bm{x}$ denotes all the deterministic inputs whereas, the vector $\bm{\xi}$ denotes all the stochastic parameters. It is assumed that the stochastic variable follows a known probability distribution: $\bm{\xi} \sim  f(\bm{\xi})$. The aim is to estimate the stochastic mean of the output  $w(\bm{x},\bm{\xi})$ defined as:
\begin{equation}
w_f (\bm{x}) = \int w(\bm{x},\bm{\xi}) f(\bm{\xi}) d \bm{\xi} 
\label{Eq:stoch_mean_integral}
\end{equation}
Since the function $w$ cannot be expressed in a closed form, above integral has to be approximated numerically. Importance sampling based Monte Carlo approximation of the integral is given by \cite{caflisch1998monte}:
\begin{equation}
w_f (\bm{x}) \approx \frac{1}{n} \sum_{i=1}^{i=n} w(\bm{x},\bm{\xi_i}) \frac{f(\bm{\xi_i})}{p(\bm{\xi_i})}
\label{Eq:stoch_mean_integral_MC_1}
\end{equation}
where, $\bm{\xi_i}$ are $n$ samples drawn from the probability distribution $p(\bm{\xi})$. It is effective to set $p(\bm{\xi})$ as $f(\bm{\xi})$ to reduce variance i.e., $\bm{\xi_i} \sim  f(\bm{\xi})$ \cite{caflisch1998monte}. Thus, \cref{Eq:stoch_mean_integral_MC_1} is simplified to:
\begin{equation}
w_f (\bm{x}) \approx \frac{1}{n} \sum_{i=1}^{i=n} w(\bm{x},\bm{\xi_i})
\label{Eq:stoch_mean_integral_MC_2}
\end{equation}
\par The error in the integral estimated using the Monte Carlo method is $\mathcal{O}(\sfrac{1}{\sqrt{n}})$. Thus, the number of samples ($n$ in \cref{Eq:stoch_mean_integral_MC_2}) can be of the order of thousands or more depending on the integrand. It is practically not possible to have so many samples as each sample corresponds to one deterministic simulation which is typically computationally costly. Thus, a surrogate model which is trained and tested using deterministic simulations is further used to estimate the outputs at multiple sample inputs. Since the surrogate model evaluation is cheap, there is practically no limit on the number of input samples.
\par The input parameters affecting the simulation of natural convection are boundary temperature, domain length and material properties like viscosity, thermal diffusivity and coefficient of thermal expansion. Two separate cases of input uncertainties are analyzed in this paper. Case A assumes that the boundary temperature is uniform. Thus, the governing equations (\ref{Eq:Continuity}-\ref{Eq:Energy}) show that the physics can be parametrized just using Rayleigh and Prandtl numbers. The uniform boundary temperature, domain length and the material properties all are included in the Rayleigh and Prandtl numbers. Hence, from the perspective of parameter uncertainty propagation, this is a two dimensional problem. Case B considers the possibility of uncertainty in the non-uniform boundary temperature with deterministic material properties. Case A is analyzed using the polynomial chaos expansion (Section~\ref{Sec:Polynomial Chaos Expansion}) whereas, deep neural networks (Section~\ref{Sec:Deep Neural Networks}) are used as surrogate models for the case B. %One of the benefits of polynomial chaos method is that the output variable statistics and sensitivity can be estimated directly from the coefficients of the expansion without any numerical approximations \cite{marelli2014uqlab}. On the other hand, the Monte-Carlo method is used to estimate the statistics of the outputs once the neural network is trained.
\par Such an analysis is practically important as there are stochastic variations in the boundary conditions due to inaccuracy in measurement and control. The fluid material properties also vary stochastically due to the presence of impurities. Hence, the effect of these uncertainties on the temperature and velocity distribution and Nusselt number is studied in this work. The following sections summarize the surrogate modeling strategies.
\subsection{Case A: Polynomial Chaos Expansion} \label{Sec:Polynomial Chaos Expansion}
Literature on uncertainty quantification describes various methods to estimate the relationship between stochastic inputs and outputs. Most of these methods rely on the idea of expanding the outputs as a linear combination of polynomial basis functions in the stochastic dimension. Orthogonal polynomial is a popular choice as basis since orthogonality helps in convergence. \citet{xiu2002wiener} showed that the Wiener's polynomial chaos \cite{wiener1938homogeneous} with orthogonal polynomials of Askey family lead to optimal convergence of the interpolation error. They identified which orthogonal polynomial basis is suitable depending on the probability distribution function followed by the stochastic variable. For example the Hermite polynomials are orthogonal to each other when the standard normal distribution is used as weighting function. Thus, Hermite polynomial basis function is recommended when the stochastic variable follows normal distribution. A polynomial chaos series (\cref{Eq:Polynomial_chaos}) is used to expand a second order random field \cite{xiu2002wiener}. The series is truncated to order $n$ for all practical purposes.
\begin{equation} 
w(\bm{x},\bm{\xi}(\theta))= \sum_{i=0}^{\infty} w_i(\bm{x})\bm{\Psi_i}(\bm{\xi}(\theta))\approx\sum_{i=0}^{n} w_i(\bm{x})\bm{\Psi_i}(\bm{\xi}(\theta))
\label{Eq:Polynomial_chaos}
\end{equation}
where, $\bm{x}$ is the spatial variable, $\bm{\xi}=(\xi_1,\xi_2,\dots \xi_r)$ is the random variable vector, $\theta$ is an elementary event, $\bm{\Psi_i}$ is a multi-dimensional orthogonal polynomial of order $i$, $w$ is the output to be estimated and $w_i$ are the coefficients of the series. In this case, $\bm{\xi}$ is a two dimensional vector following normal distribution and $\bm{\Psi_i}$ is a two dimensional Hermite polynomial.
\par In this paper, stochastic collocation method is used to estimate the deterministic coefficients ($w_i$) of the polynomial chaos expansion. Collocation acts as a wrapper over the existing deterministic software since it is a non-intrusive method. Thus, no modification of the deterministic software is required. After deterministic simulations at $M$ sample points ($\bm{\xi}^m$), a constraint $w(\bm{x},\bm{\xi}^m)=w_{sim}(\bm{x},\bm{\xi}^m)$ is imposed. The left hand side is estimated from polynomial chaos expansion (\cref{Eq:Polynomial_chaos}) and right hand side from each deterministic simulation. These $M$ constraints can be written in the matrix vector form \cite{smith2013uncertainty}. For accuracy, it is recommended to have more samples than the number of basis functions ($M>n+1$) and thus the Vandermonde system (\cref{Eq:Collocation_vandermonde}) is overdetermined.
\begin{equation} 
\renewcommand\arraystretch{1.8}
\begin{bmatrix} \bm{\Psi_0}\left(\bm{\xi}^1\right) & \cdots & \bm{\Psi_n}\left(\bm{\xi}^1\right)  \\ \vdots  && \vdots \\ \bm{\Psi_0}\left(\bm{\xi}^M\right) & \cdots & \bm{\Psi_n}\left(\bm{\xi}^M\right) \end{bmatrix}
\begin{bmatrix}
w_0(\bm{x})\\
\vdots \\
w_n(\bm{x})
\end{bmatrix} = 
\begin{bmatrix}
w_{sim}\left(\bm{x},\bm{\xi}^1\right)\\
\vdots \\
w_{sim}\left(\bm{x},\bm{\xi}^M\right)
\end{bmatrix} 
\label{Eq:Collocation_vandermonde}
\end{equation}
\par Sampling strategy plays a vital role in the accuracy and stability of stochastic collocation. Uniformly distributed samples lead to highly oscillatory interpolation and hence, poor convergence. Thus, for one dimensional stochastic problems, the roots of the basis orthogonal polynomials is a popular choice of sample points ($\bm{\xi}^m$) \cite{smith2013uncertainty}. For higher stochastic dimensions, tensor product of the single dimensional samples can be used. The sample size grows exponentially with dimensions if tensor product is used. This is a problem as each sample corresponds to a deterministic simulation and thus, the cost of computations grows exponentially. An algorithm to reduce the number of samples in high dimensions without spoiling the interpolation accuracy was published by \citet{smolyak1963quadrature}. It is found that for stochastic dimensions less than three, the Smolyak algorithm is not effective in reducing the sample size \cite{heiss2008likelihood}. Thus, in this work, a tensor product of single dimensional samples (roots of the Hermite polynomials) are used \cite{heiss2008likelihood}. UQLab, a MATLAB based tool developed by \citet{marelli2014uqlab} is used for estimation of polynomial chaos coefficients and response surfaces. The Polynomial Chaos-Kriging module of UQLab is used as it is found to be more effective than the basic polynomial chaos method.
\subsection{Case B: Deep Neural Network}\label{Sec:Deep Neural Networks}
Practically, a heat exchanger setup with a closed loop feedback system is used to maintain the cold and hot wall temperatures. Due to the errors in measurement and control, there are stochastic variations in the set temperature. For a large wall, there would be multiple heat exchangers in contact with the wall. Since the objective is to maintain a uniform temperature, the design and operation of all the heat exchangers would be similar. Thus, it is safe to assume that the temperature achieved by each of them is a random variable following the normal distribution with mean as the expected temperature and error modeled as the standard deviation. It is also assumed that the heat exchangers are independently controlled. Thus, the set temperatures are independent and identically distributed (i.i.d) random variables.
\par
\begin{figure}[H] 
	\centering
	\begin{subfigure}[t]{0.32\textwidth}
		\includegraphics[width=\textwidth]{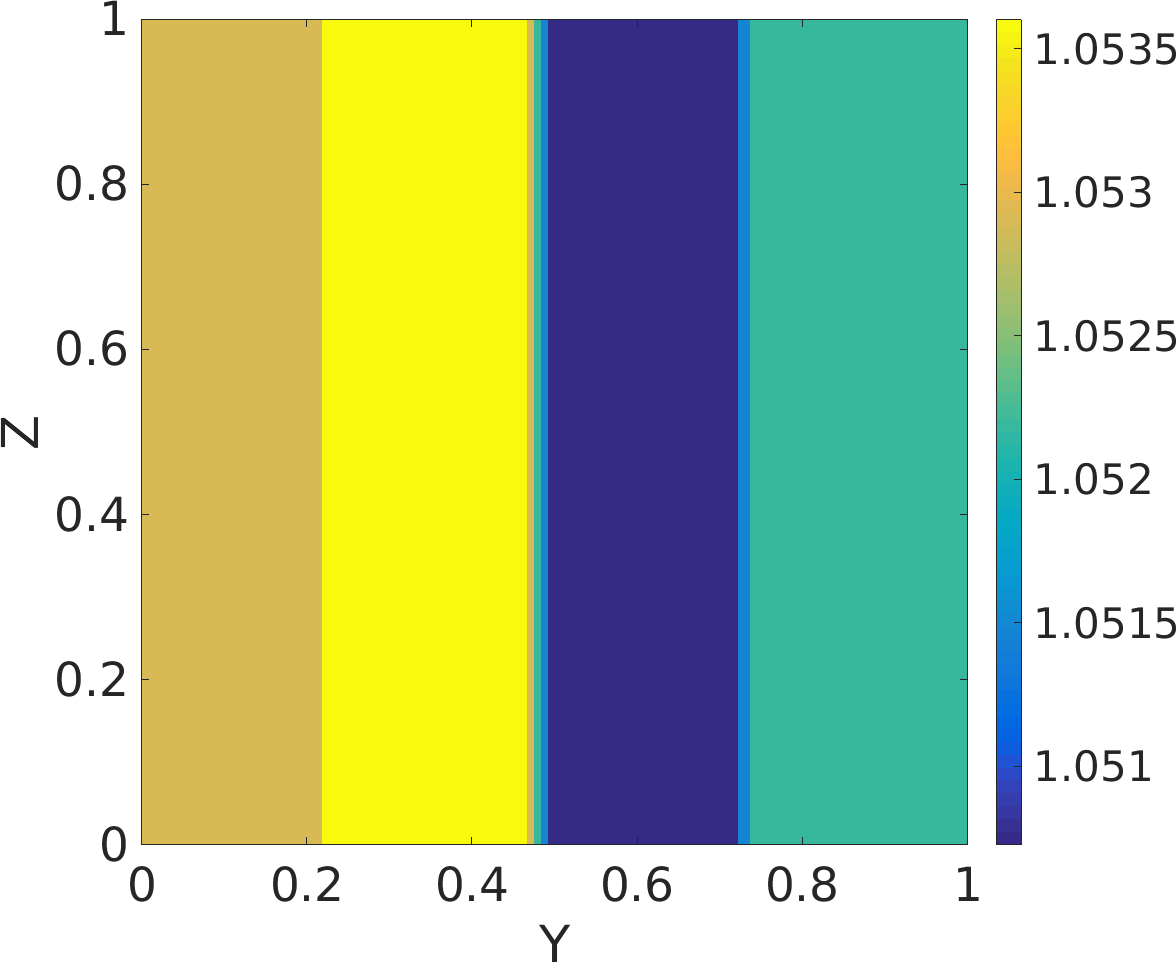}
		\caption{Number of Strips: 4}
		\label{Fig:Temperature_BC_4}
	\end{subfigure}
	\begin{subfigure}[t]{0.32\textwidth}
		\includegraphics[width=\textwidth]{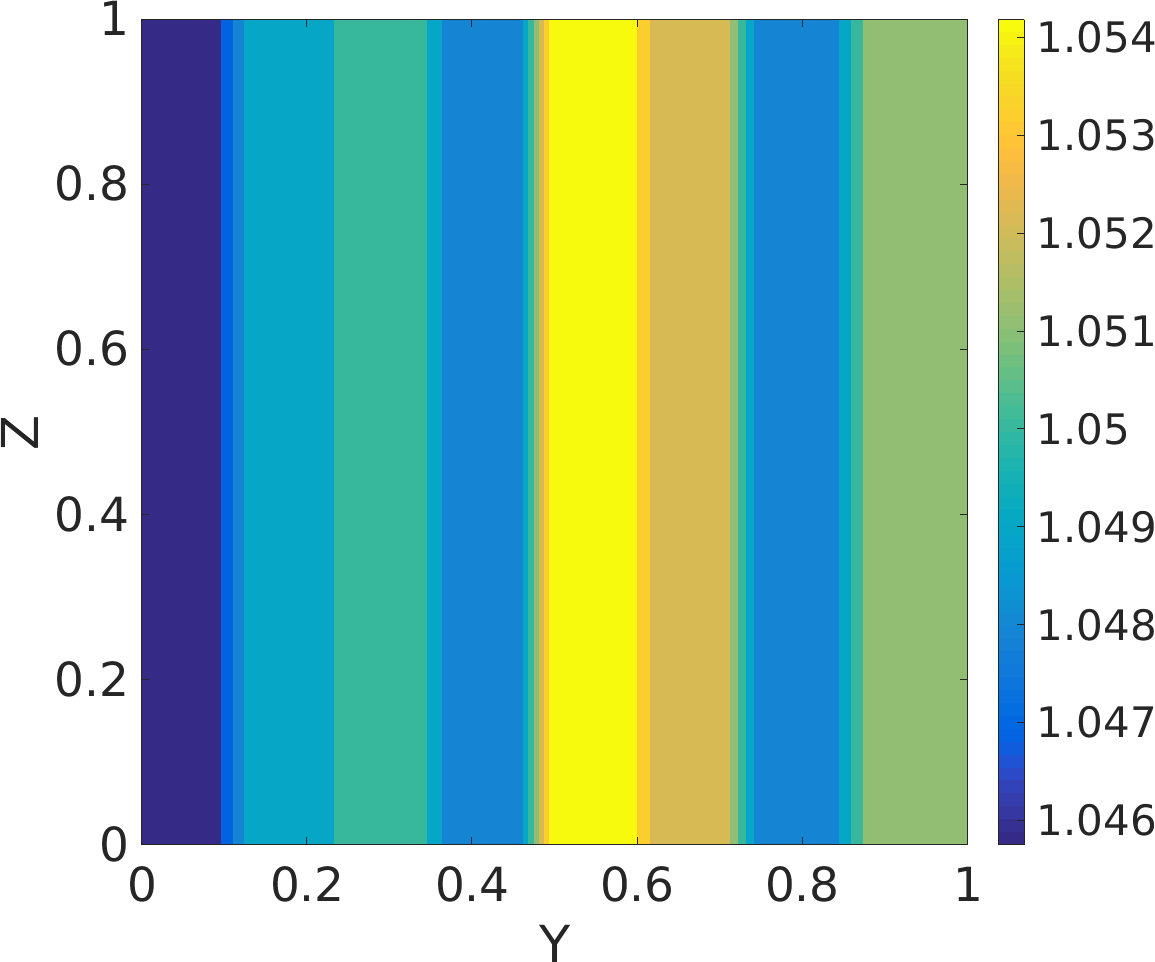}
		\caption{Number of Strips: 8}
		\label{Fig:Temperature_BC_8}
	\end{subfigure}
	\begin{subfigure}[t]{0.32\textwidth}
		\includegraphics[width=\textwidth]{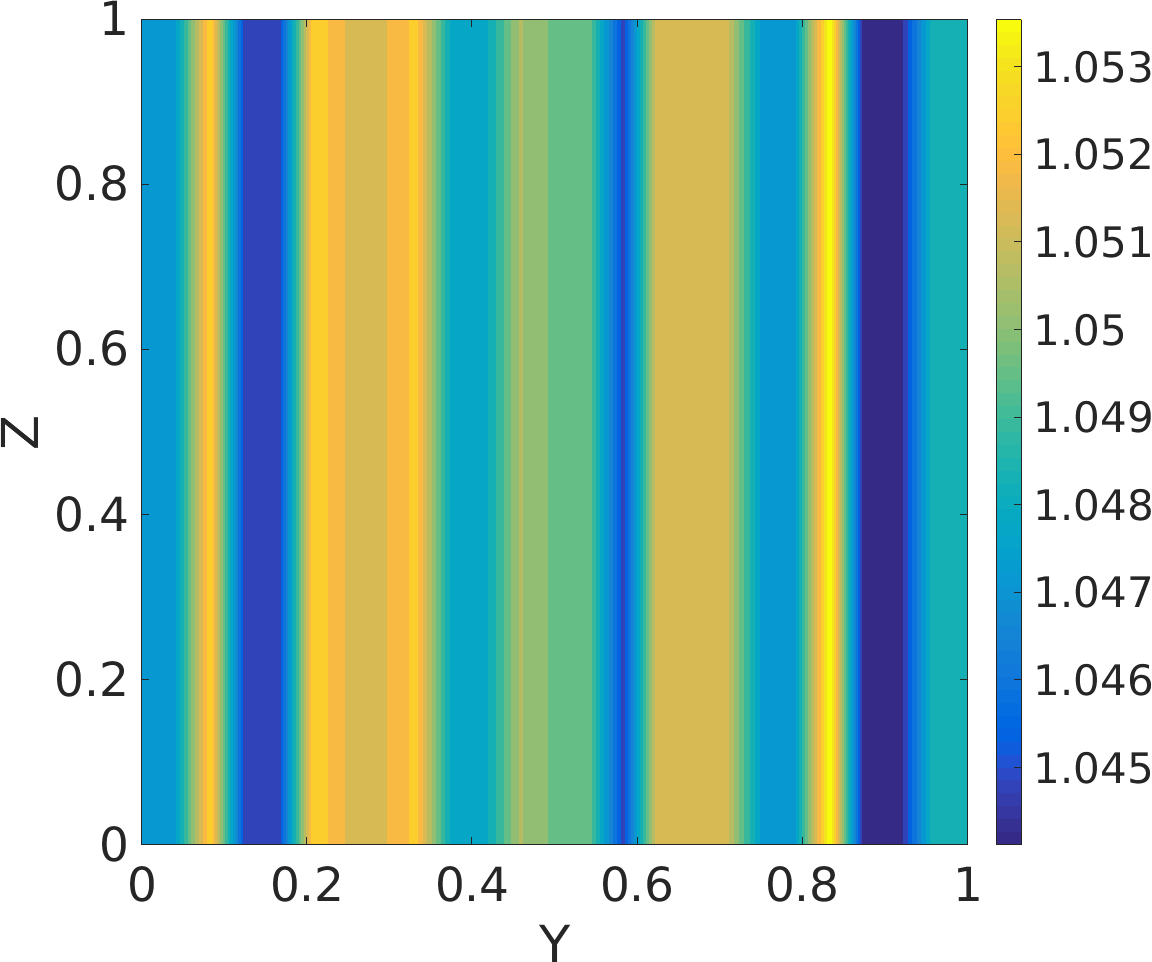}
		\caption{Number of Strips: 16}
		\label{Fig:Temperature_BC_16}
	\end{subfigure}
	\caption{Samples of Temperature Boundary Condition}
	\label{Fig:Temperature_BC}
\end{figure}

 In the literature, there are examples of using Gaussian process model with the Karhunen-Lo$\grave{e}$ve expansion to model uncertainty in the boundary condition \cite{ganapathysubramanian2007sparse,venturi2012supercritical,le2004natural}. In the present work, the wall is subdivided into multiple domains with different values of temperatures imposed as boundary conditions. From \cref{Fig:Nu_1E6}, it can be seen that the variation of Nusselt number is stronger in the direction of gravity ($Y$) compared to the orthogonal direction ($Z$). Thus, the wall is divided into strips along gravity ($Y$). \Cref{Fig:Temperature_BC} shows samples of the temperature boundary condition with 4, 8 and 16 number of strips. Each strip temperature is assumed to be an i.i.d. random variable following a normal distribution $(\mu=1.05,3 \sigma=0.01)$.
%Polynomial chaos expansion is a popular surrogate modeling strategy for low dimensional stochastic problems. The statistics of the output and the Sobol sensitivity indices can be computed from the polynomial coefficients directly without using any numerical approximation \cite{marelli2014uqlab}. Thus, the polynomial chaos method should be used whenever feasible. However, for higher stochastic dimensions, the number of samples needed to estimate the coefficients of the expansion is quite high even if Smolyak sparse grid is used. For instance, an eight and sixteen dimensional problem needs 3905 and 51073 samples respectively, for the accuracy level of five \cite{heiss2008likelihood}. Practically, it is computationally expensive to simulate the deterministic software thousands of times. Thus, a deep neural network is used as a surrogate model for case B.
%\subsubsection{Deep Neural Network Representation}
\par
\begin{figure}[H]
	\centering
	\begin{subfigure}[t]{0.35\textwidth}
		\includegraphics[width=\textwidth]{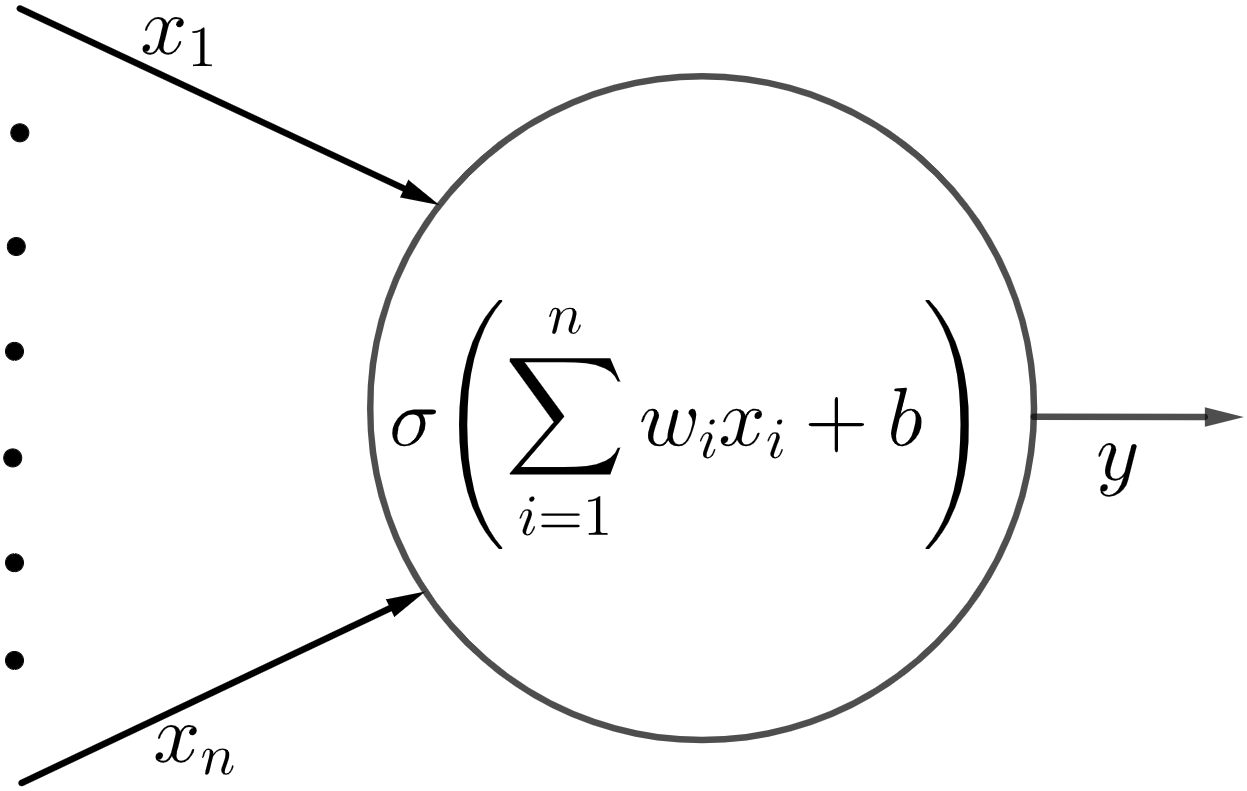}
		\caption{Single Neuron}
		\label{Fig:NN_schematics_1}
	\end{subfigure}
	\begin{subfigure}[t]{0.64\textwidth}
		\includegraphics[width=\textwidth]{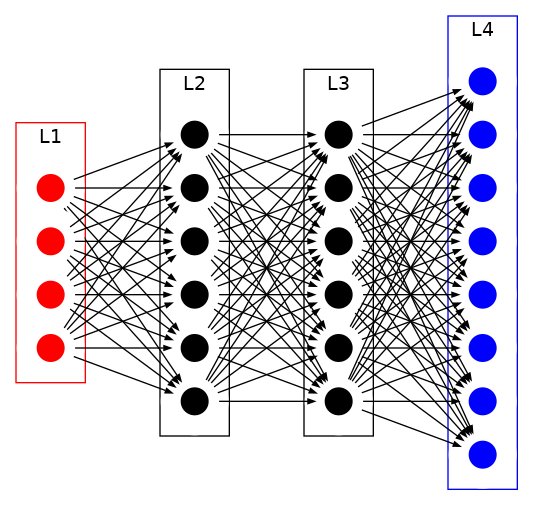}
		\caption{Deep Neural Network}
		\label{Fig:NN_schematics_2}
	\end{subfigure}
	\caption{Neural Network Schematics}
	\label{Fig:NN_schematics}
\end{figure}

 A neural network is a set of interconnected nodes such that the information flows from inputs to outputs. Each node is known as a neuron. \Cref{Fig:NN_schematics_1} shows a single neuron which has n scalar inputs $(x_1,x_2,...,x_n)$ and single output ($y$). Each neuron performs the following two operations in sequence:
\begin{enumerate}
	\item Linear transformation: $a=\sum_{i=1}^{n}{w_i x_i } + b$; where, $w_i$ are the weights and $b$ is a bias term
	\item Element-wise nonlinear transformation: $y=\sigma(a)$; where, $\sigma$ is the activation function
\end{enumerate}
A neural network is formed by stacking single neurons in a layer and connecting multiple layers as shown in \cref{Fig:NN_schematics_2}. It depicts an input (layer L1), an output (layer L4) and two hidden layers (layers L2 and L3). The arrows indicate the direction of information flow from input to output layer through the hidden layers. A deep neural network (DNN) is essentially a neural network with multiple hidden layers. Adding multiple hidden layers increases the nonlinearity of the network and thus, the network can approximate more complex functions successfully. The number of neurons in the input and output layers is specified by the problem definition whereas, number of hidden layers and neurons has to be fine tuned. 
The linear transformation followed by the nonlinear activation function of each neuron can be written in a single matrix vector equation:
\begin{equation}
\bm{y}^{(j)} = 
\begin{cases}
\bm{x} & \text{if } j=1\\
\sigma (\bm{W}^{(j)} \bm{y}^{(j-1)} + \bm{b}^{(j)}) & \forall j \in \{2,3,...,L\}
\end{cases}
\label{Eq:NN_forward_prop}
\end{equation}
where, $\bm{x} \in \mathbb{R}^{l_1}$ is the input vector, $\bm{W}^{(j)} \in \mathbb{R}^{l_j \times l_{j-1}}$ is the matrix of the weights, $\bm{y}^{(j)} \in \mathbb{R}^{l_j}$ is the activation produced by the $j^{th}$ layer and $\bm{b}^{(j)} \in \mathbb{R}^{l_j}$ is the bias. $L$ is the total number of layers including input, output and hidden layers. Number of neurons in the $j^{th}$ layer is denoted by $l_j$. For instance, in \cref{Fig:NN_schematics_2}, $L=4$, $l_1=4$, $l_4=8$ and $l_2=l_3=6$. Applying \cref{Eq:NN_forward_prop} sequentially starting from the input layer is known as forward propagation. This operation estimates the output vector ($\bm{y}^{(L)}$) from the input vector ($\bm{y}^{(1)}=\bm{x}$) if the weights and bias are known. Logistic sigmoid, hyperbolic tangent and rectified linear unit (ReLU) are some of the popular activation functions \cite{goodfellow2016deep}. In this work, the ReLU function defined by $\sigma(y)=max\{0,y\}$ is used for all the hidden layers. For output layer, in order to allow negative values, the identity function $\sigma(y)=y$ is used.
%\subsubsection{Deep Neural Network Training}
\par The process of estimation of weights and bias using a given set of inputs $(\bm{x}_i, 1\leq i \leq m)$ and the corresponding outputs $(\bm{z}_i,1\leq i \leq m)$ is known as training. For the given set of $m$ training samples, mean squared error between the neural network estimate ($\bm{\hat{z}}_i$) and the true value ($\bm{z}_i$) of the output is defined as the loss function: 
\begin{equation}
\mathcal{L} (\bm{W},\bm{b}; \bm{x}_i,\bm{z}_i) = \frac{1}{m} \sum_{i=1}^{i=m} ||\bm{z}_i - \bm{\hat{z}}_i||_2^2
\label{Eq:NN_loss_function}
\end{equation}
It is commonly seen that the neural network performs well on the training data but performs poorly on the unseen test data. This phenomenon is known as overfitting and is controlled with regularization. \citet{goodfellow2016deep} discuss various regularization methods in detail. Here, the L2 weight regularization with parameter $\lambda$ is used in which, the loss function (\cref{Eq:NN_loss_function}) is modified:
\begin{equation}
\mathcal{L} (\bm{W},\bm{b}; \bm{x}_i,\bm{z}_i) = \frac{1}{m} \sum_{i=1}^{i=m} ||\bm{z}_i - \bm{\hat{z}}_i||_2^2 + \lambda \sum_{l=1}^{L} ||\bm{W}^{(l)}||_2^2
\label{Eq:NN_reg_loss_function}
\end{equation}
The weights and bias which minimize the loss function are estimated using a numerical optimization algorithm. The gradient of the loss function with respect to the weights and bias is required in the optimization algorithms like gradient descent. The gradient is estimated by the backpropagation algorithm \cite{chauvin2013backpropagation}. The Adam optimizer is used in this work with the parameters $\beta_1$ and $\beta_2$ as suggested by \citet{kingma2014adam}. Other hyper-parameters like learning rate, regularization constant, number of hidden layers and number of hidden units are tuned using a validation set. All the implementation details are given in \cref{Sec:Deep Neural Network Training and Testing}. %Python library Tensorflow \cite{tensorflow2015-whitepaper} with a high level API Keras \cite{chollet2015keras} is used to implement the DNN.
\section{Uncertainty Propagation Results}
\subsection{Deterministic Results: Output Values at Input Mean}
\begin{figure}[H] 
	\centering
	\begin{subfigure}[t]{0.49\textwidth}
		\includegraphics[width=\textwidth]{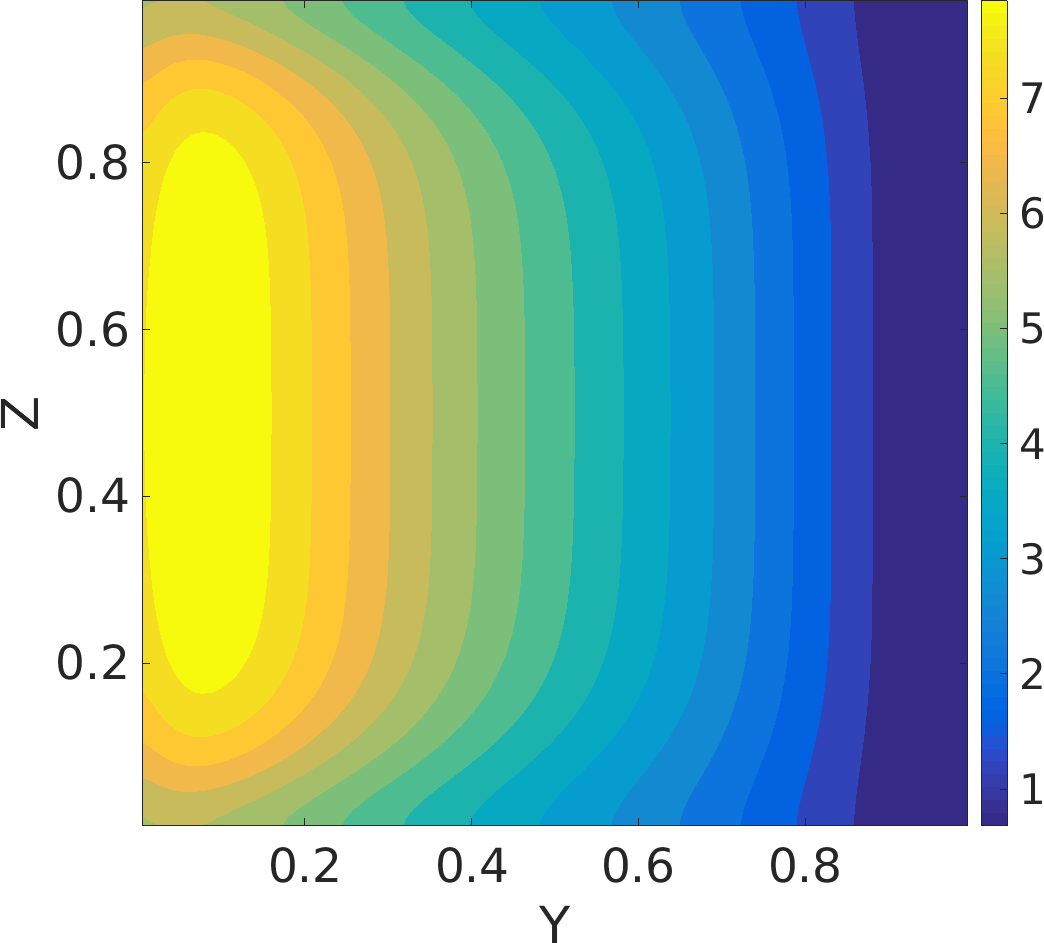}
		\caption{Hot Wall Nusselt Number}
		\label{Fig:Ra_1E5_value_at_input_mean_Nu}
	\end{subfigure}
	\begin{subfigure}[t]{0.49\textwidth}
		\includegraphics[width=\textwidth]{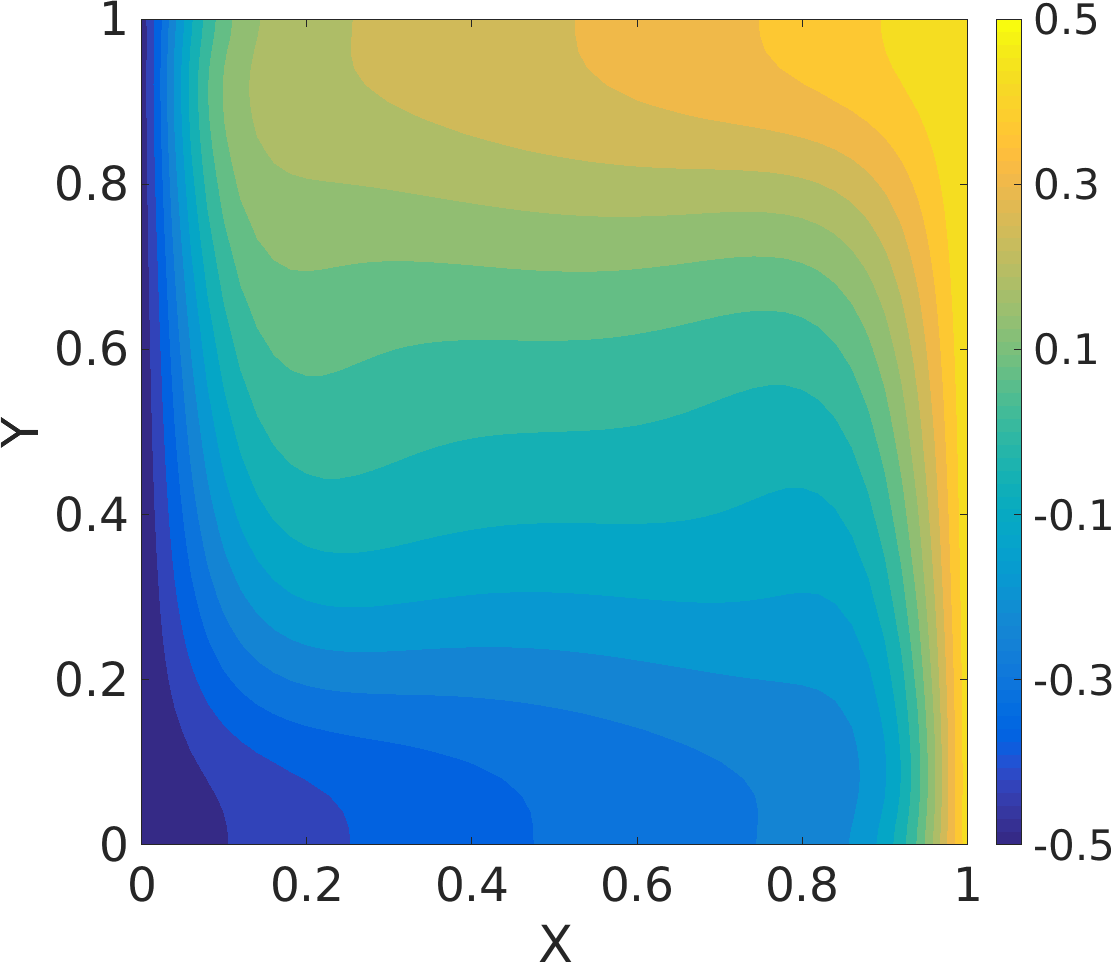}
		\caption{Temperature at $Z=0.5$}
		\label{Fig:Ra_1E5_value_at_input_mean_T}
	\end{subfigure}
	\begin{subfigure}[t]{0.49\textwidth}
		\includegraphics[width=\textwidth]{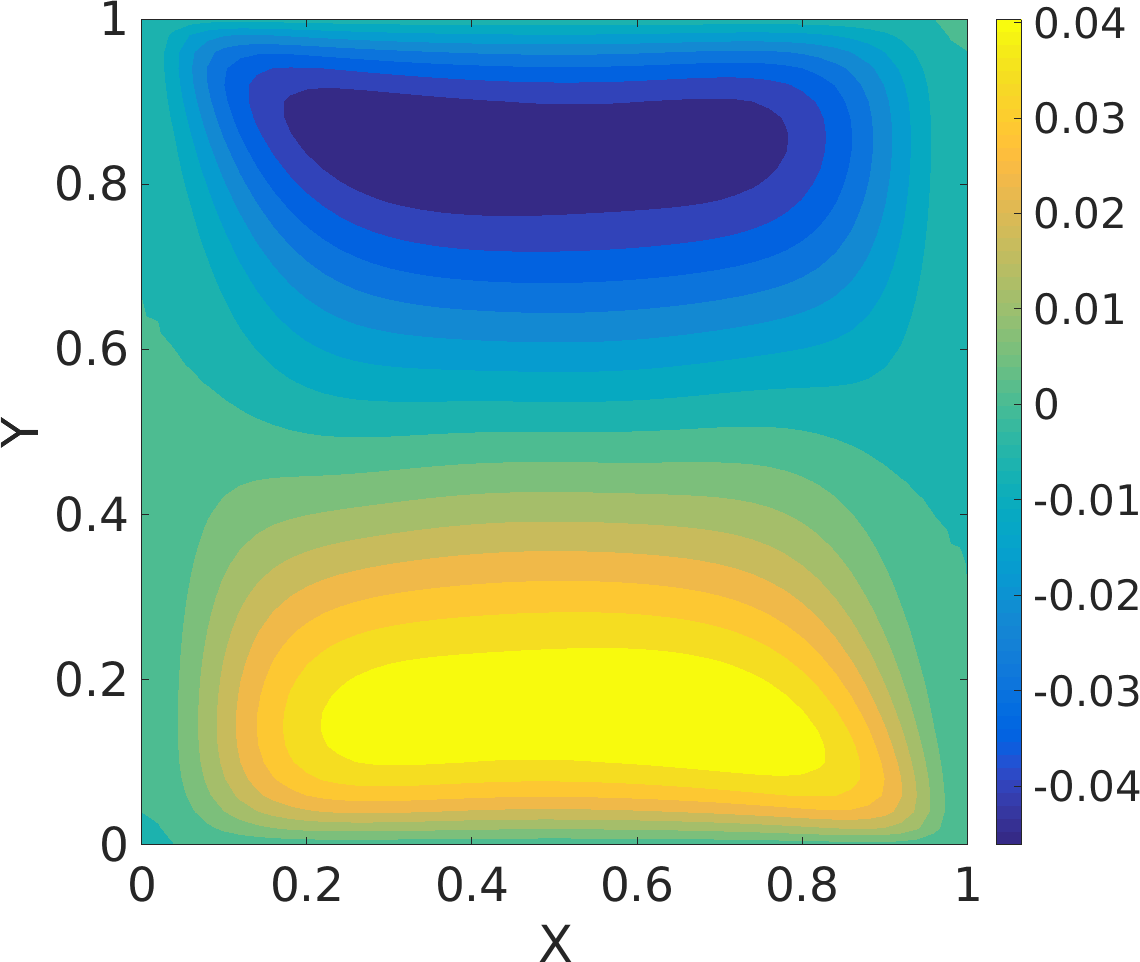}
		\caption{X-Velocity at $Z=0.5$}
		\label{Fig:Ra_1E5_value_at_input_mean_u}
	\end{subfigure}
	\begin{subfigure}[t]{0.49\textwidth}
		\includegraphics[width=\textwidth]{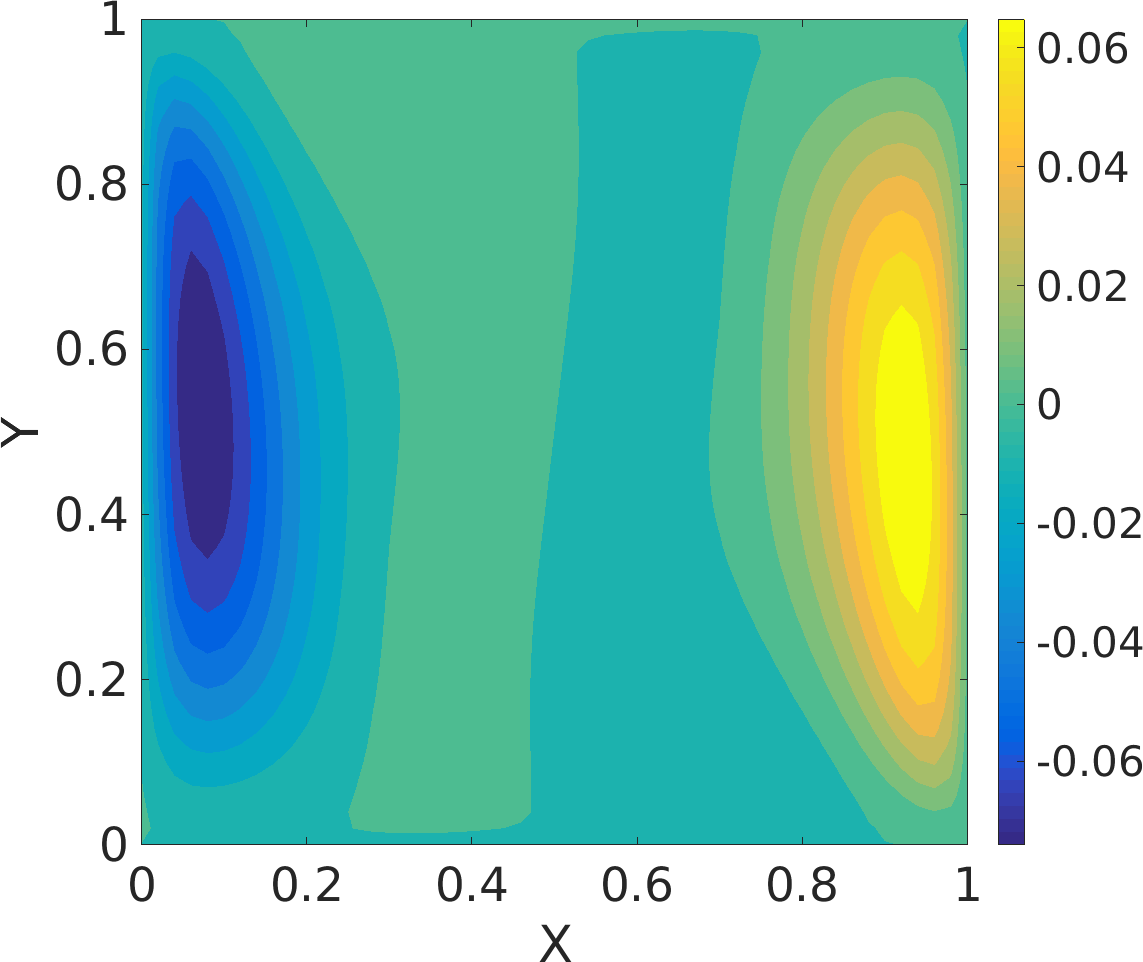}
		\caption{Y-Velocity at $Z=0.5$}
		\label{Fig:Ra_1E5_value_at_input_mean_v}
	\end{subfigure}
	\caption{Deterministic Results for $Ra=10^5$ and $Pr=7.5$}
	\label{Fig:Ra_1E5_value_at_input_mean}
\end{figure}

\begin{figure}[H] 
	\centering
	\begin{subfigure}[t]{0.49\textwidth}
		\includegraphics[width=\textwidth]{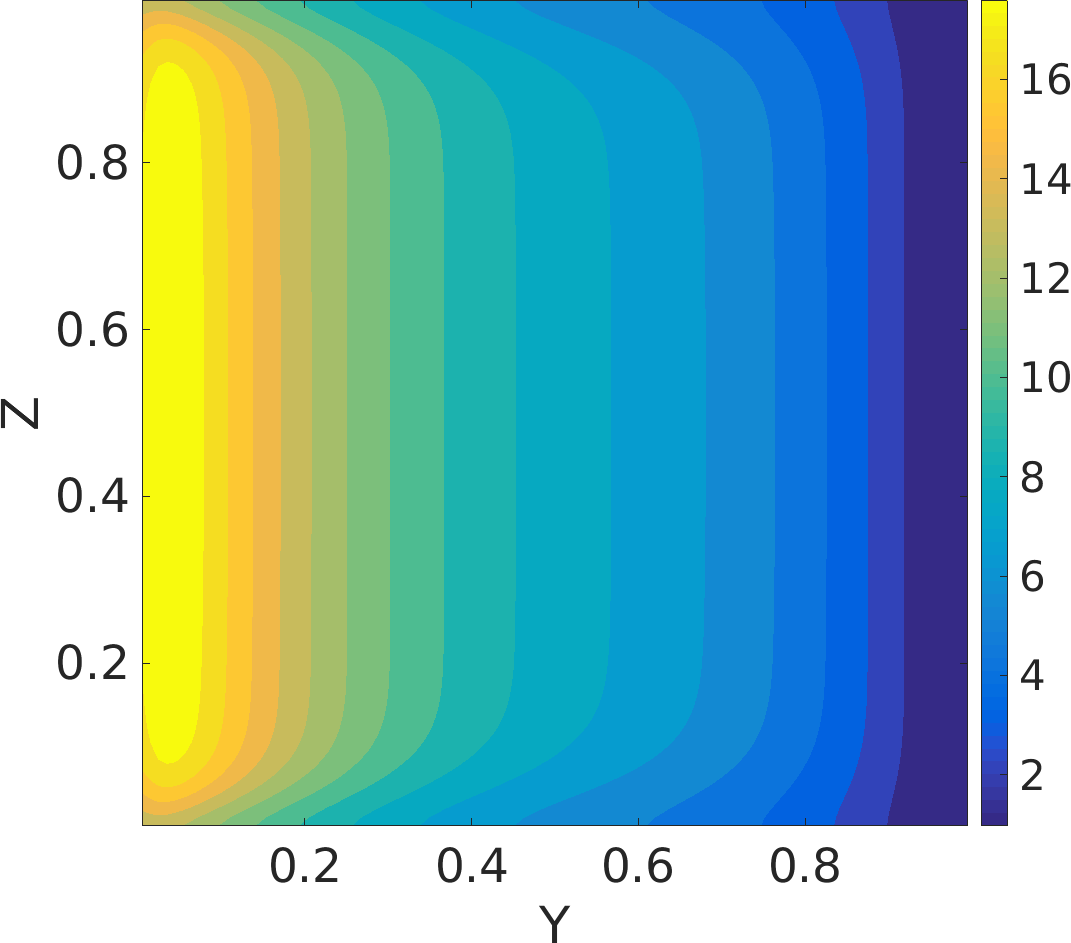}
		\caption{Hot Wall Nusselt Number}
		\label{Fig:Ra_1E6_value_at_input_mean_Nu}
	\end{subfigure}
	\begin{subfigure}[t]{0.49\textwidth}
		\includegraphics[width=\textwidth]{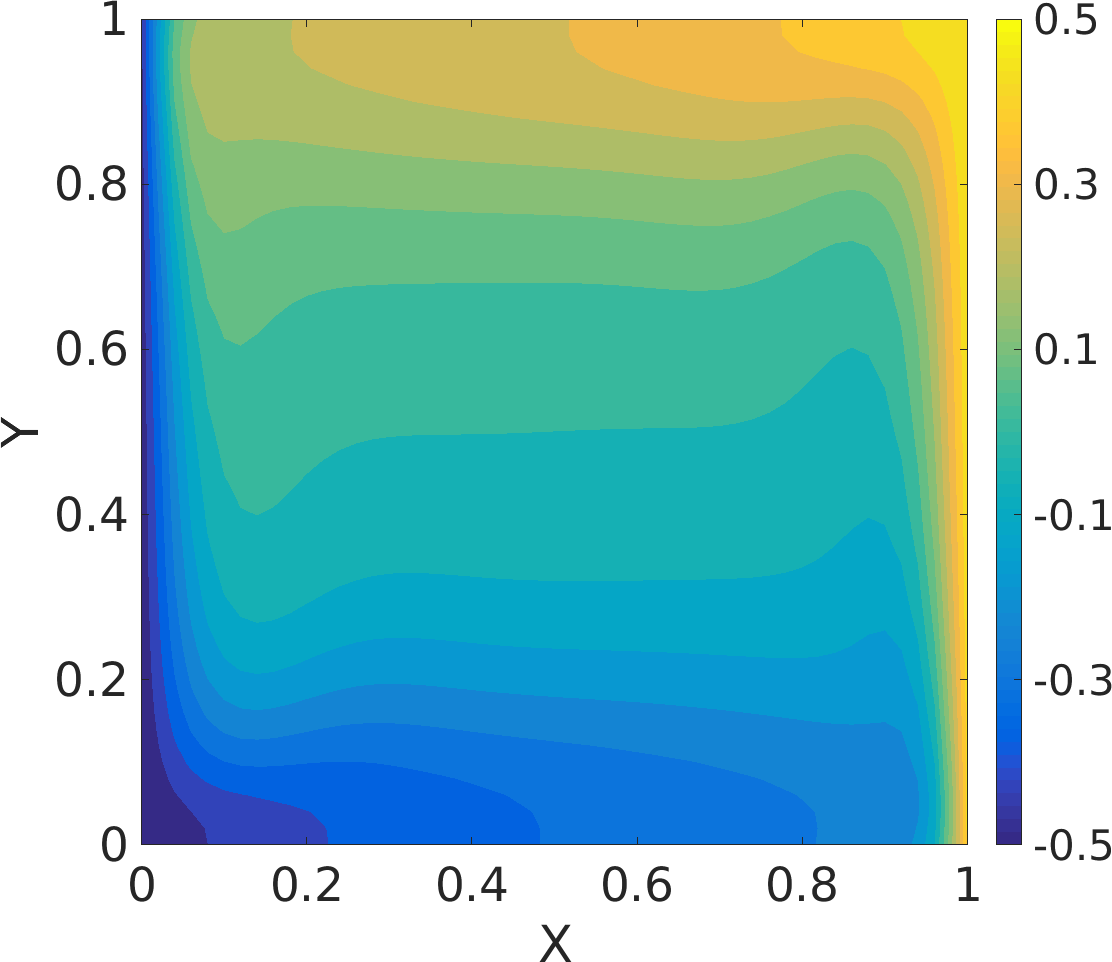}
		\caption{Temperature at $Z=0.5$}
		\label{Fig:Ra_1E6_value_at_input_mean_T}
	\end{subfigure}
	\begin{subfigure}[t]{0.49\textwidth}
		\includegraphics[width=\textwidth]{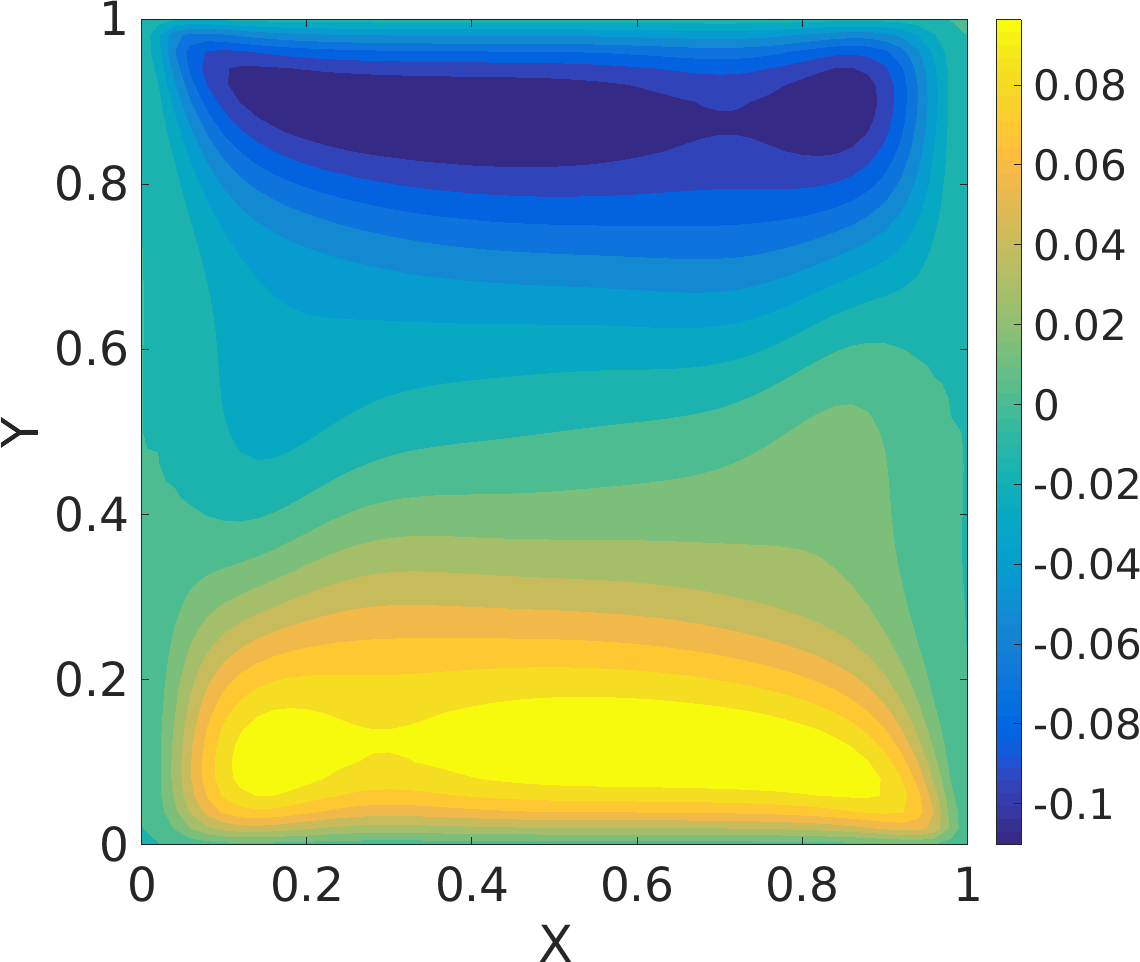}
		\caption{X-Velocity at $Z=0.5$}
		\label{Fig:Ra_1E6_value_at_input_mean_u}
	\end{subfigure}
	\begin{subfigure}[t]{0.49\textwidth}
		\includegraphics[width=\textwidth]{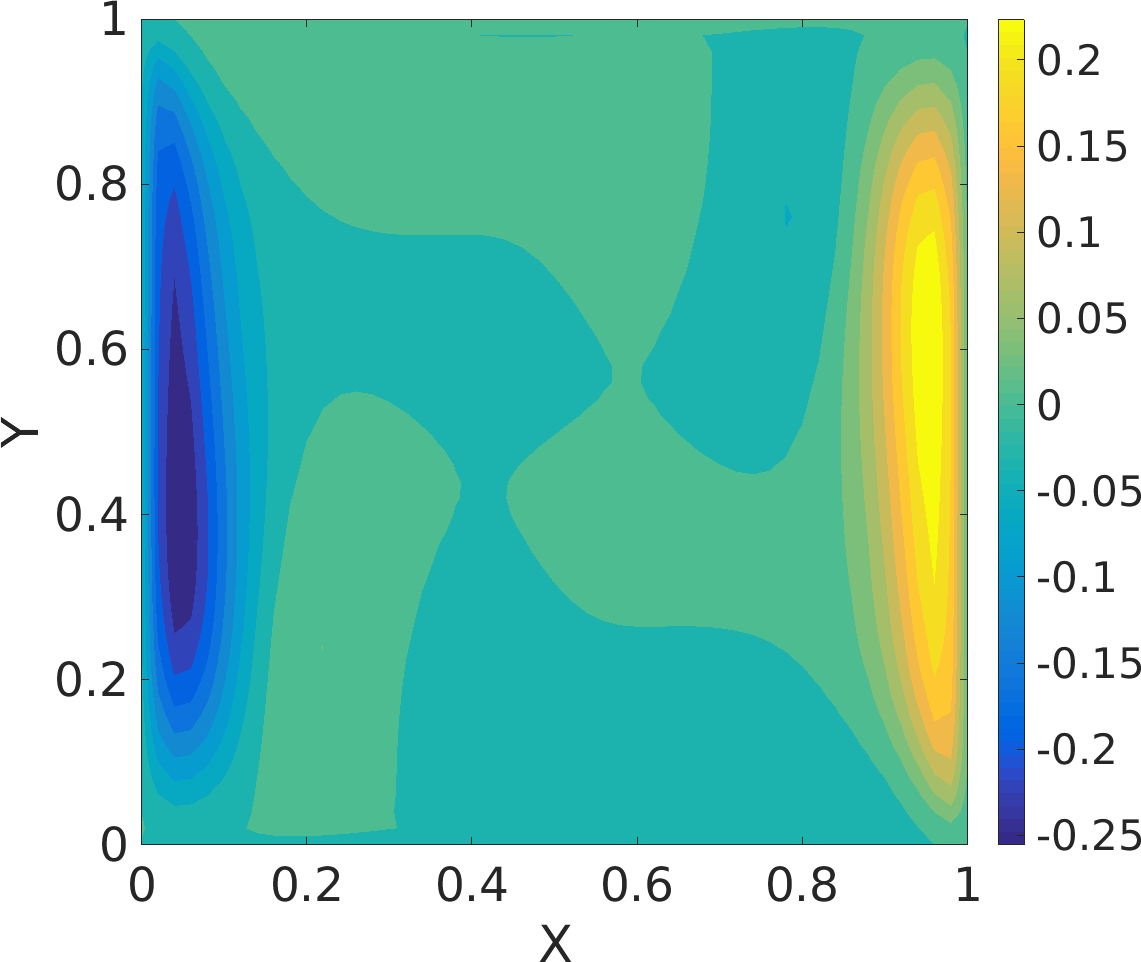}
		\caption{Y-Velocity at $Z=0.5$}
		\label{Fig:Ra_1E6_value_at_input_mean_v}
	\end{subfigure}
	\caption{Deterministic Results for $Ra=10^6$ and $Pr=7.5$}
	\label{Fig:Ra_1E6_value_at_input_mean}
\end{figure}

Uncertainty quantification analyzes the effects of small stochasticity in the input on the output. Since the stochasticity in the input is of the order of a small percentage of its mean, a similar variation is expected in the output. Thus, in this paper, a comparison of the stochastic mean of each output variable is done with the value of that variable at the input mean. The temperature and velocities are non-dimensionalized as discussed in section~\ref{Sec:Governing Equation}. \Cref{Fig:Ra_1E5_value_at_input_mean,Fig:Ra_1E6_value_at_input_mean} plot the Nusselt number at hot wall and temperature, X and Y velocities along the Z midplane for Rayleigh number of $10^5$ and $10^6$, respectively. The stochastic means of the Nusselt number, temperature and velocities are expected to follow trends similar to the \cref{Fig:Ra_1E5_value_at_input_mean,Fig:Ra_1E6_value_at_input_mean}. Thus, in the following sections, contour plots of the difference between the output stochastic mean and the deterministic values are plotted for comparison. The difference gives an estimate of the effect of uncertainty in the input on the outputs compared to the case with deterministic inputs. 
\subsection{Case A}
The non-dimensionalized governing equations (\ref{Eq:Continuity}-\ref{Eq:Energy}) show that the natural convection problem is parametrized by two parameters viz. the Rayleigh number and the Prandtl number. It is assumed that both of them follow a normal distribution with a $2\%$ standard deviation with respect to mean:
\begin{itemize}
	\item $Ra \sim \mathcal{N}(\mu=10^5,\sigma=0.02 \mu)$ or $Ra \sim \mathcal{N}(\mu=10^6,\sigma=0.02 \mu)$
	\item $Pr \sim \mathcal{N}(\mu=7.5,\sigma=0.02 \mu)$
\end{itemize}
Results from two different Rayleigh numbers ($10^5$ and $10^6$) for which the flow is known to remain laminar are presented here \cite{fusegi1991numerical}. High Rayleigh number implies higher buoyancy compared to viscous forces and thus, higher velocities and Nusselt number are observed. The fluid inside the cube is assumed to be water and hence, the Prandtl number is set to 7.5. As mentioned in section \ref{Sec:Polynomial Chaos Expansion}, a tensor product of roots of Hermite polynomial scaled with mean and standard deviation is chosen as samples.
\subsubsection{Convergence of the Stochastic Method}
Hundred uniform Latin hypercube samples are used as test points to verify the convergence of the stochastic collocation method. Two independent estimates of the same output parameter are obtained using polynomial chaos expansion and deterministic simulation. The root mean square of the difference between these two estimates normalized by the maximum value of the parameter is defined as the non-dimensional error estimate. Spatial mean Nusselt number (\cref{Eq:Mean_Nu}) over the hot face is used to estimate the collocation error.
\begin{equation} 
Nu_{mean} = \int Nu(y,z)\Bigr|_{\substack{x=1}} dy dz = \int \frac{\partial T(y,z)}{\partial x}\Bigr|_{\substack{x=1}} dy dz
\label{Eq:Mean_Nu}
\end{equation}
\par First column of \cref{Table:Stochastic Collocation Error Analysis} is the accuracy level of the sample points used for interpolation. Accuracy level $l$ integrates polynomials up to degree $2l-1$ exactly \cite{heiss2008likelihood}. Second column is the number of sample points i.e., the number of deterministic simulations required ($M$ in \cref{Eq:Collocation_vandermonde}). The last two columns list the non-dimensional RMS error in computation of the spatial mean Nusselt number for both the Rayleigh numbers. Although the error increases slightly at the first two levels, eventually the error drops with higher accuracy level. This proves the convergence of the stochastic method. At the highest accuracy level, the error is of order $10^{-4}$ or $10^{-5}$ which shows that the polynomial chaos is reasonably accurate and can be used for further analysis.
\begin{table}[]
	\centering
	\begin{tabular}{|c|c|c|c|}
		\hline
		Accuracy Level & \# Samples & Ra 1E5   & Ra 1E6   \\ \hline
		4              & 16         & --       & 1.09E-03 \\ \hline
		5              & 25         & 3.45E-05 & 1.32E-03 \\ \hline
		6              & 36         & 6.76E-05 & 5.24E-04 \\ \hline
		7              & 49         & 1.79E-05 & 2.50E-04 \\ \hline
		8              & 64         & 1.50E-05 & --       \\ \hline
	\end{tabular}
	\caption{Stochastic Collocation Error Analysis}
	\label{Table:Stochastic Collocation Error Analysis}
\end{table}
\begin{figure}[H]
	\centering
	\begin{subfigure}[t]{0.49\textwidth}
		\includegraphics[width=\textwidth]{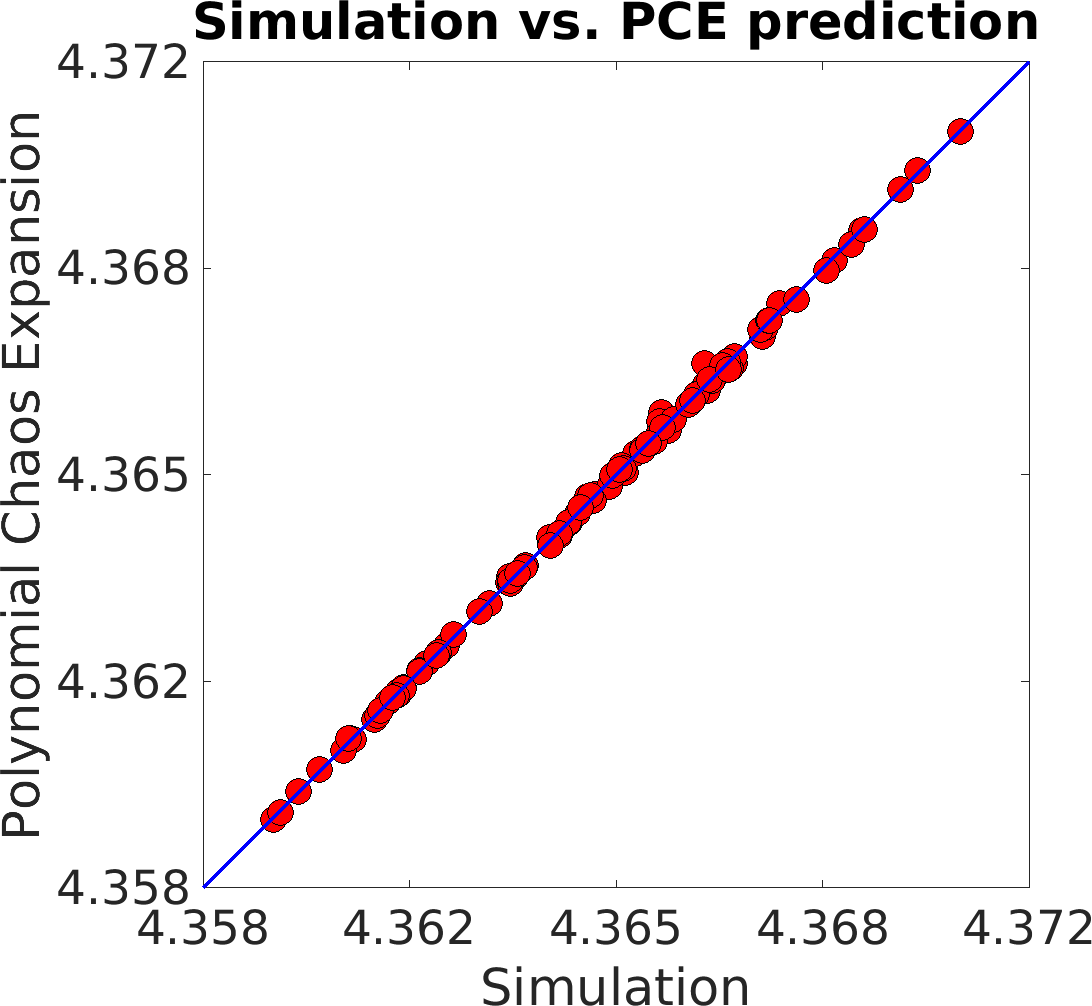}
		\caption{$Ra=10^5$}
	\end{subfigure}
	\begin{subfigure}[t]{0.49\textwidth}
		\includegraphics[width=\textwidth]{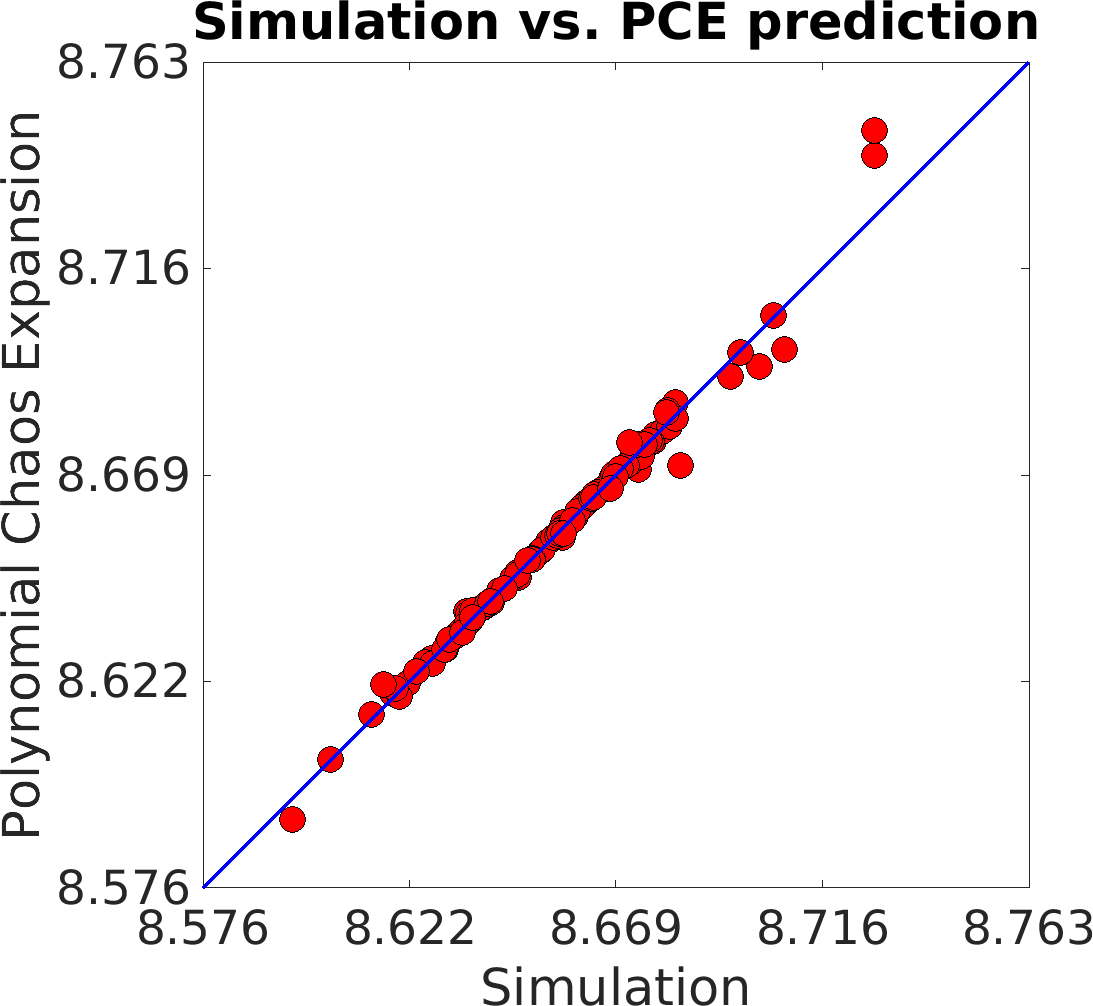}
		\caption{$Ra=10^6$}
	\end{subfigure}
	\caption{Mean Nusselt Number Estimate from Numerical Simulation and Polynomial Chaos Expansion}
	\label{Fig:PCE_vs_Det_Sim}
\end{figure}
\par For visual inspection, spatial mean Nusselt number estimates from the numerical simulation and the polynomial chaos expansion are plotted together in \cref{Fig:PCE_vs_Det_Sim} for the hundred test points. Ideally, all the points should lie on the $Y=X$ line but due to the stochastic interpolation error, some points are off the line. Since most of the points follow the expected trend of the $Y=X$ line, it can be concluded that the polynomial chaos is accurate.
\subsubsection{Nusselt Number}
\begin{figure}[H]
	\centering
	\begin{subfigure}[t]{0.49\textwidth}
		\includegraphics[width=\textwidth]{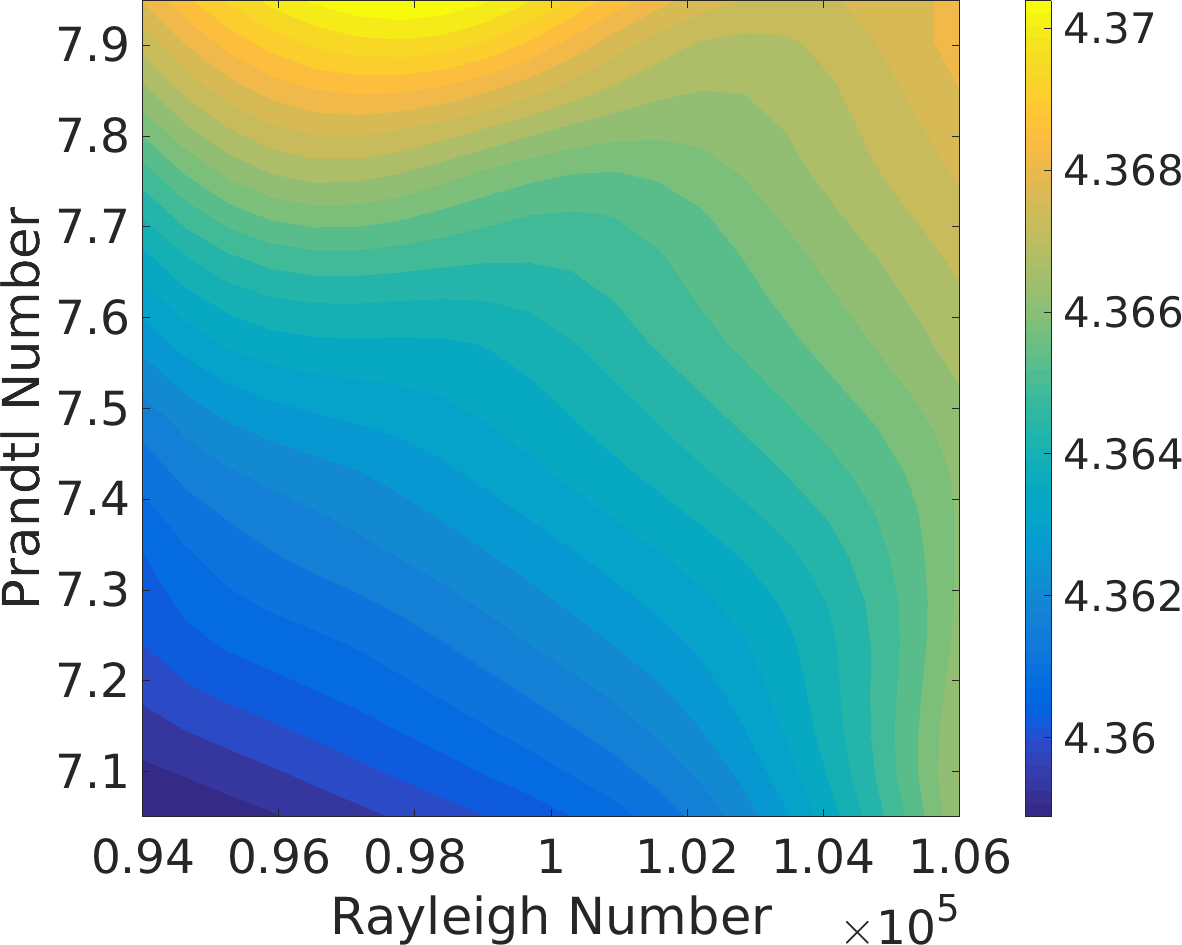}
		\caption{$Ra=10^5$}
		\label{Fig:Nu_Resp_Surf_1E5}
	\end{subfigure}
	\begin{subfigure}[t]{0.49\textwidth}
		\includegraphics[width=\textwidth]{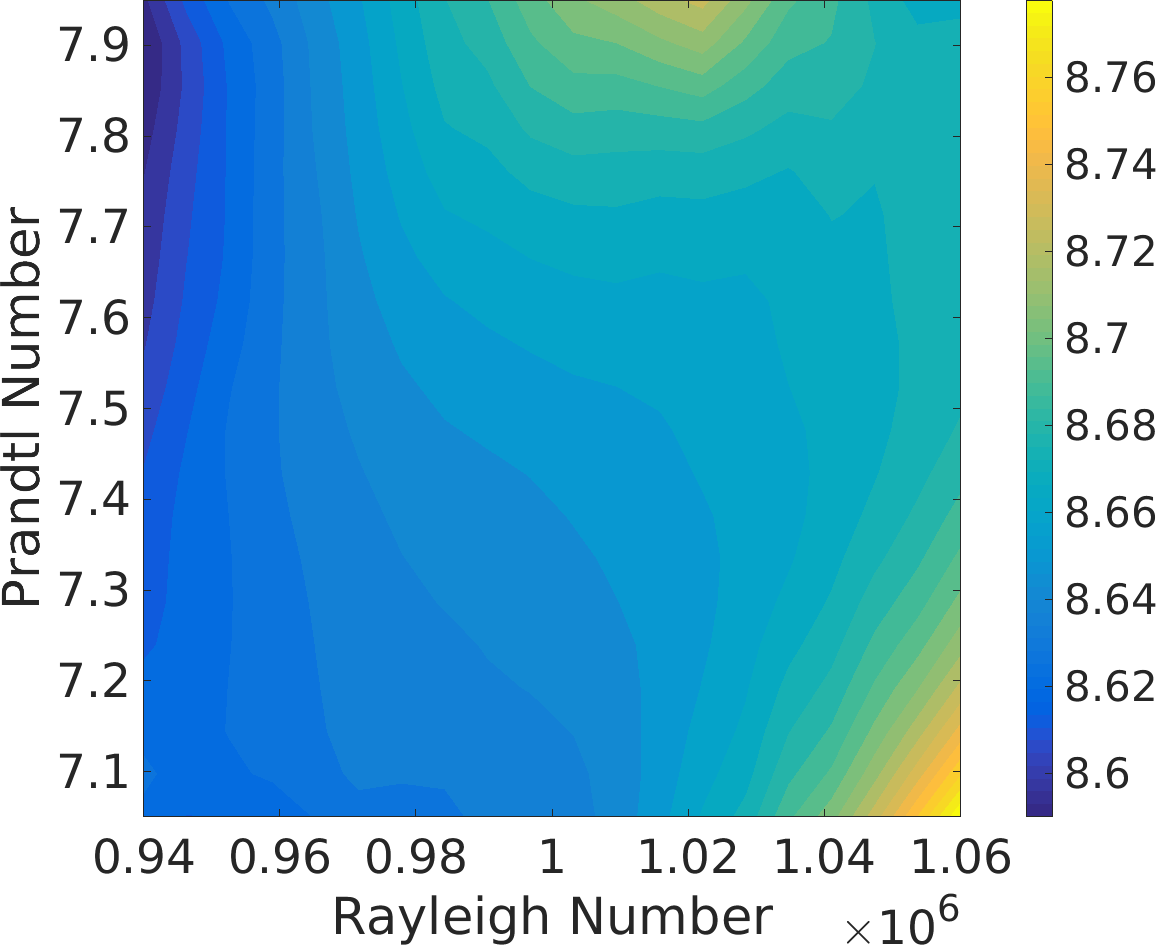}
		\caption{$Ra=10^6$}
		\label{Fig:Nu_Resp_Surf_1E6}
	\end{subfigure}
	\caption{Spatial Mean Nusselt Number Response Surface}
	\label{Fig:Nu_Resp_Surf}
\end{figure}

Response surface gives a visual representation of the variation of an output parameter with input stochastic parameters. For a two dimensional stochastic problem, the response surface can be plotted as a contour. \Cref{Fig:Nu_Resp_Surf} plots the response surfaces of Nusselt number averaged over the hot wall (\cref{Eq:Mean_Nu}) for both the Rayleigh numbers. In each plot, X and Y axes denote Rayleigh and Prandtl numbers, respectively which are the stochastic input parameters. Since both the input variables are assumed to follow normal distribution, they are plotted in the range $( \mu-3 \sigma,\mu+3 \sigma )$. For example, Rayleigh and Prandtl numbers are plotted in the range $(10^5-3\times2000,10^5+3\times2000)$ and $(7.5-3\times0.15,7.5+3\times0.15)$, respectively in \cref{Fig:Nu_Resp_Surf_1E5}. The contour lines represent the value of mean Nusselt number. The slope of the contour line can be used to estimate the local sensitivity of the output with respect to a particular input. For example, the contour lines are nearly vertical in the left region of \cref{Fig:Nu_Resp_Surf_1E6} which implies that the local sensitivity of the mean Nusselt number in the left region is high towards the input plotted on the X axis i.e., the Rayleigh number.
\begin{figure}[H] 
	\centering
	\begin{subfigure}[t]{0.32\textwidth}
		\includegraphics[width=\textwidth]{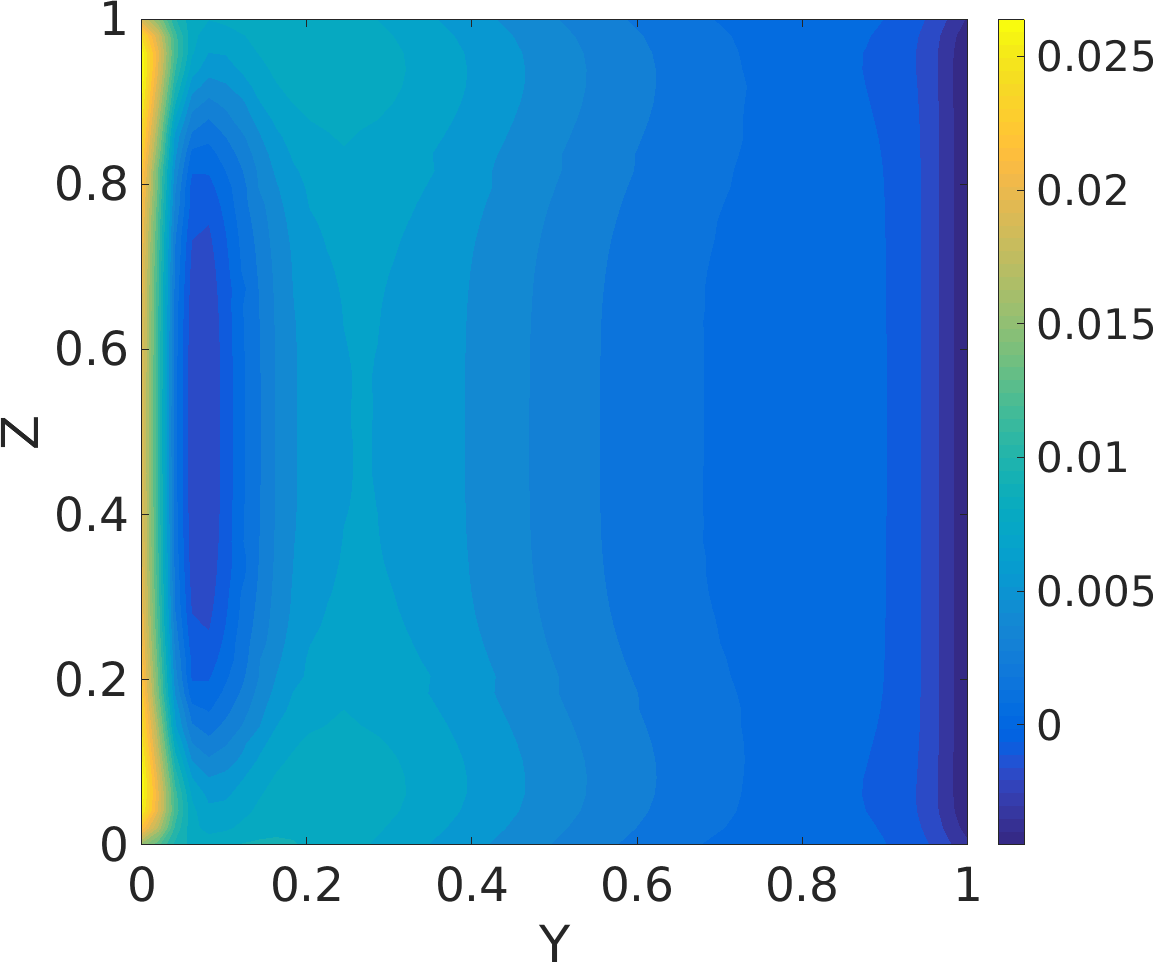}
		\caption{Diff. of Stochastic Mean and Deterministic Value}
		\label{Fig:Nu_1E5_a}
	\end{subfigure}
	\begin{subfigure}[t]{0.32\textwidth}
		\includegraphics[width=\textwidth]{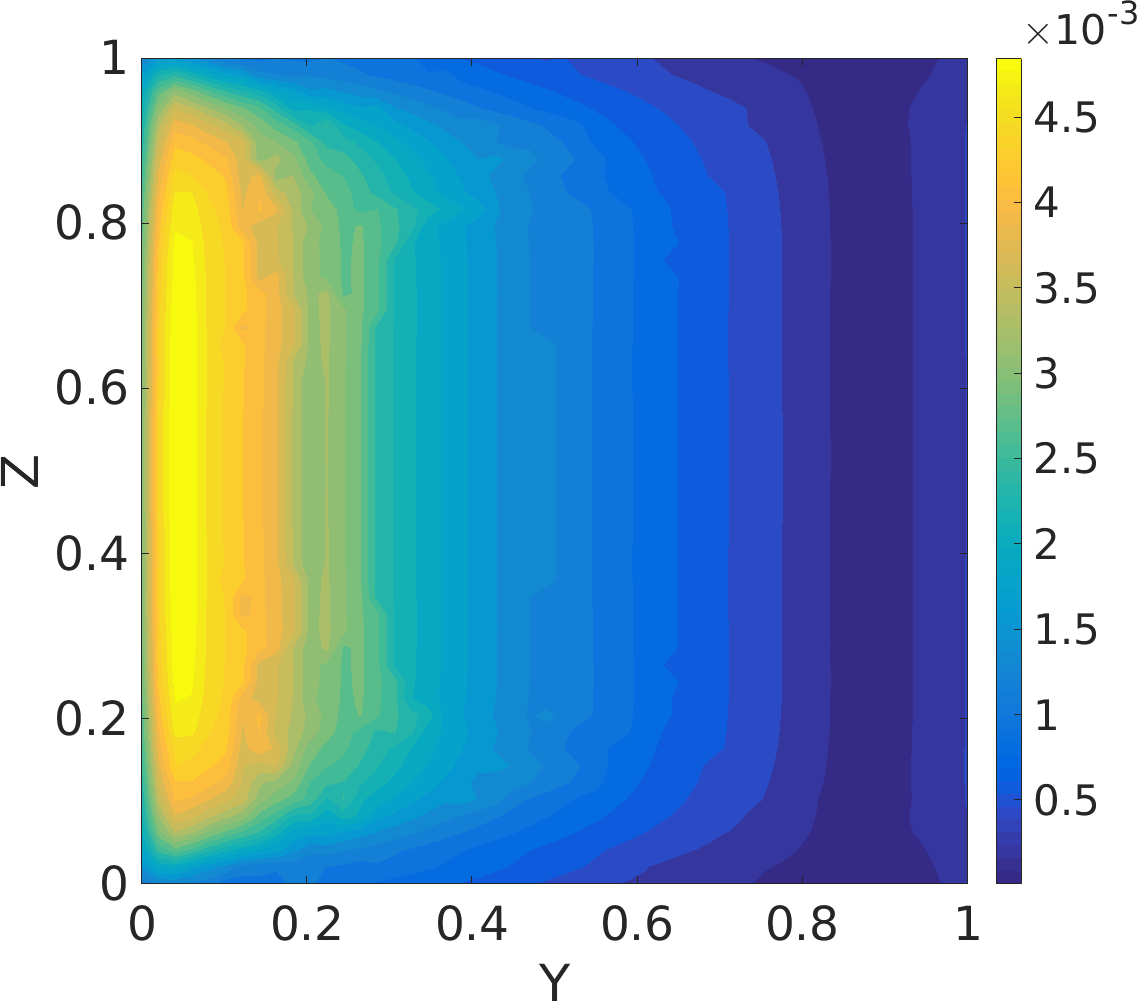}
		\caption{Stochastic Std. Dev.}
		\label{Fig:Nu_1E5_b}
	\end{subfigure}
	\begin{subfigure}[t]{0.32\textwidth}
		\includegraphics[width=\textwidth]{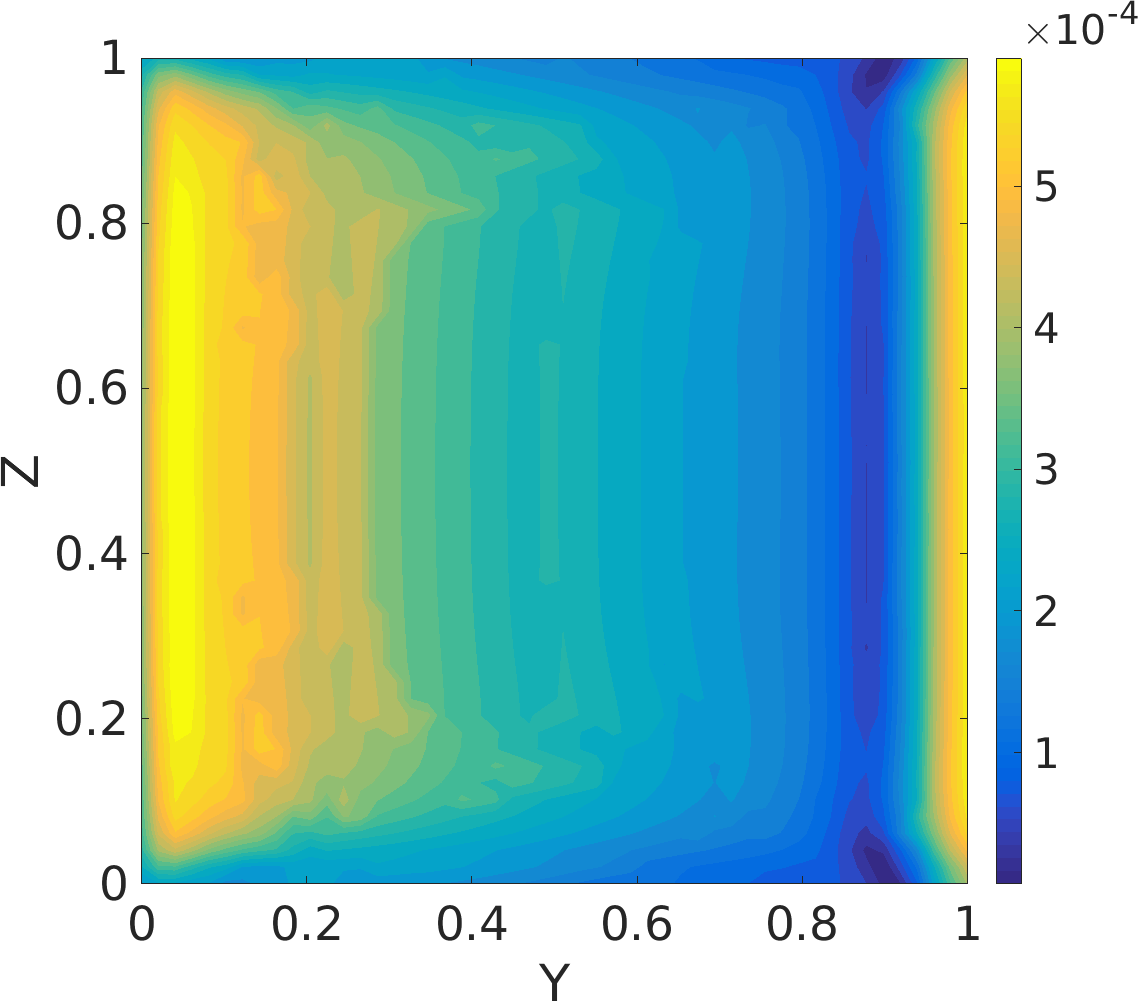}
		\caption{Stochastic Std. Dev. by Stochastic Mean}
		\label{Fig:Nu_1E5_c}
	\end{subfigure}
	\caption{Local Nusselt Number at Hot Wall for $Ra=10^5$}
	\label{Fig:Nu_1E5}
\end{figure}
\begin{figure}[h]
	\centering
	\begin{subfigure}[t]{0.32\textwidth}
		\includegraphics[width=\textwidth]{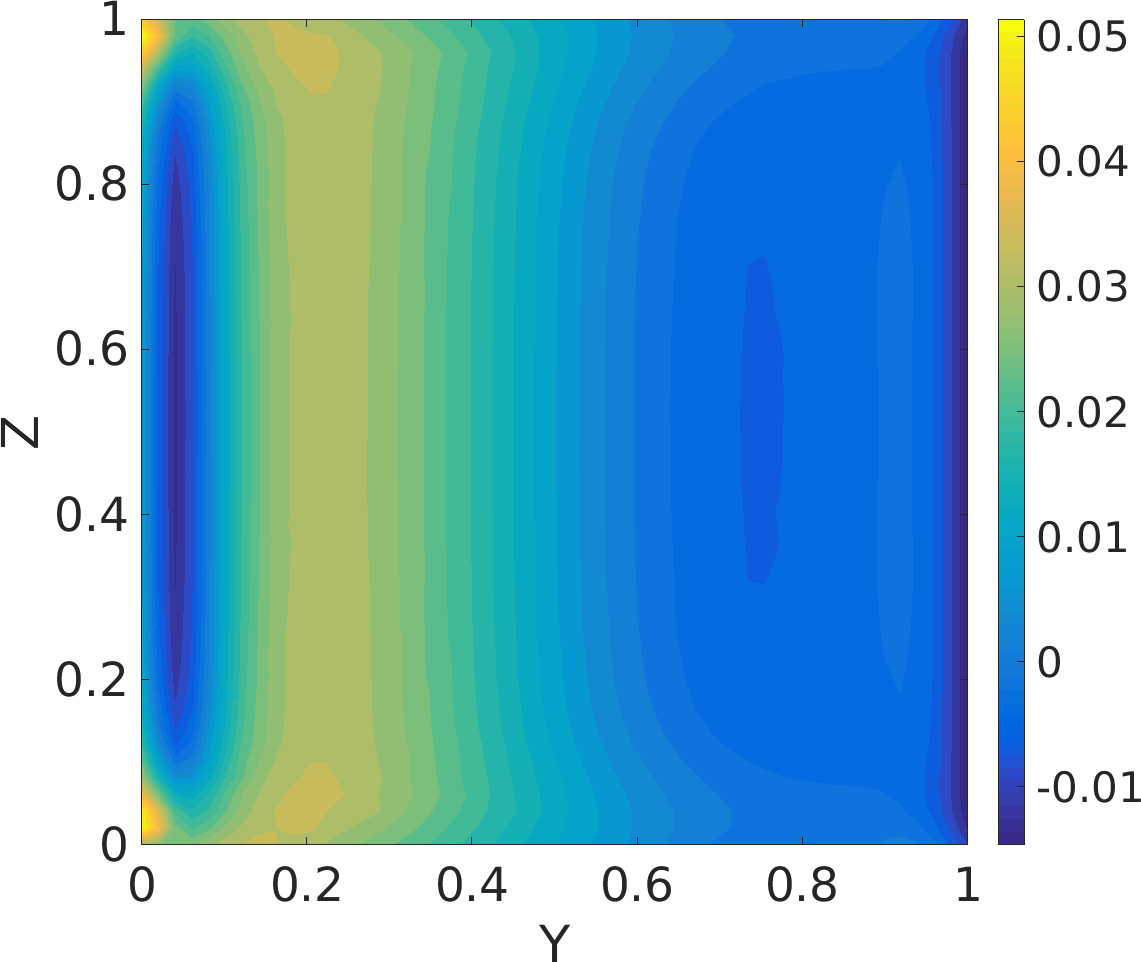}
		\caption{Diff. of Stochastic Mean and Deterministic Value}
		\label{Fig:Nu_1E6_a}
	\end{subfigure}
	\begin{subfigure}[t]{0.32\textwidth}
		\includegraphics[width=\textwidth]{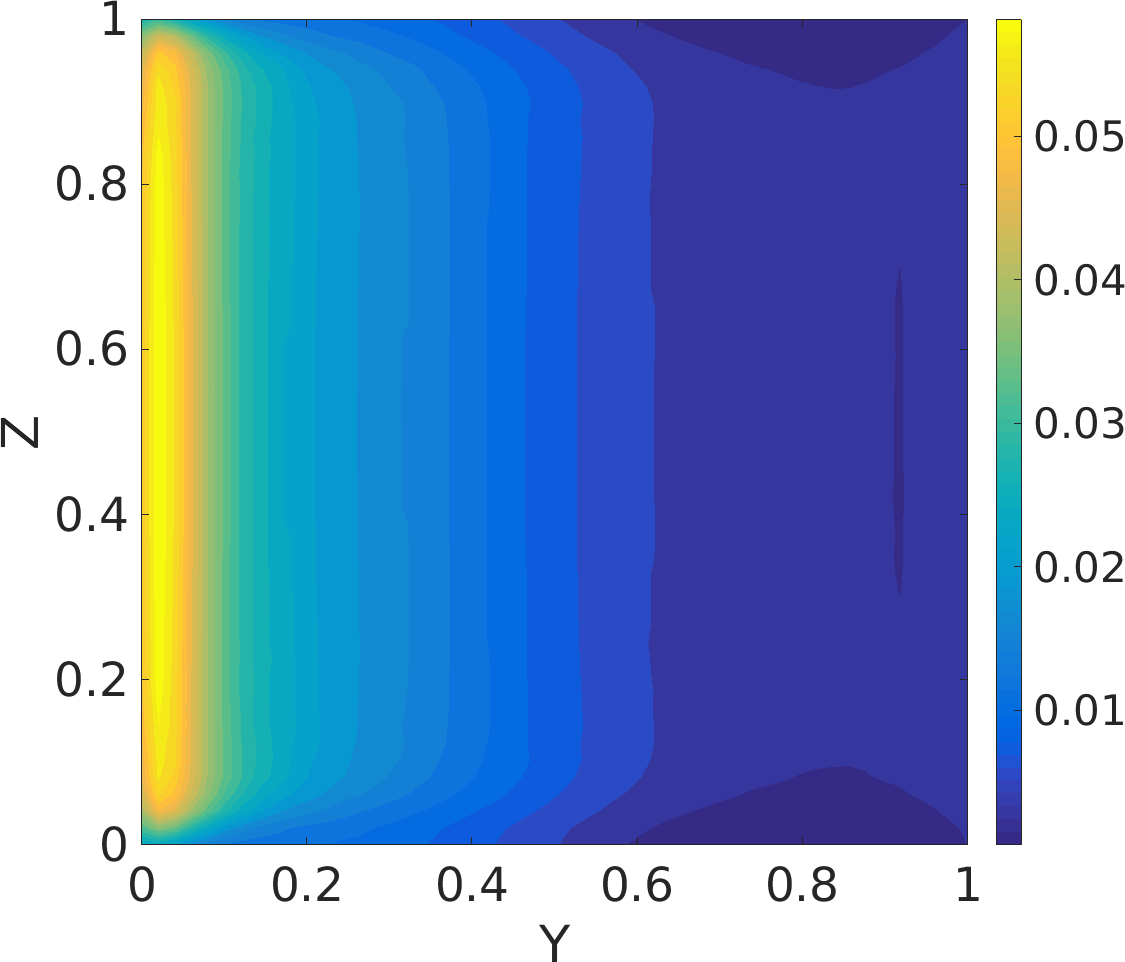}
		\caption{Stochastic Std. Dev.}
		\label{Fig:Nu_1E6_b}
	\end{subfigure}
	\begin{subfigure}[t]{0.32\textwidth}
		\includegraphics[width=\textwidth]{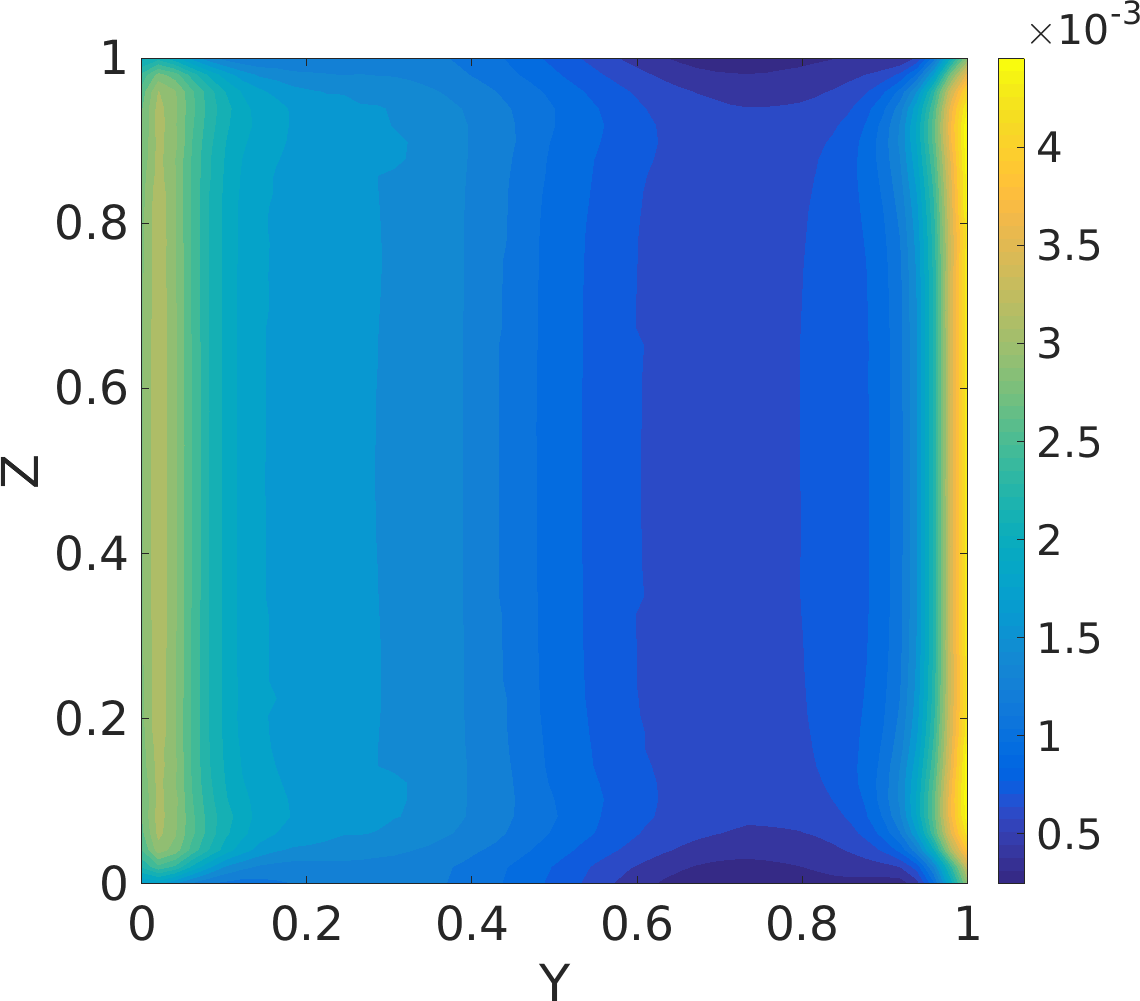}
		\caption{Stochastic Std. Dev. by Stochastic Mean}
		\label{Fig:Nu_1E6_c}
	\end{subfigure}
	\caption{Local Nusselt Number at Hot Wall for $Ra=10^6$}
	\label{Fig:Nu_1E6}
\end{figure}
The local Nusselt number on the wall ($Nu(y,z)$) varies due to the stochasticity in the Rayleigh and Prandtl numbers. \Cref{Fig:Nu_1E5,Fig:Nu_1E6} plot its statistics for both the Rayleigh numbers. The stochastic mean plots look visually similar to the plots at mean Rayleigh and Prandtl numbers (\cref{Fig:Ra_1E5_value_at_input_mean_Nu,Fig:Ra_1E6_value_at_input_mean_Nu}). Thus, \cref{Fig:Nu_1E5_a,Fig:Nu_1E6_a} plot contours of difference between the stochastic mean and deterministic value of the Nusselt number at hot wall (\cref{Fig:Ra_1E5_value_at_input_mean_Nu,Fig:Ra_1E6_value_at_input_mean_Nu}). \Cref{Fig:Nu_1E5_b,Fig:Nu_1E6_b} plot the standard deviation due to the stochastic input variation. Higher standard deviation is observed in the region of higher mean. Thus, \cref{Fig:Nu_1E5_c,Fig:Nu_1E6_c} are plotted to annihilate the effect of mean. The local ratio of standard deviation to mean shows higher values on the left and right sides near $Y=0$ and $Y=1$ as gravity is acting in $Y$ direction. It can be seen that the difference between the stochastic mean and the deterministic values (\cref{Fig:Nu_1E5_a}) is one order of magnitude higher than the stochastic standard deviation (\cref{Fig:Nu_1E5_b}) for Rayleigh number of $10^5$. This implies that the input stochasticity shifts the deterministic mean of the output more than its standard deviation. On the other hand, for Rayleigh number of $10^6$, both the difference (\cref{Fig:Nu_1E6_a}) and stochastic standard deviation (\cref{Fig:Nu_1E6_b}) are of similar orders of magnitude. Thus, the shifting of mean and the standard deviation on the shift are of similar orders of magnitude for higher Rayleigh number. The difference between the stochastic mean and the deterministic value is of the order of $0.3\%$ of the deterministic value for both the Rayleigh numbers.
%\par From figs.~\ref{Fig:Nu_1E5_c} and \ref{Fig:Nu_1E6_c}, it can be seen that the standard deviation is four and three orders of magnitude smaller compared to the mean for Rayleigh number $10^5$ and $10^6$, respectively. On the other hand, the difference between the stochastic means (figs.~\ref{Fig:Nu_1E5_a} and \ref{Fig:Nu_1E6_a}) is three orders of magnitude smaller than the deterministic values (figs.~\ref{Fig:Ra_1E5_value_at_input_mean_Nu} and \ref{Fig:Ra_1E6_value_at_input_mean_Nu}) for both the Rayleigh numbers. 
\subsubsection{Velocity and Temperature}
%\begin{figure}[H] 
%	\centering
%	\begin{subfigure}[t]{0.49\textwidth}
%		\includegraphics[width=\textwidth]{Figures/Ra_1E5_UQ/Temp_difference.png}
%		\caption{Difference between Stochastic Mean and Deterministic Value}
%		\label{Fig:Case_A_Temp_1E5_a}
%	\end{subfigure}
%	\begin{subfigure}[t]{0.49\textwidth}
%		\includegraphics[width=\textwidth]{Figures/Ra_1E5_UQ/Temp_std.png}
%		\caption{Stochastic Standard Deviation}
%		\label{Fig:Case_A_Temp_1E5_b}
%	\end{subfigure}
%	\caption{Temperature at $Z=0.5$ Mid-plane for $Ra=10^5$}
%	\label{Fig:Case_A_Temp_1E5}
%\end{figure}
\begin{figure}[H] 
	\centering
	\begin{subfigure}[t]{0.32\textwidth}
		\includegraphics[width=\textwidth]{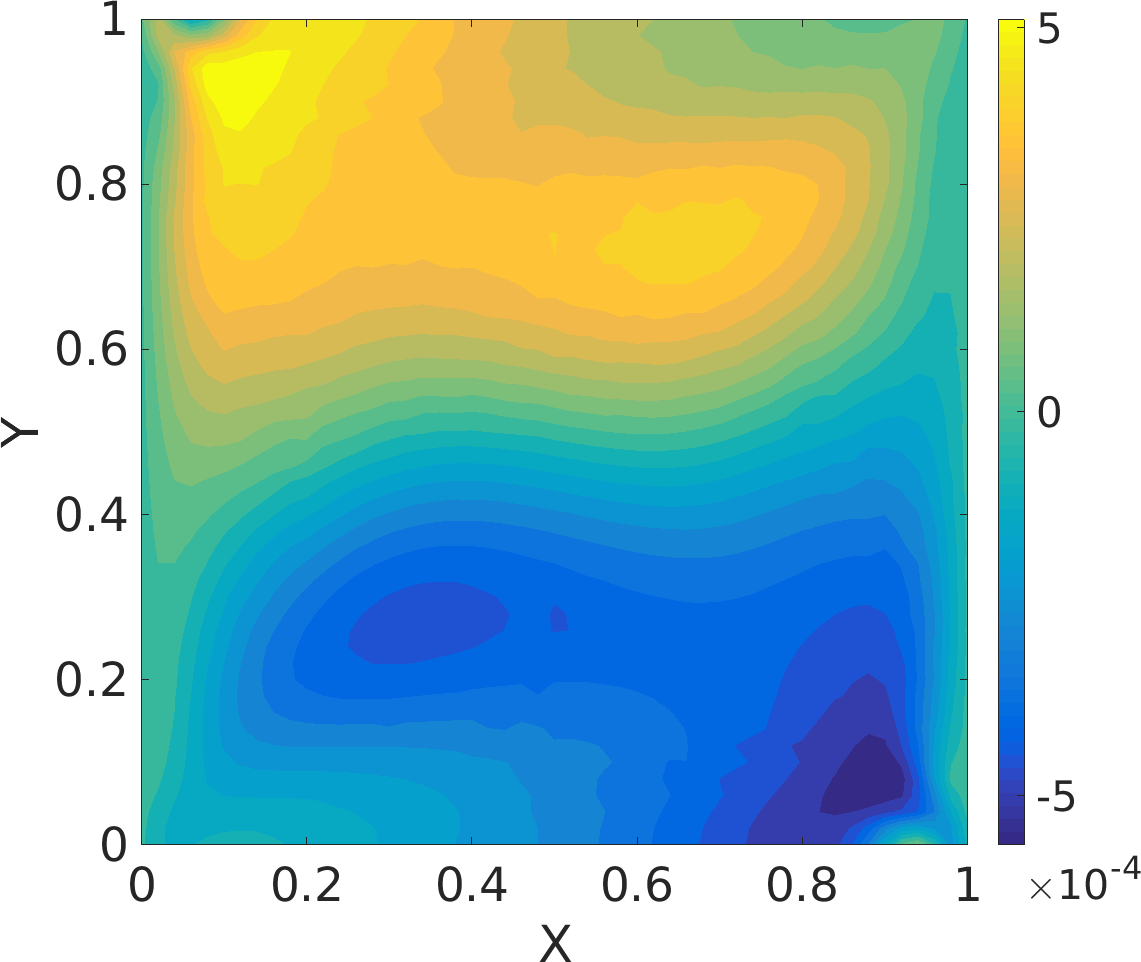}
		\caption{Temperature}
	\end{subfigure}
	\begin{subfigure}[t]{0.32\textwidth}
		\includegraphics[width=\textwidth]{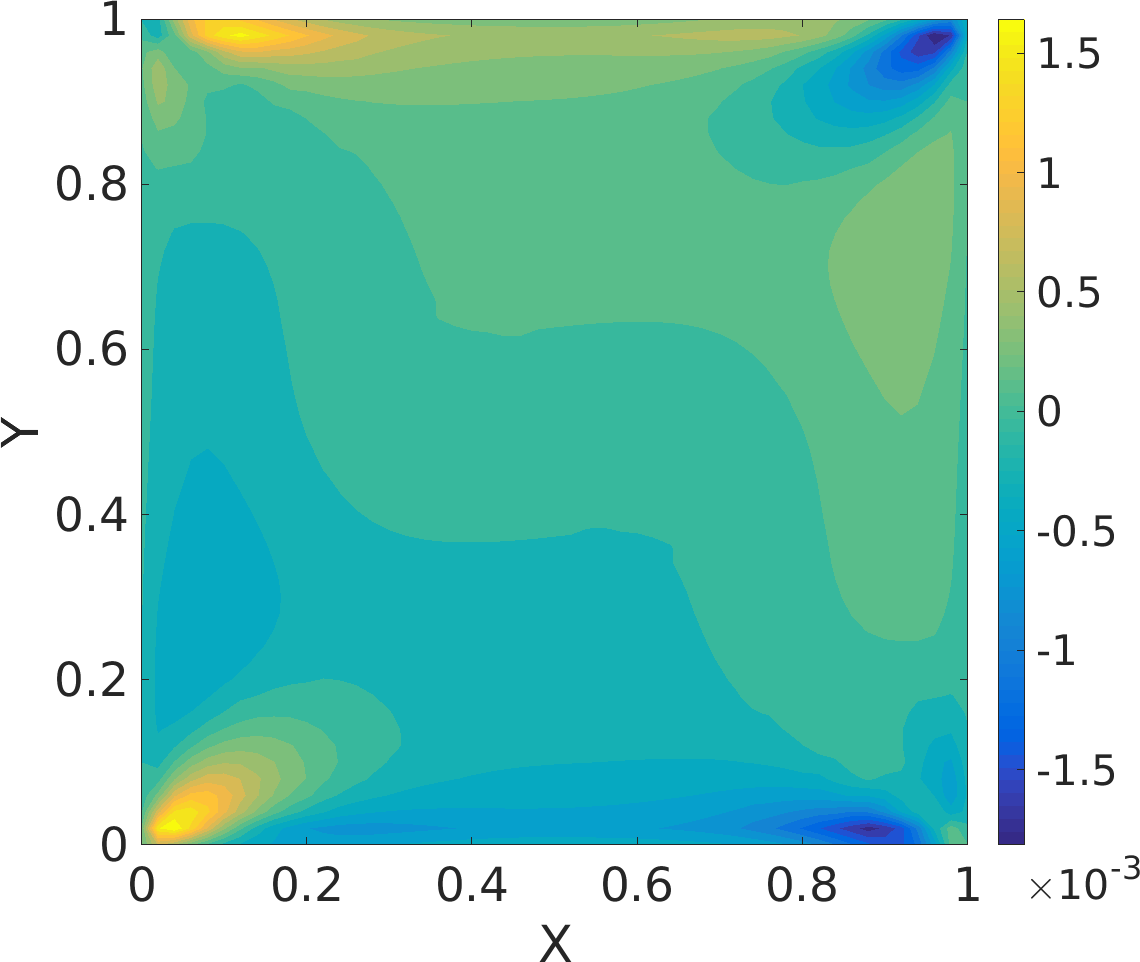}
		\caption{X Velocity}
	\end{subfigure}
	\begin{subfigure}[t]{0.32\textwidth}
		\includegraphics[width=\textwidth]{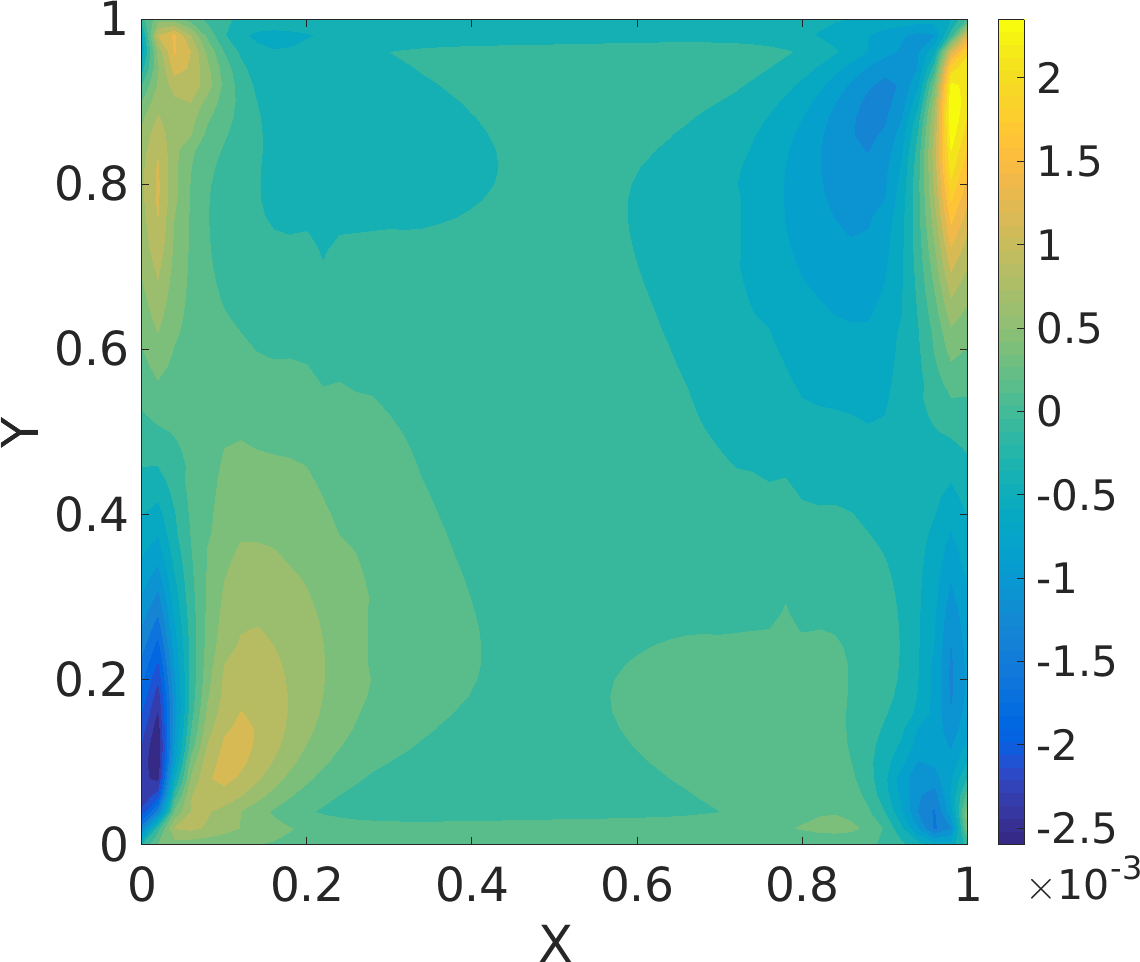}
		\caption{Y Velocity}
	\end{subfigure}
	\caption{Difference between Stochastic Mean and Deterministic Value at $Z=0.5$ Mid-plane for $Ra=10^5$}
	\label{Fig:Case_A_1E5_Diff}
\end{figure}
\begin{figure}[H] 
	\begin{subfigure}[t]{0.32\textwidth}
		\includegraphics[width=\textwidth]{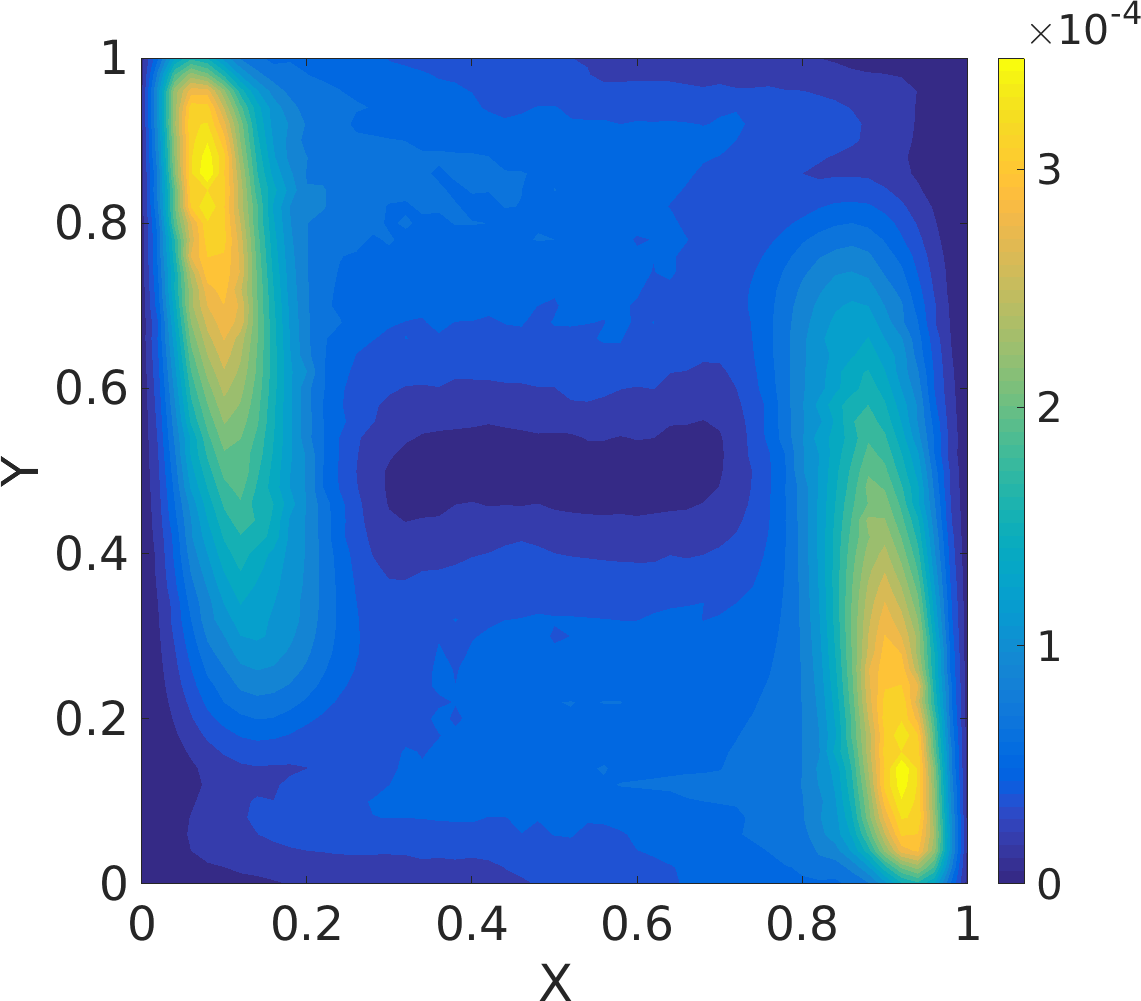}
		\caption{Temperature}
	\end{subfigure}
	\begin{subfigure}[t]{0.32\textwidth}
		\includegraphics[width=\textwidth]{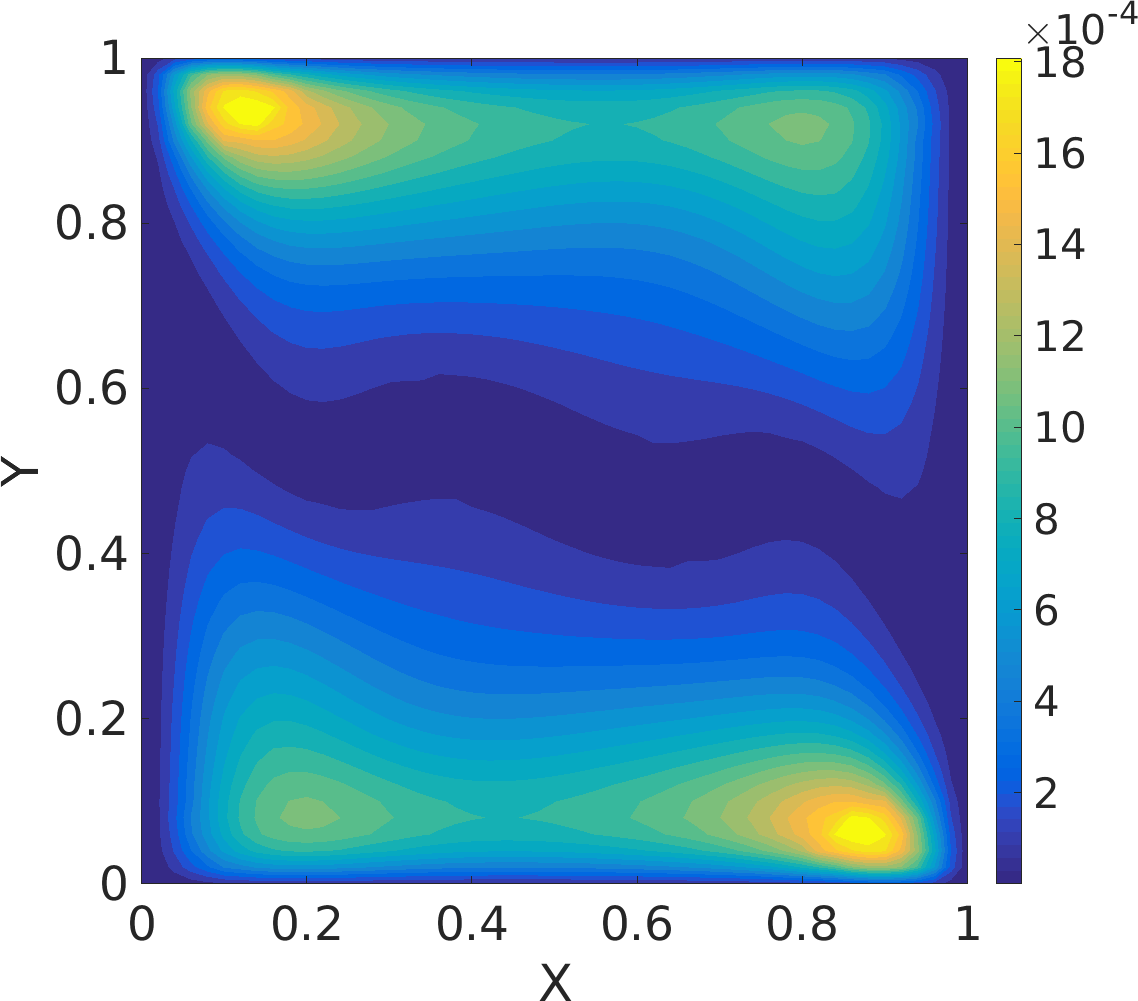}
		\caption{X Velocity}
	\end{subfigure}
	\begin{subfigure}[t]{0.32\textwidth}
		\includegraphics[width=\textwidth]{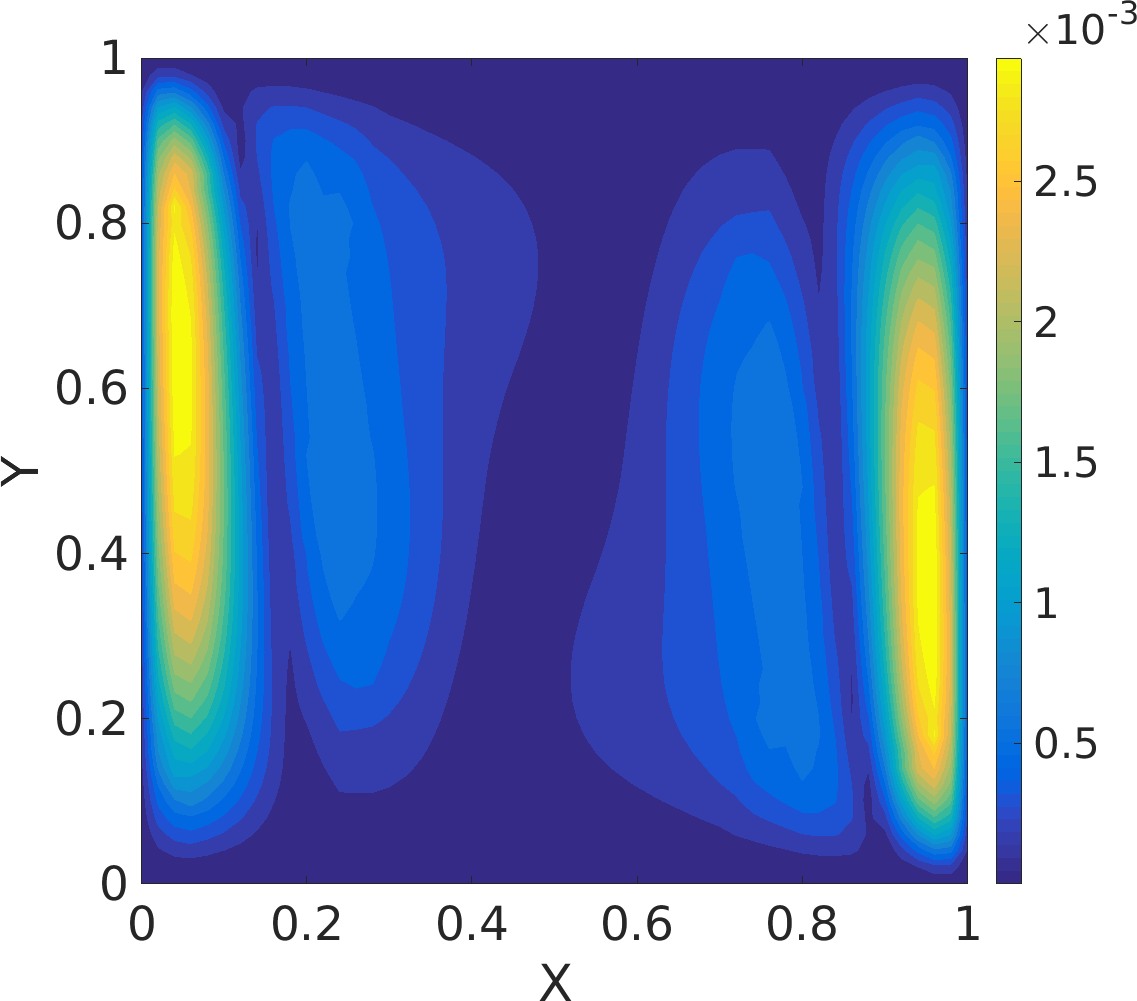}
		\caption{Y Velocity}
	\end{subfigure}
	\caption{Stochastic Standard Deviation at $Z=0.5$ Mid-plane for $Ra=10^5$}
	\label{Fig:Case_A_1E5_Std}
\end{figure}
\begin{figure}[H] 
	\centering
	\begin{subfigure}[t]{0.32\textwidth}
		\includegraphics[width=\textwidth]{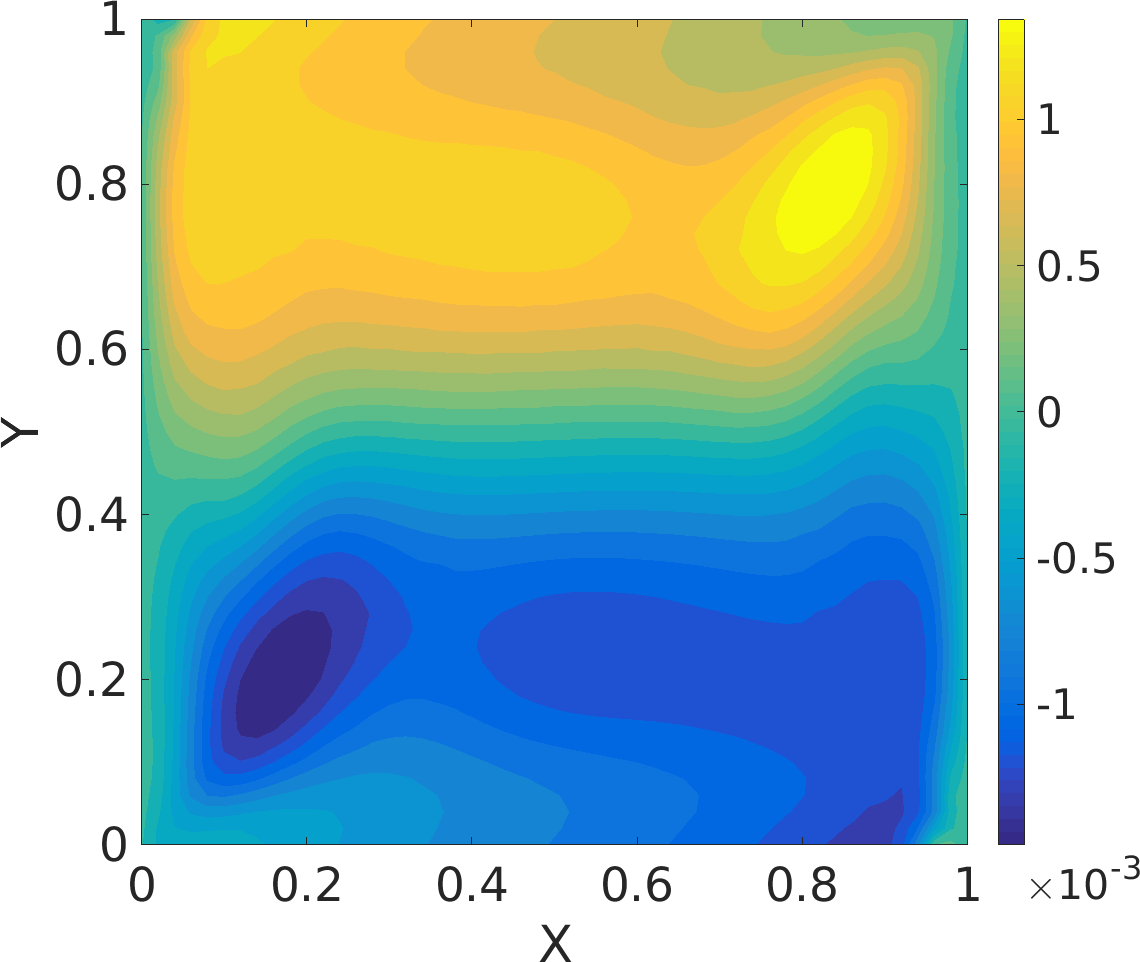}
		\caption{Temperature}
	\end{subfigure}
	\begin{subfigure}[t]{0.32\textwidth}
		\includegraphics[width=\textwidth]{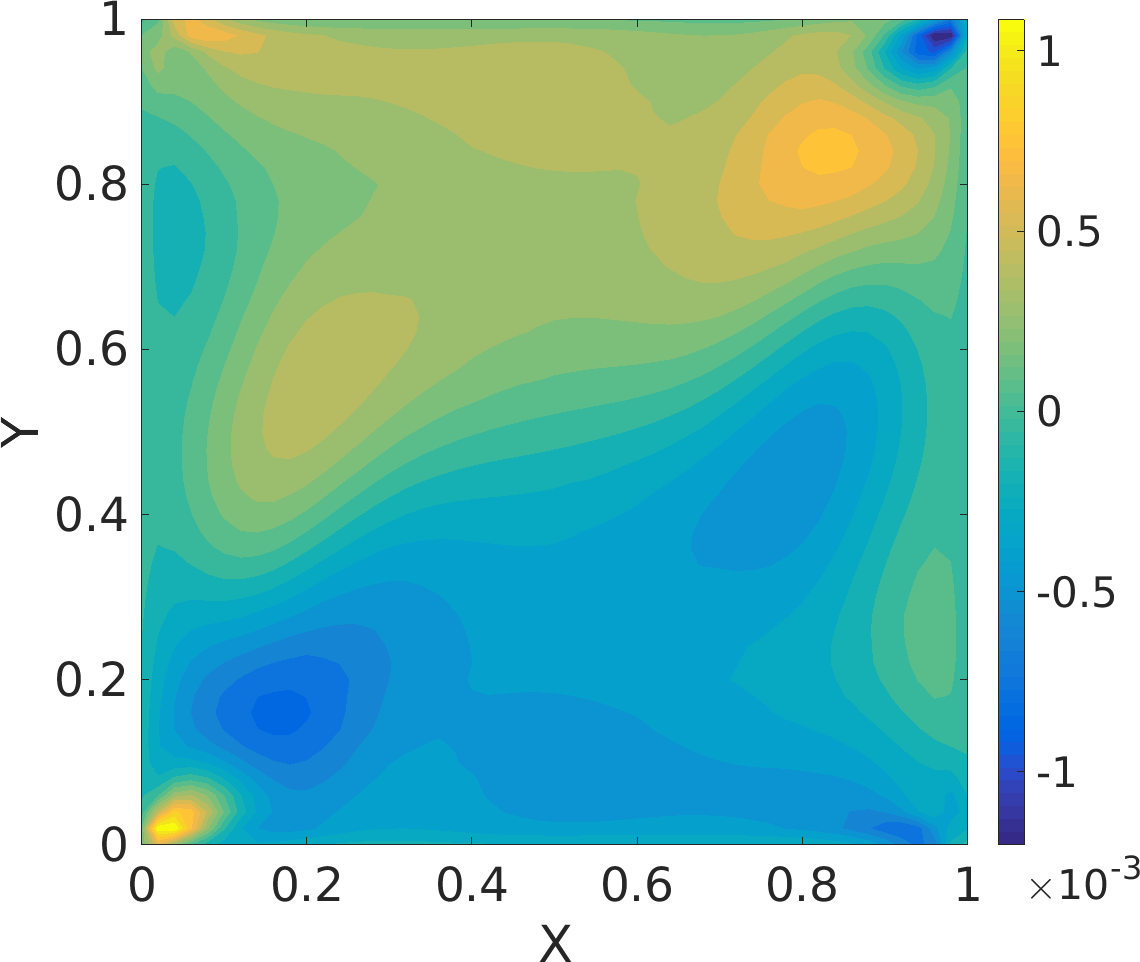}
		\caption{X Velocity}
	\end{subfigure}
	\begin{subfigure}[t]{0.32\textwidth}
		\includegraphics[width=\textwidth]{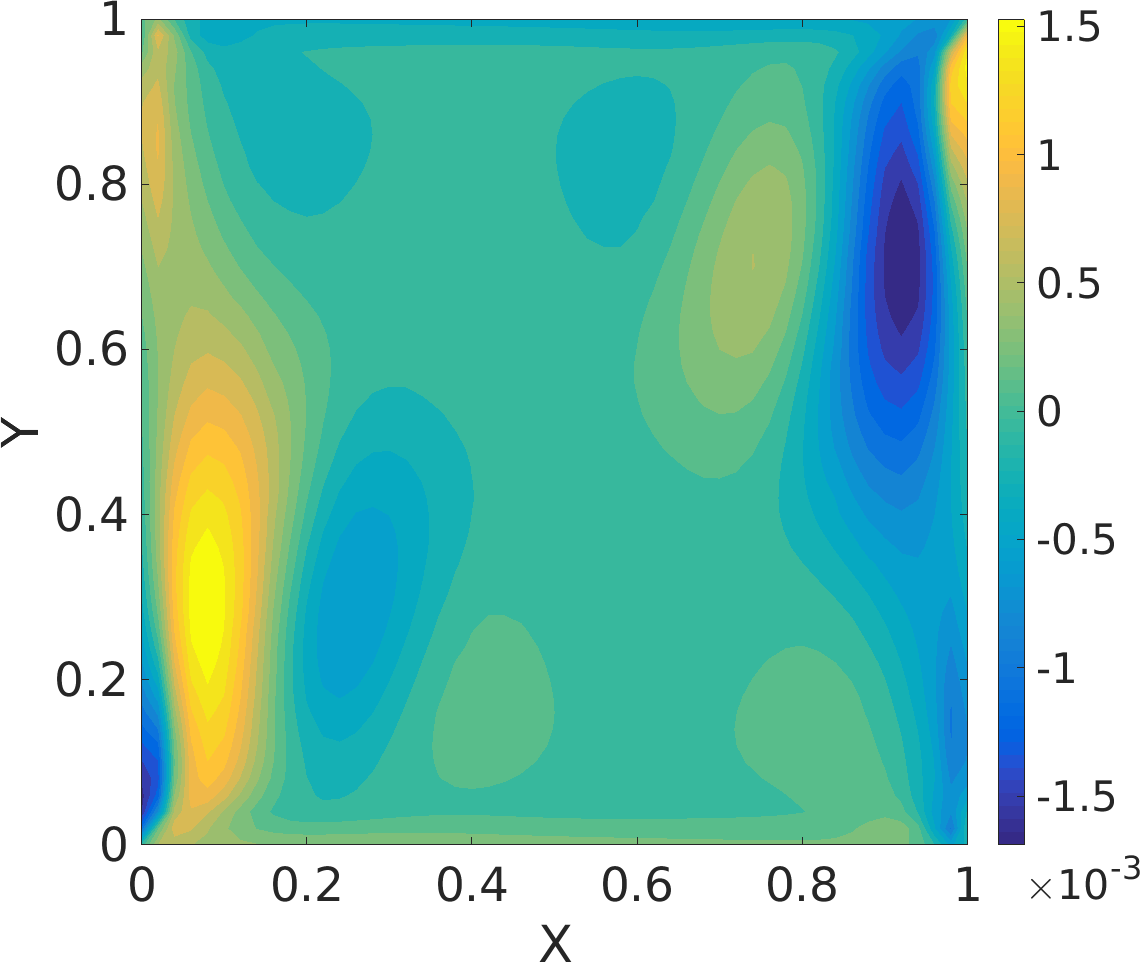}
		\caption{Y Velocity}
	\end{subfigure}
	\caption{Difference between Stochastic Mean and Deterministic Value at $Z=0.5$ Mid-plane for $Ra=10^6$}
	\label{Fig:Case_A_1E6_Diff}
\end{figure}
\begin{figure}[H] 
	\begin{subfigure}[t]{0.32\textwidth}
		\includegraphics[width=\textwidth]{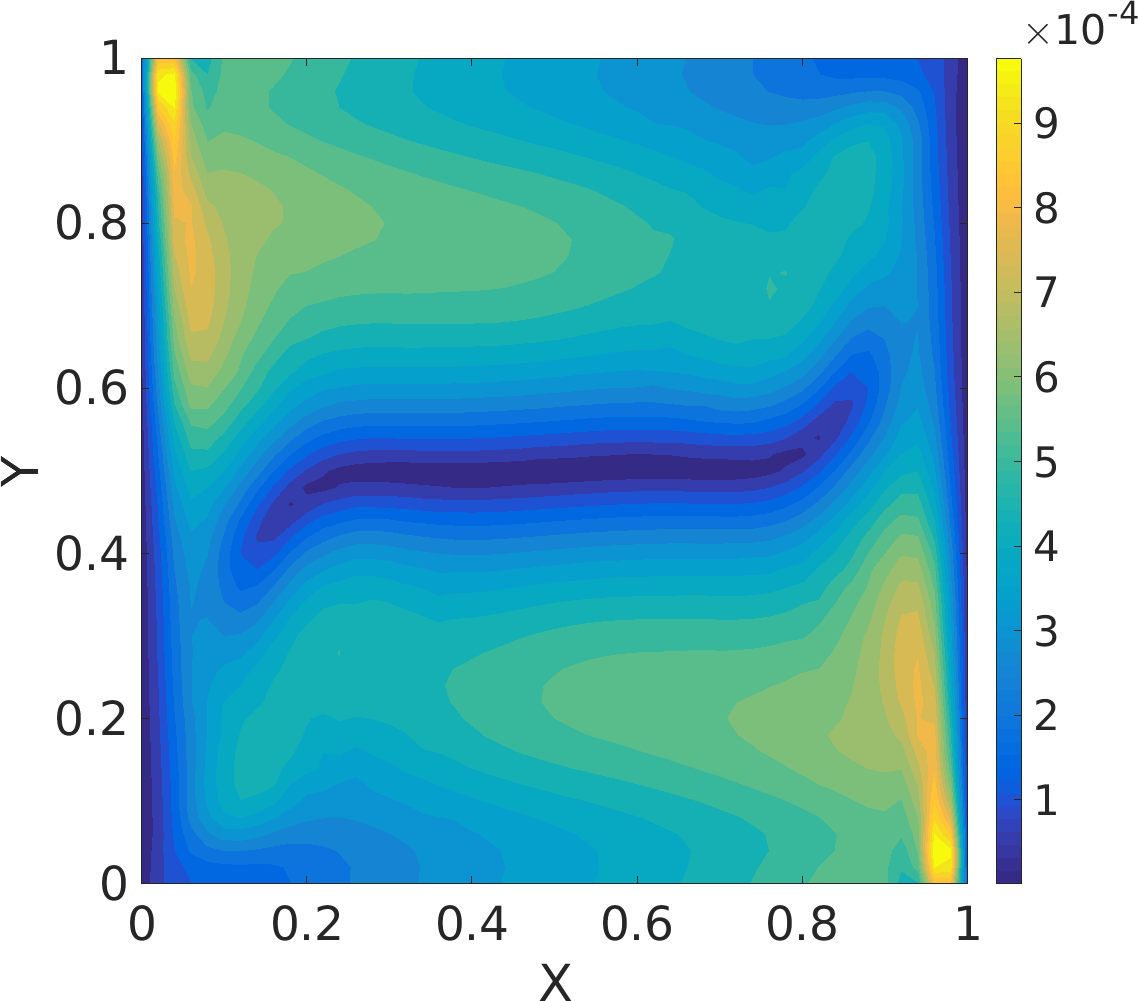}
		\caption{Temperature}
	\end{subfigure}
	\begin{subfigure}[t]{0.32\textwidth}
		\includegraphics[width=\textwidth]{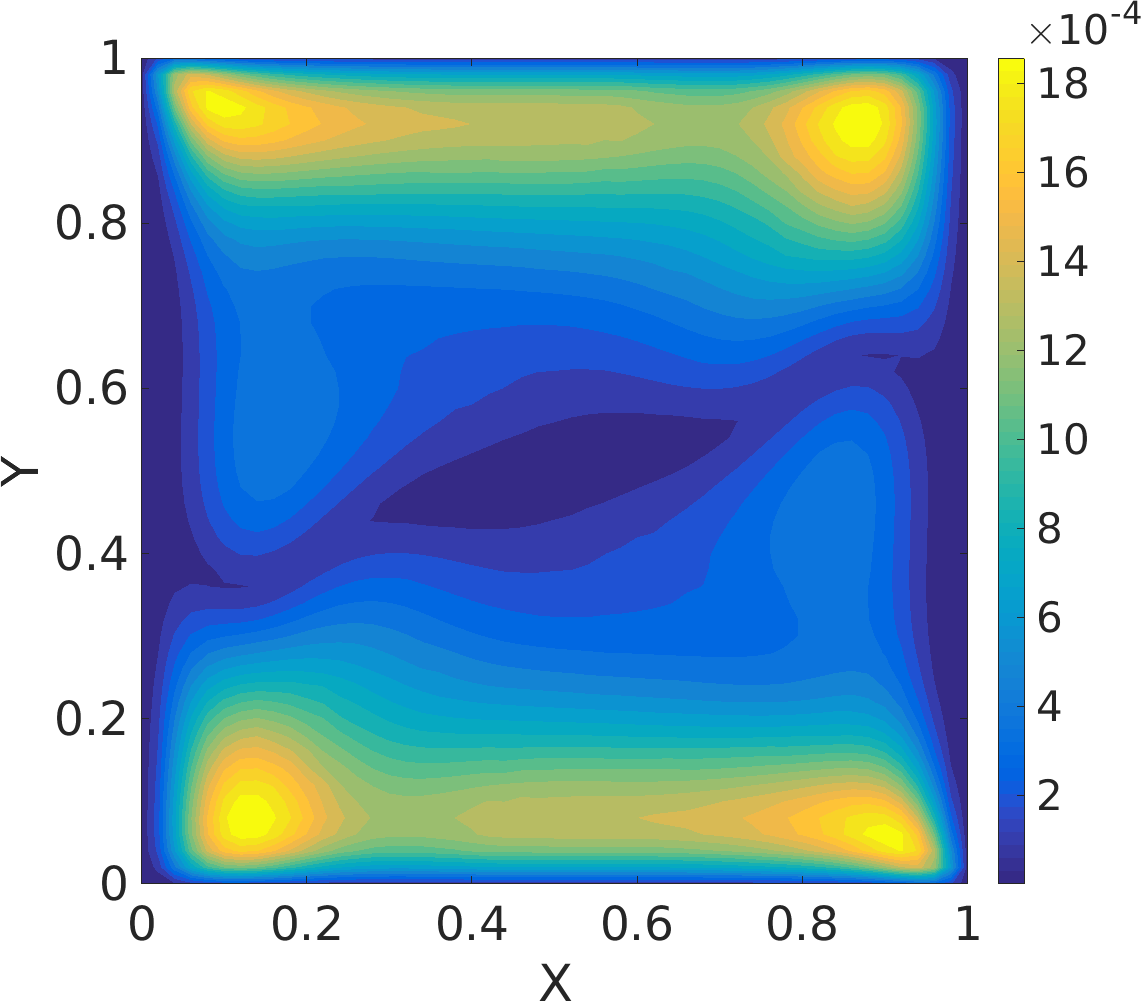}
		\caption{X Velocity}
	\end{subfigure}
	\begin{subfigure}[t]{0.32\textwidth}
		\includegraphics[width=\textwidth]{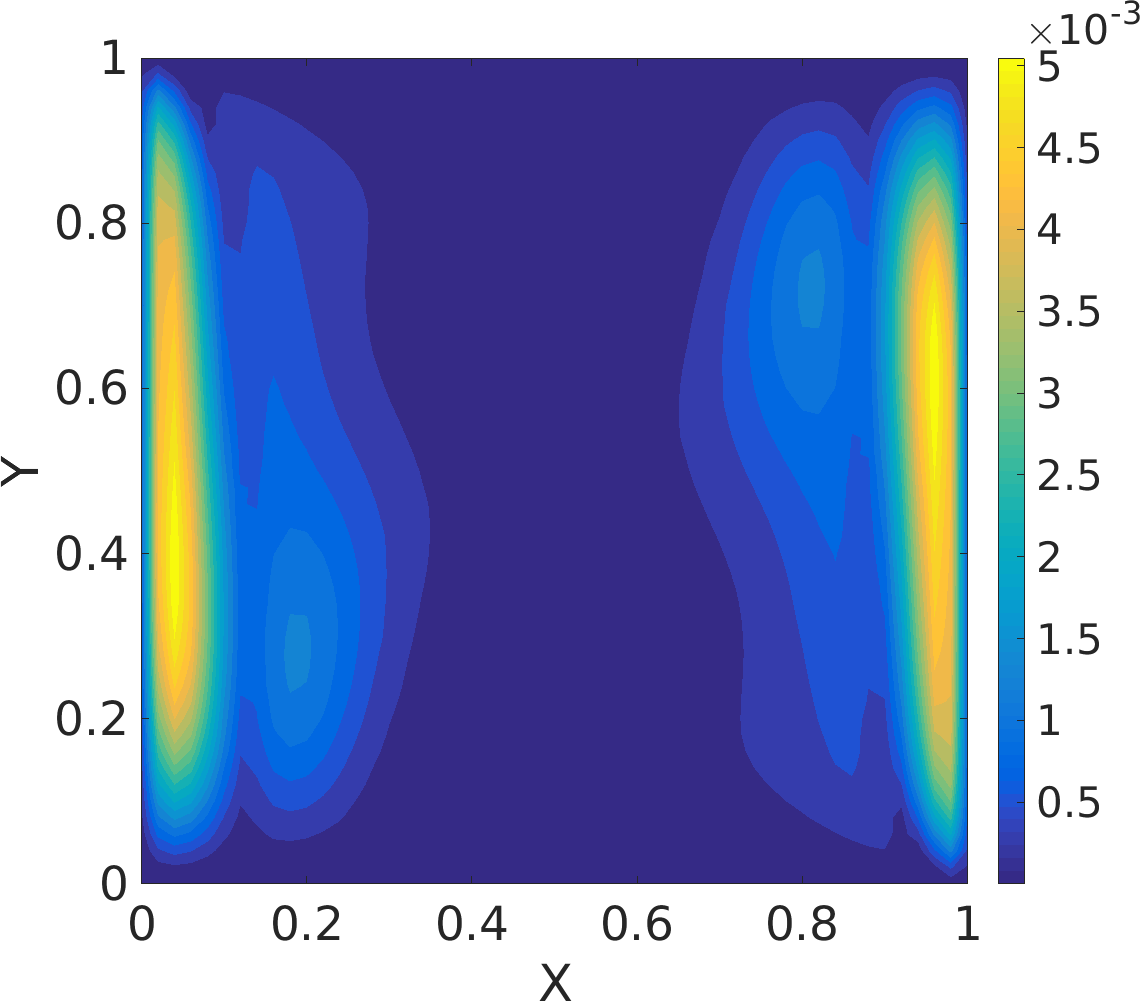}
		\caption{Y Velocity}
	\end{subfigure}
	\caption{Stochastic Standard Deviation at $Z=0.5$ Mid-plane for $Ra=10^6$}
	\label{Fig:Case_A_1E6_Std}
\end{figure}

For the case of temperature gradient in the $X$ direction and gravity in the $Y$ direction, temperature and velocity contours at the $Z$ mid-plane are quite informative. All the stochastic mean contour plots look visually similar to those of the deterministic natural convection problem (\cref{Fig:Ra_1E5_value_at_input_mean,Fig:Ra_1E6_value_at_input_mean}) and thus, are not plotted. Instead, the difference between the stochastic mean and deterministic values is plotted here. \Cref{Fig:Case_A_1E5_Diff,Fig:Case_A_1E5_Std,Fig:Case_A_1E6_Diff,Fig:Case_A_1E6_Std} plot the difference and standard deviation of temperature and velocities for the Rayleigh numbers of $10^5$ and $10^6$, respectively. It is observed that both the difference and standard deviations are of similar orders of magnitude. The difference and standard deviation in the temperature are three orders of magnitude smaller than the mean. On the other hand, the difference and standard deviation in velocities are two orders of magnitude smaller than the mean. Thus, it can be concluded that the effect of uncertainty is significant on the velocities than on the temperature. The  uncertainty has higher impact on the temperature at $Ra=10^6$ simulation than at $Ra=10^5$. Similar to the Nusselt number, the standard deviation is higher when the mean value is higher. The difference between the stochastic mean and the deterministic value for temperature is of the order of $0.1-0.3\%$ of the deterministic value for low and high Rayleigh numbers respectively. The difference is of the order of $1-3\%$ for both the velocities of low and high Rayleigh numbers respectively. Thus, the effect of input stochasticity is higher on the velocities compared to temperature.
\subsubsection{Sensitivity Analysis}
\begin{figure}[H]
	\centering
	\begin{subfigure}[t]{0.49\textwidth}
		\includegraphics[width=\textwidth]{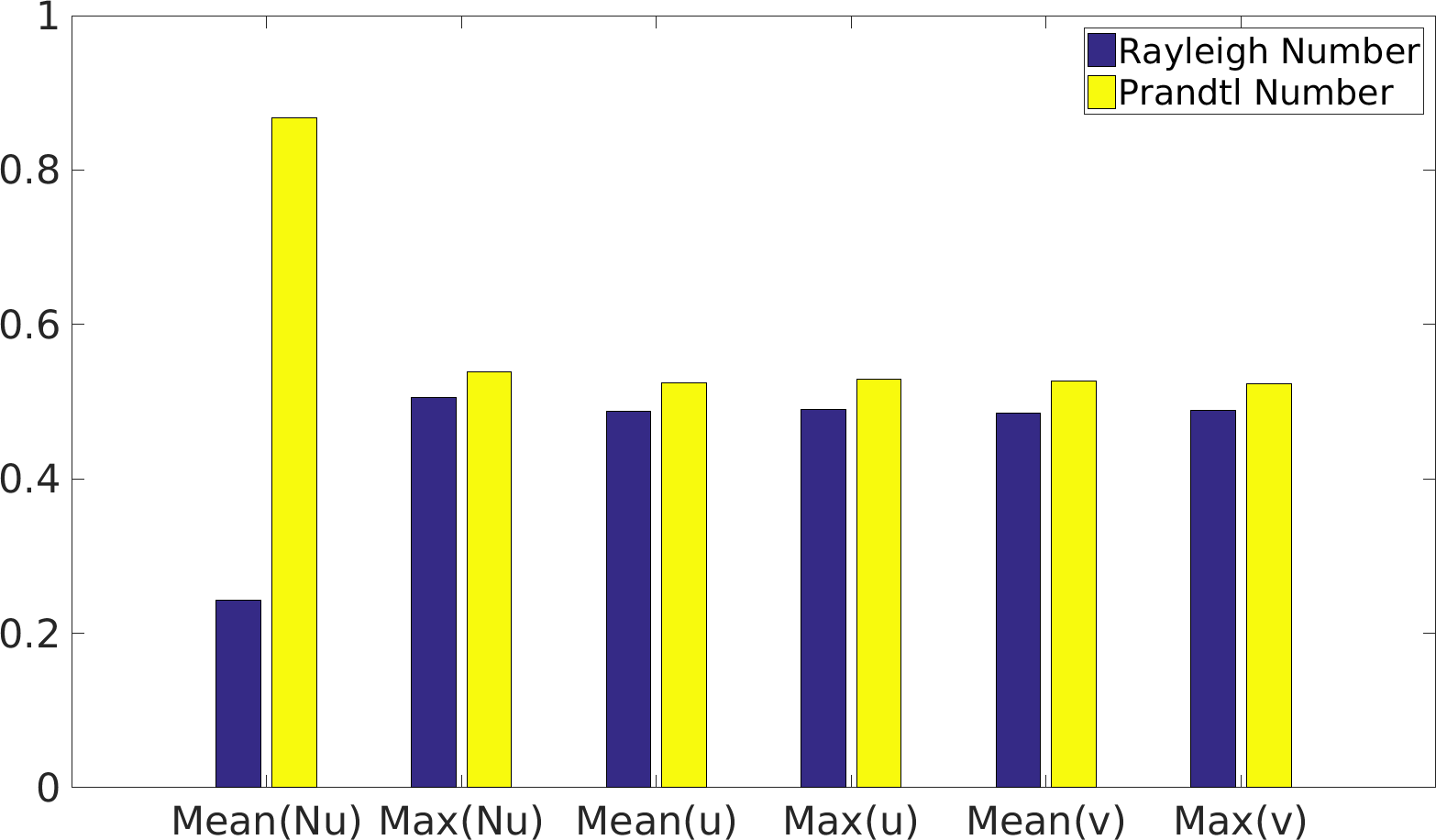}
		\caption{$Ra=10^5$}
		\label{Fig:sensitivity_1E5}
	\end{subfigure}
	\begin{subfigure}[t]{0.49\textwidth}
		\includegraphics[width=\textwidth]{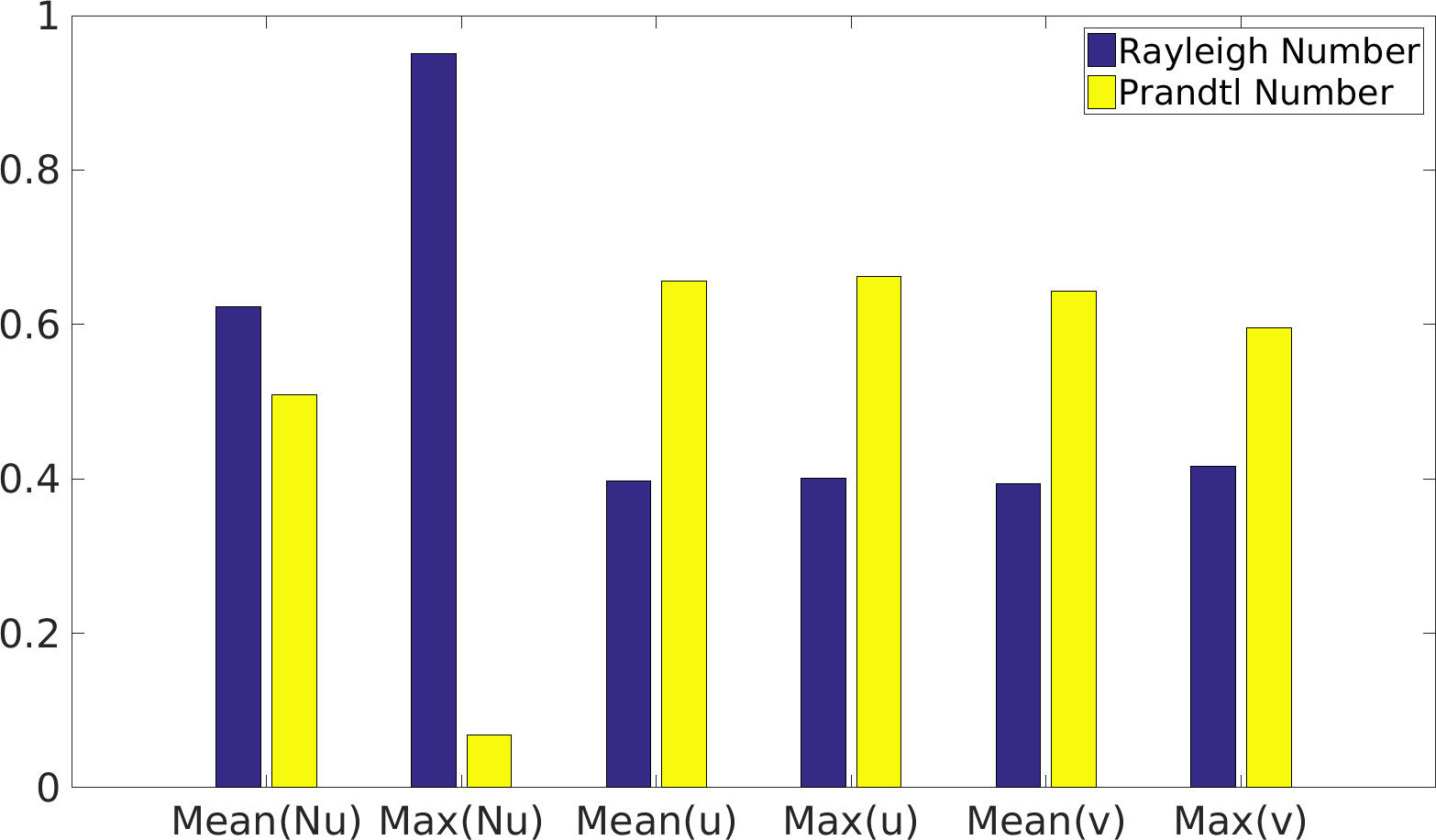}
		\caption{$Ra=10^6$}
		\label{Fig:sensitivity_1E6}
	\end{subfigure}
	\caption{Sensitivity: 6 Outputs, 2 Inputs}
	\label{Fig:sensitivity}
\end{figure}
Sensitivity analysis quantifies the variation in an output due to the variation in a particular input. In this work, the global sensitivity is estimated using the Sobol indices based on the Sobol decomposition \cite{sobol1993sensitivity}. Partial Sobol index measures the contribution of a subset of inputs to the total variance. This includes the variance due to a coupling between the inputs. Total Sobol index for each input is defined as the sum of all the partial indices involving that input parameter. Here, the total Sobol indices are estimated from the polynomial chaos coefficients using the sensitivity analysis tool of the software UQLab \cite{marelli2014uqlab}.
\par Wall Nusselt number and X and Y velocities (u,v) have been chosen as representative outputs of the natural convection problem to study the sensitivity. The mean and maximum are taken over the hot wall and entire cube for the Nusselt number and velocities, respectively. \Cref{Fig:sensitivity} plots the sensitivity (total Sobol index) of the mean and maximum of Nusselt number and velocities with respect to each stochastic input parameter viz. Rayleigh and Prandtl numbers. It can be seen that the Nusselt number is more sensitive to Prandtl number for the case of $Ra=10^5$ whereas, it is more sensitive to Rayleigh number for $Ra=10^6$. On the other hand, the velocities are more sensitive to the Prandtl number in both the cases. But for the higher Rayleigh number case, their sensitivity towards the Prandtl number increases. 
\subsection{Case B}
\subsubsection{Deep Neural Network Training and Testing}\label{Sec:Deep Neural Network Training and Testing}
As described briefly in \Cref{Sec:Deep Neural Networks}, the hot wall is divided into strips in a direction orthogonal to gravity (\cref{Fig:Temperature_BC}). Different boundary temperatures are prescribed on each strip. It is assumed that each strip temperature is an  independent and identically distributed (i.i.d) random variable and follows a normal distribution with $\mu=1.05$ and $3 \sigma=0.01$. Since the cold wall is held at a constant temperature of 0.95, the mean temperature difference driving the natural convection flow is still 0.1. This implies that a $3 \sigma$ error of 10\% is specified in the input stochasticity. The number of strips in this study is varied from 4 to 32 in multiples of 2. The material properties of the fluid are kept constant. In order to estimate the statistics of the outputs, a deep neural network (DNN) surrogate model is used. The boundary temperatures on each strip are inputs to the DNN. A $99^3$ finite volume mesh is used for the numerical simulation. Four separate DNNs are trained with the following outputs:
\begin{enumerate}
	\item Nusselt number along the hot wall ($99^2=9801$ outputs)
	\item Temperature along the $Z=0.5$ midplane ($99^2=9801$ outputs)
	\item $X$ velocity along the $Z=0.5$ midplane ($99^2=9801$ outputs)
	\item $Y$ velocity along the $Z=0.5$ midplane ($99^2=9801$ outputs)
\end{enumerate}
Latin hypercube samples (LHS) are generated using the python package pyDOE \cite{pydoe}. The uniformly distributed LHS are transformed into normal distribution using the inverse cumulative distribution function (ppf) of the statistical functions module of scipy \cite{scipy_stats}. Separate sets of LHS are generated for training, validation and testing. Cases with 4, 8 and 16 strips are trained with 500 samples whereas, for 32 strips, 1000 samples are required. For each case, two different sets of 100 samples are used for validation and testing.
\par
\begin{figure}[H] 
	\centering
	\begin{subfigure}[t]{0.49\textwidth}
		\includegraphics[width=\textwidth]{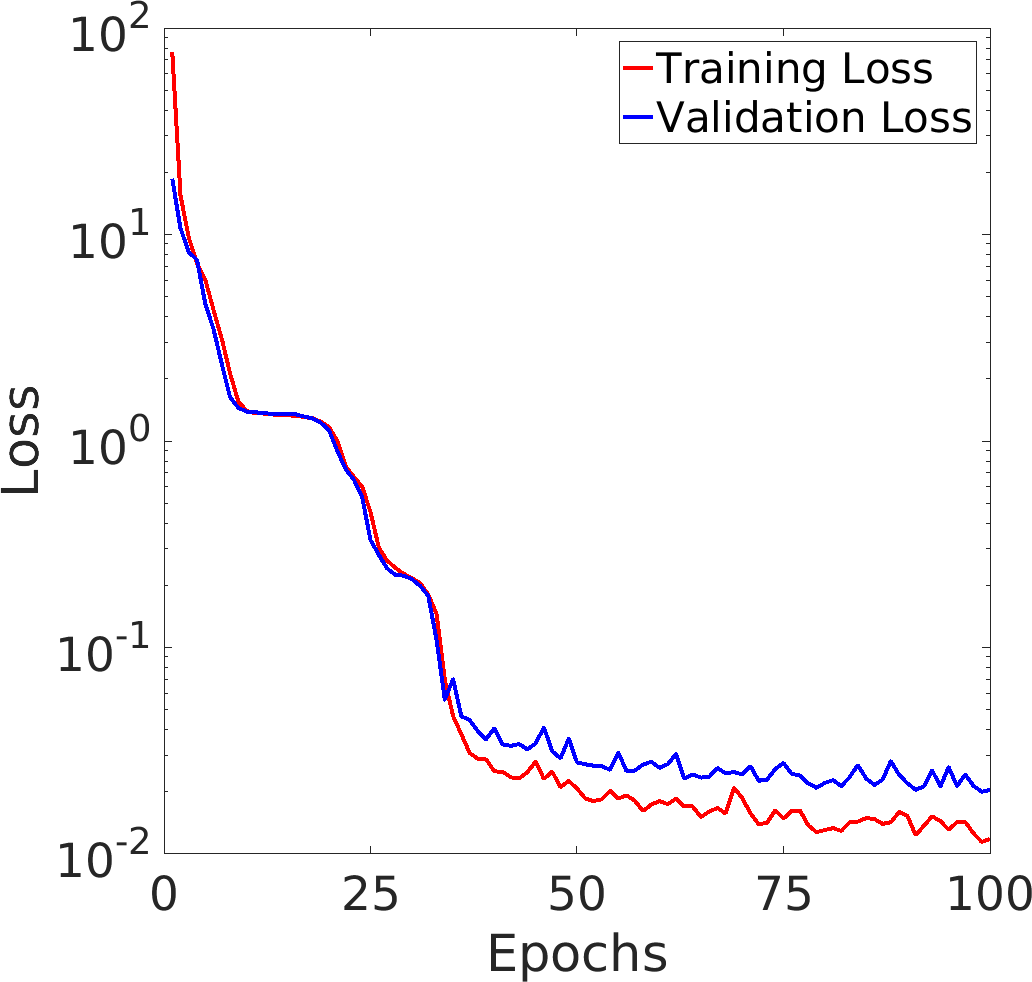}
		\caption{Hot Wall Nusselt Number}
		\label{Fig:DNN Training and Validation Loss Nu}
	\end{subfigure}
	\begin{subfigure}[t]{0.49\textwidth}
		\includegraphics[width=\textwidth]{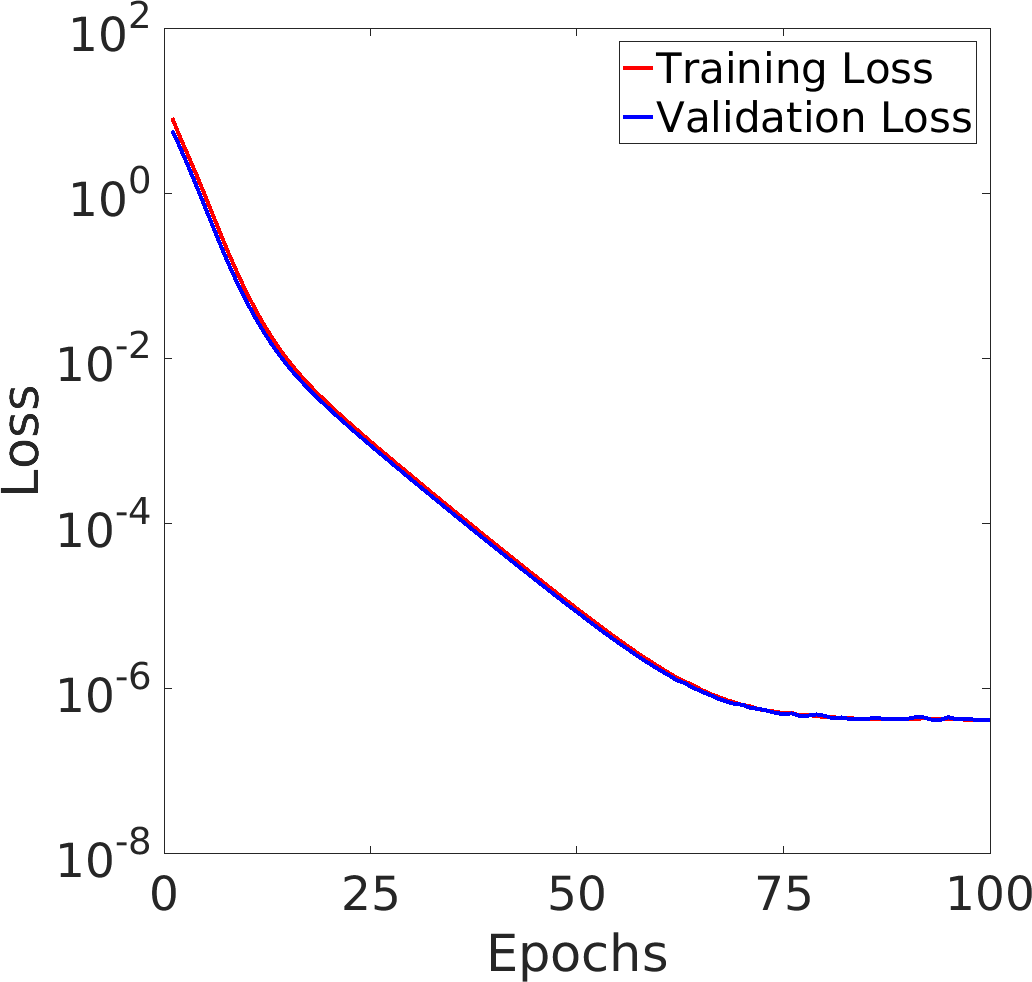}
		\caption{Temperature at $Z=0.5$}
		\label{Fig:DNN Training and Validation Loss T}
	\end{subfigure}
	\begin{subfigure}[t]{0.49\textwidth}
		\includegraphics[width=\textwidth]{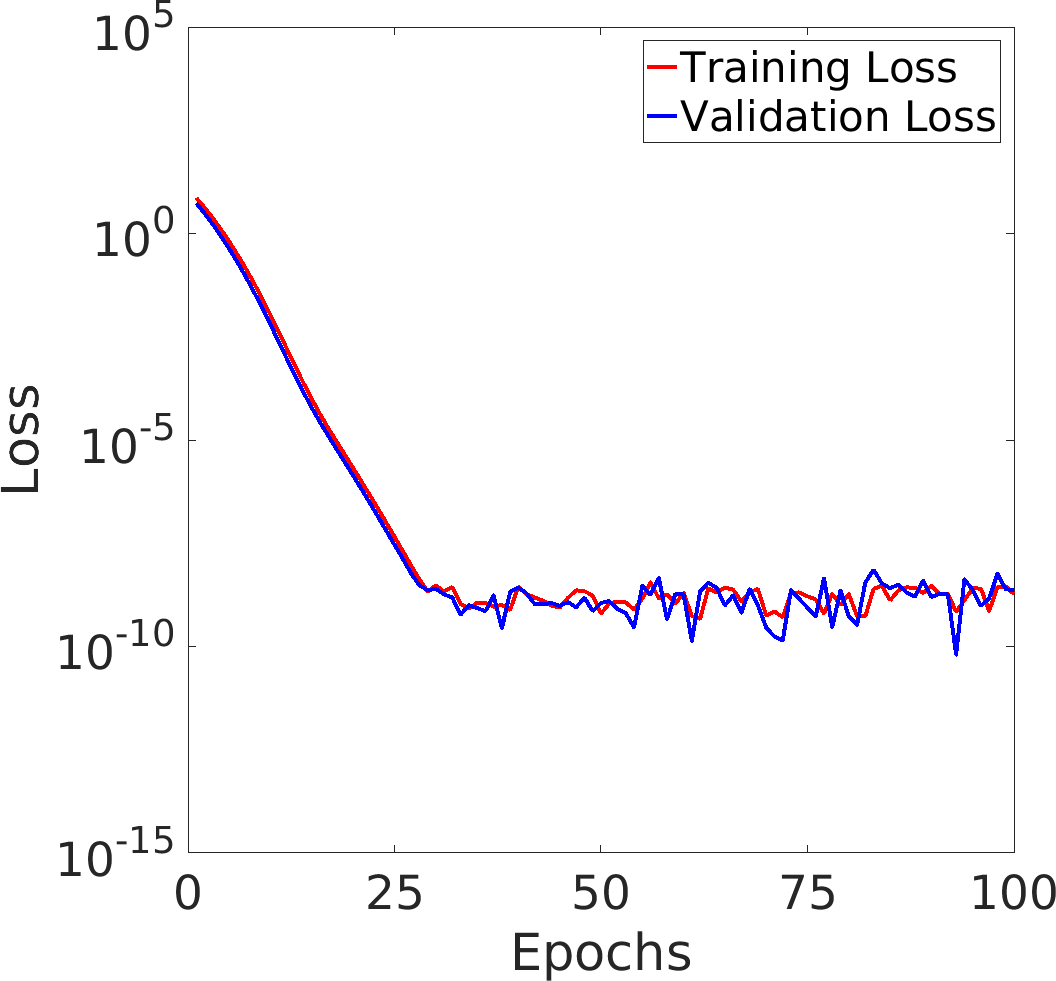}
		\caption{X-Velocity at $Z=0.5$}
		\label{Fig:DNN Training and Validation Loss u}
	\end{subfigure}
	\begin{subfigure}[t]{0.49\textwidth}
		\includegraphics[width=\textwidth]{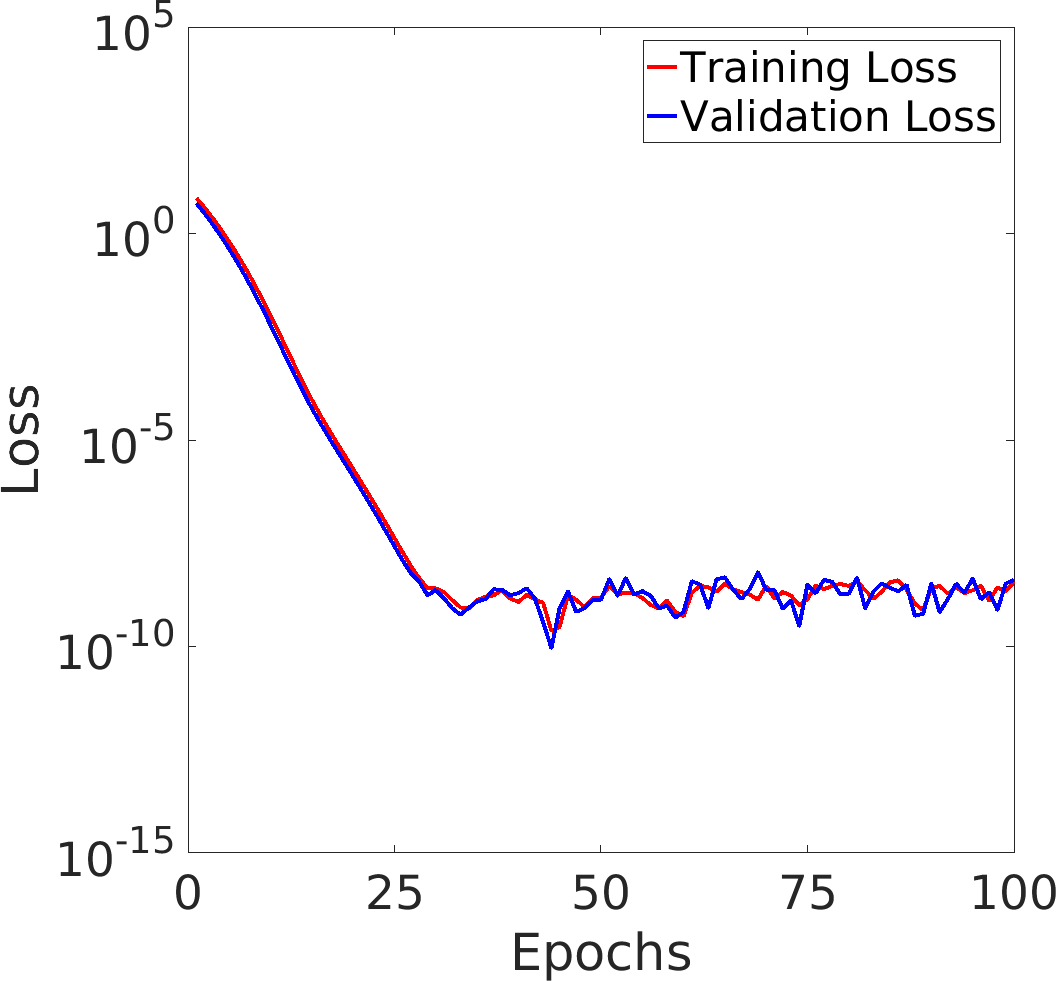}
		\caption{Y-Velocity at $Z=0.5$}
		\label{Fig:DNN Training and Validation Loss v}
	\end{subfigure}
	\caption{DNN Training and Validation Loss for 4 Strips}
	\label{Fig:DNN Training and Validation Loss}
\end{figure}

 The number of neurons in the input and output layers is specified by the number of inputs and outputs. The learning rate, regularization constant, optimizer and the number of hidden layers and neurons are highly problem specific and are chosen so that both the training and validation error are simultaneously minimized. The prediction accuracy is then checked on an unseen testing set. This overall procedure helps in fitting a DNN with low bias and low variance \cite{goodfellow2016deep}. The DNNs are implemented in the Python library Tensorflow \cite{tensorflow2015-whitepaper} with a high level API Keras \cite{chollet2015keras}. Among the various optimizers available in Keras, the Adam optimizer \cite{kingma2014adam} is found most suitable in this work. Settings of Adam optimizer are as follows: learning rate of $10^{-3}$, $\beta_1=0.9$, $\beta_2=0.999$ and `amsgrad' option switched on. ReLU and identity are the activation functions for all the hidden and output layers respectively. Other hyperparameters specific to each of the four DNNs are as follows:
\begin{enumerate}
	\item Wall Nusselt number: $\lambda=0.001$, $L_h=5$, $n_h=300$
	\item Temperature: $\lambda=0.001$, $L_h=4$, $n_h=300$
	\item $X$ velocity: $\lambda=0.01$, $L_h=4$, $n_h=300$
	\item $Y$ velocity: $\lambda=0.01$, $L_h=4$, $n_h=300$
\end{enumerate}
where, $\lambda$ is $L_2$ regularization constant, $L_h$ is the number of hidden layers and $n_h$ is the number of neurons in each hidden layer. All the DNNs are trained for 100 iterations on the entire dataset known as epochs. All the hyperparameters given above are tuned using the validation set with an objective to minimize bias and variance. \Cref{Fig:DNN Training and Validation Loss} plots the loss versus epochs during training for each of the four DNNs for the case with 4 strips. Losses for DNNs of 8, 12 and 16 strips are similar and hence are not plotted here. Since both the training and validation losses are close to each other, it shows that the variance is low.
\par
\begin{figure}[H] 
	\centering
	\begin{subfigure}[t]{0.49\textwidth}
		\includegraphics[width=\textwidth]{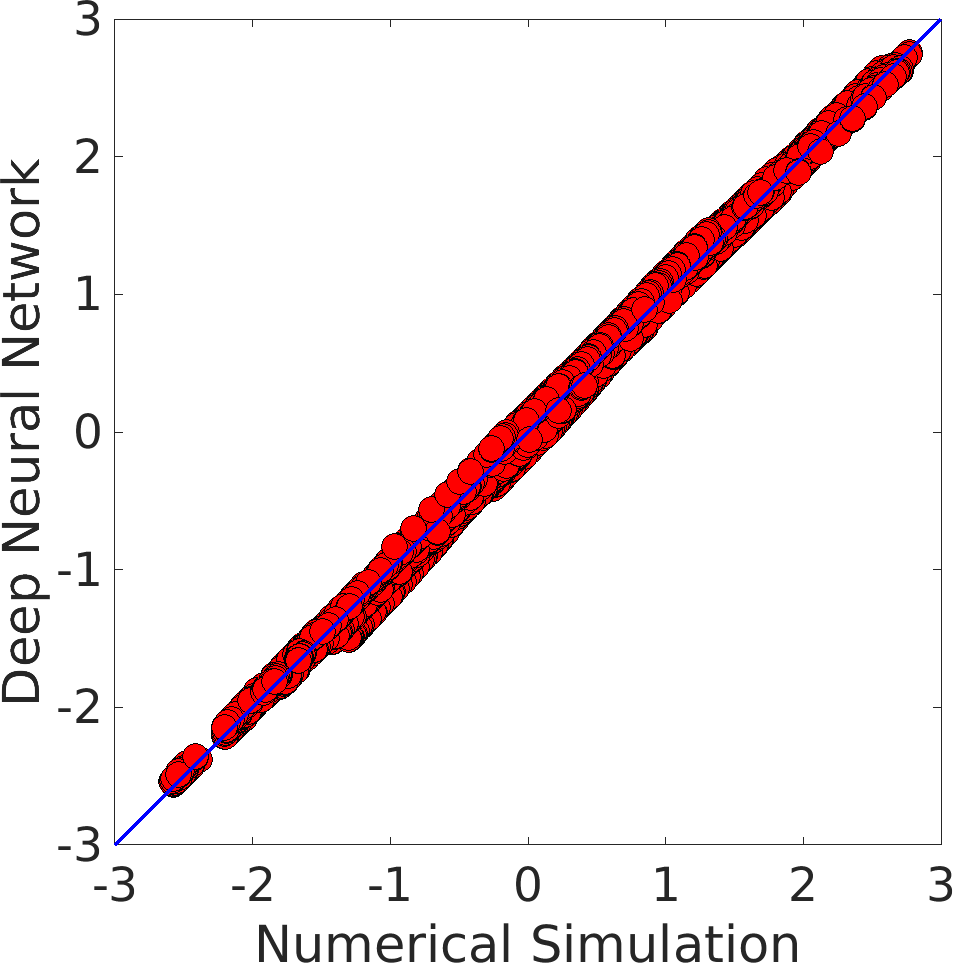}
		\caption{Hot Wall Nusselt Number}
		\label{Fig:CFD_DNN Nu}
	\end{subfigure}
	\begin{subfigure}[t]{0.49\textwidth}
		\includegraphics[width=\textwidth]{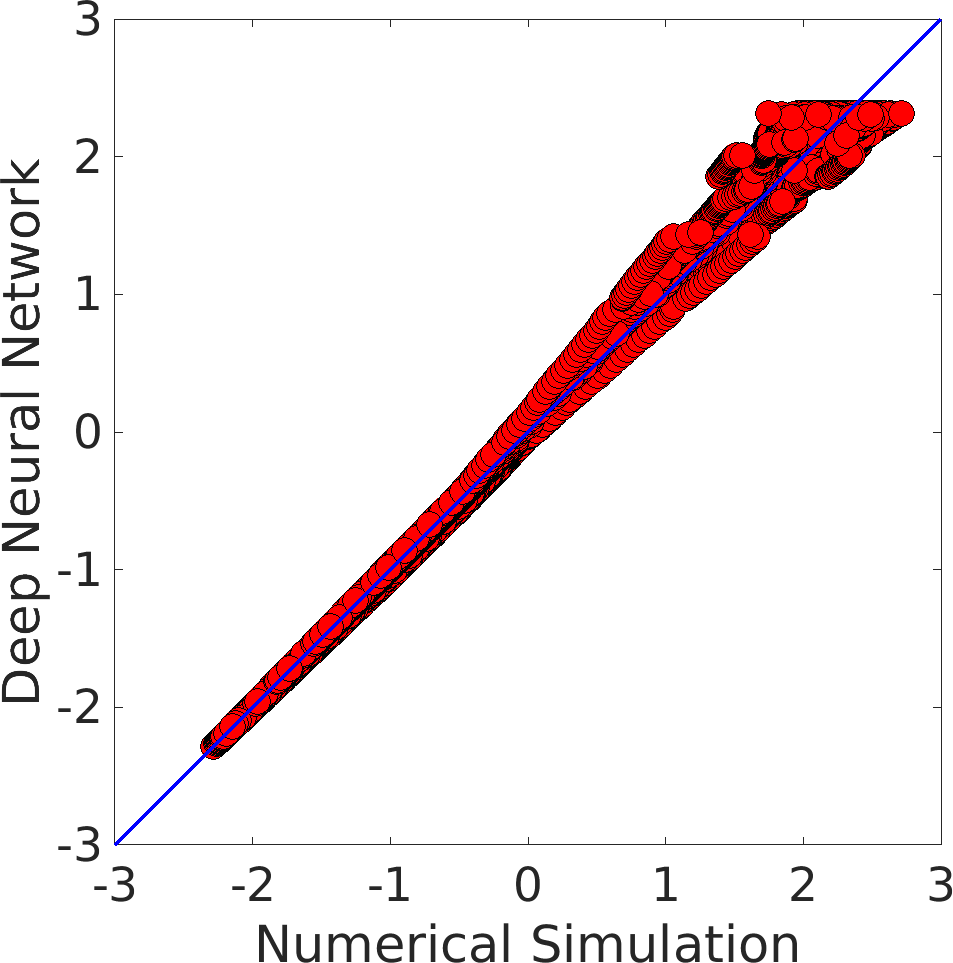}
		\caption{Temperature at $Z=0.5$}
		\label{Fig:CFD_DNN T}
	\end{subfigure}
	\begin{subfigure}[t]{0.49\textwidth}
		\includegraphics[width=\textwidth]{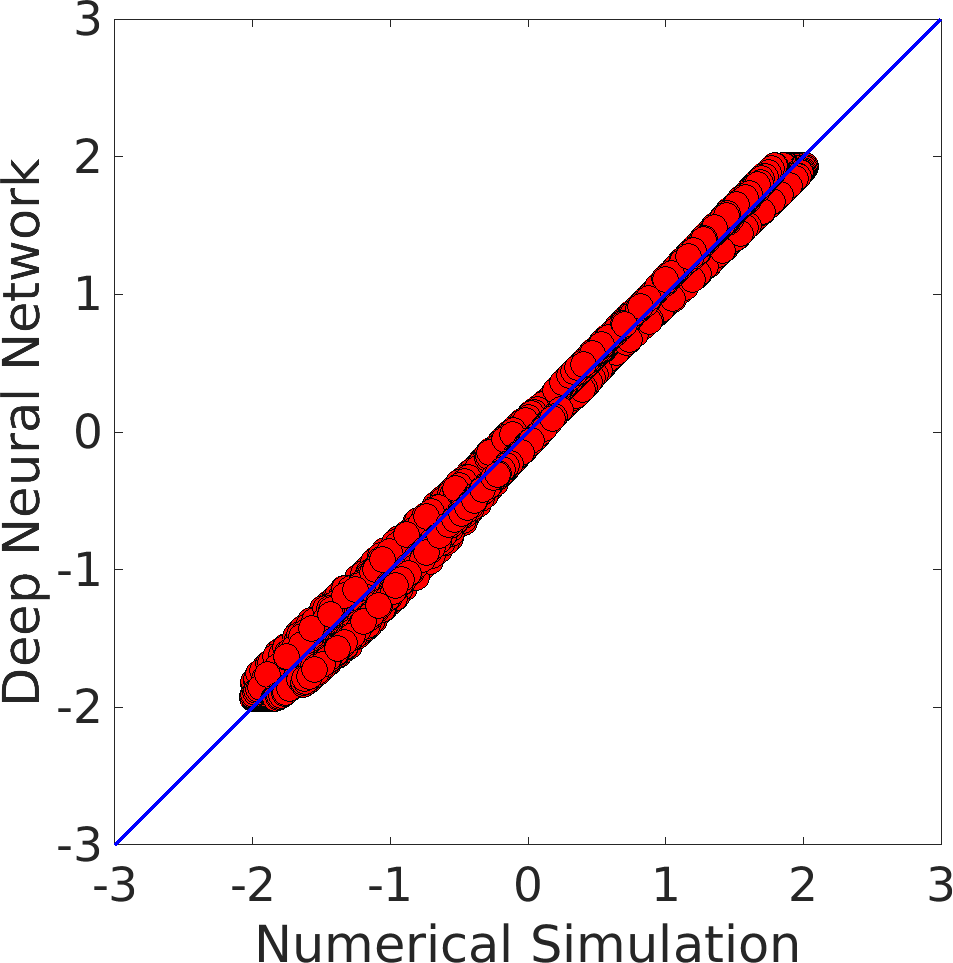}
		\caption{X-Velocity at $Z=0.5$}
		\label{Fig:CFD_DNN u}
	\end{subfigure}
	\begin{subfigure}[t]{0.49\textwidth}
		\includegraphics[width=\textwidth]{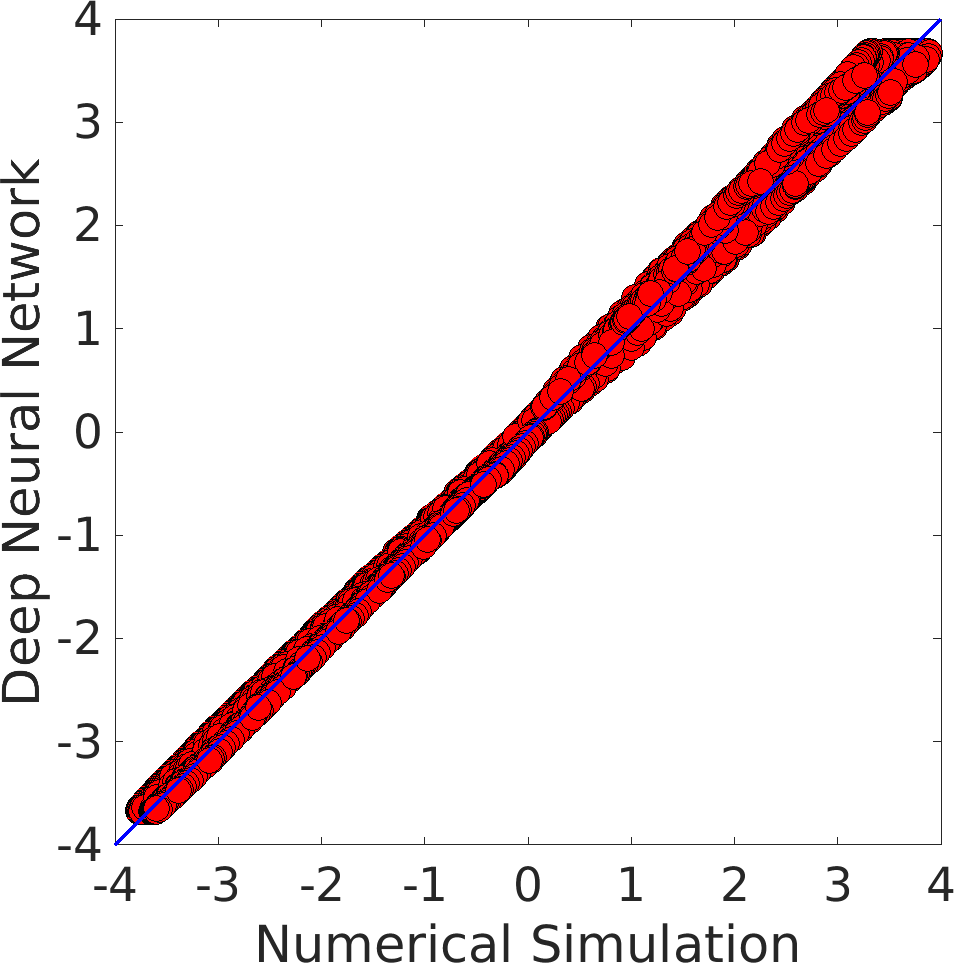}
		\caption{Y-Velocity at $Z=0.5$}
		\label{Fig:CFD_DNN v}
	\end{subfigure}
	\vspace{-0.15cm}
	\caption{Estimate from Numerical Simulation and DNN for 4 Strips}
	\label{Fig:CFD_DNN}
\end{figure}

 \Cref{Tab:DNN_Train_Test_Error} documents the training and testing errors for all the 16 DNNs: 4 DNNs each for 4, 8, 16, and 32 strips. The relative average percent error is defined as hundred times the $L_1$ norm of the difference between true values (numerical simulation) and DNN estimates divided by the maximum absolute value of the output. The training and testing errors are small and close enough thus, implying low bias and variance. For visual inspection, estimates from the numerical simulation and the DNN are plotted together in \cref{Fig:CFD_DNN} for the testing samples for the case of 4 strips. Each output is normalized by subtracting its mean and dividing by its standard deviation and thus, is non-dimensional. Ideally, all the points should lie on the $Y=X$ line but due to the interpolation error, some points are off the line. Since most of the points follow the expected trend of the $Y=X$ line, it can be concluded that the neural network surrogate is accurate.
\begin{table}[h]
	\centering
	\begin{tabular}{|c|c|c|c|c|c|c|c|c|}
		\hline
		& \multicolumn{2}{c|}{4 Strips} & \multicolumn{2}{c|}{8 Strips} & \multicolumn{2}{c|}{16 Strips} & \multicolumn{2}{c|}{32 Strips} \\ \hline
		& Train & Test & Train & Test & Train & Test & Train & Test \\ \hline
		Wall Nusselt No. & 0.302 & 0.390 & 0.205 & 0.512 & 0.602 & 2.260 & 1.393 & 2.373 \\ \hline
		Temperature & 0.029 & 0.030 & 0.024 & 0.025 & 0.022 & 0.023 & 0.022 & 0.024 \\ \hline
		X Velocity & 1.114 & 1.144 & 0.864 & 0.903 & 0.704 & 0.731 & 0.644 & 0.722 \\ \hline
		Y Velocity & 0.372 & 0.387 & 0.297 & 0.311 & 0.248 & 0.260 & 0.233 & 0.230 \\ \hline
	\end{tabular}
	\caption{Relative Average Percent Error: DNN Training and Testing}
	\label{Tab:DNN_Train_Test_Error}
\end{table}
\subsubsection{Nusselt Number}
\begin{figure}[H]
	\centering
	\begin{subfigure}[t]{0.49\textwidth}
		\includegraphics[width=\textwidth]{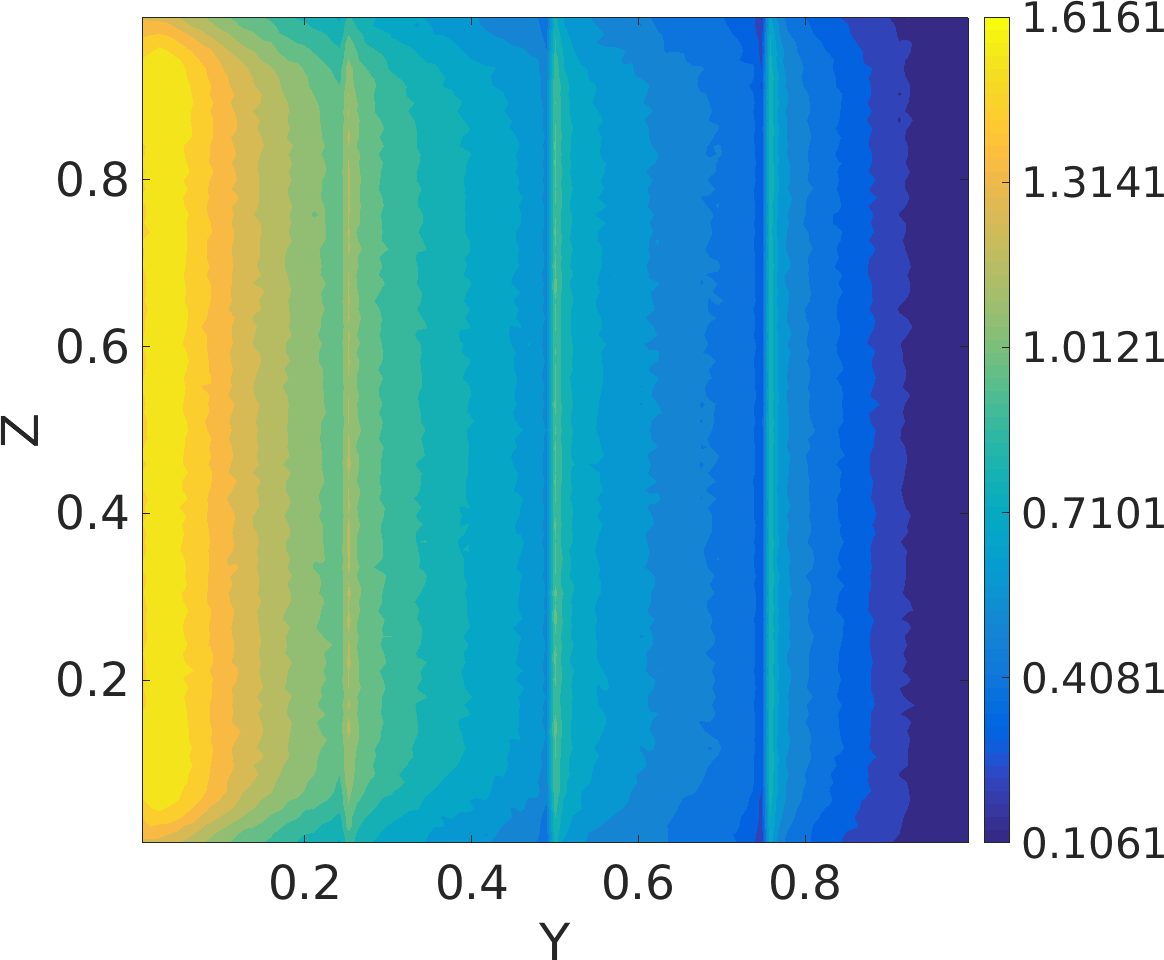}
		\caption{4 Strips}
		\label{Fig:Case B Nusselt Diff 4 Strips}
	\end{subfigure}
	\begin{subfigure}[t]{0.49\textwidth}
		\includegraphics[width=\textwidth]{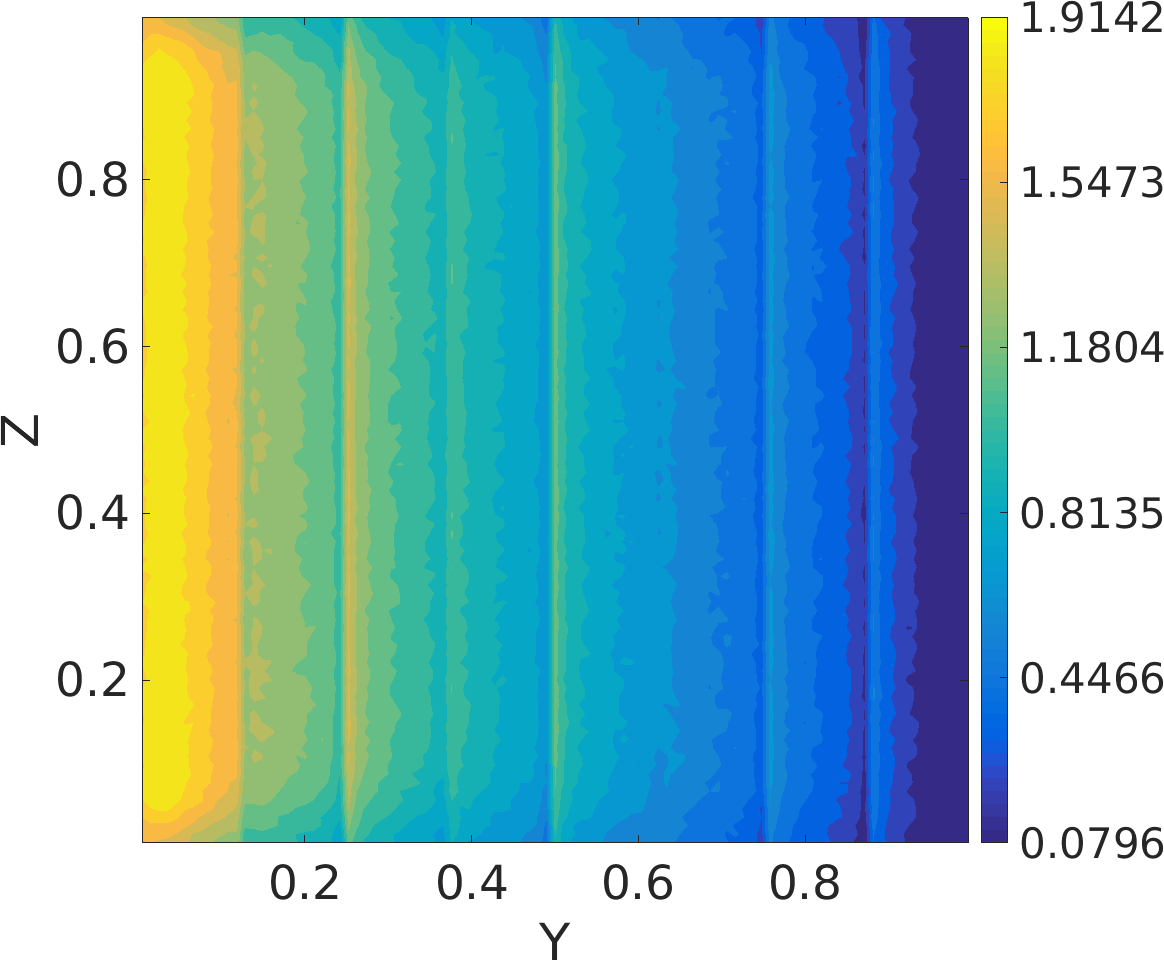}
		\caption{8 Strips}
		\label{Fig:Case B Nusselt Diff 8 Strips}
	\end{subfigure}
	\begin{subfigure}[t]{0.49\textwidth}
		\includegraphics[width=\textwidth]{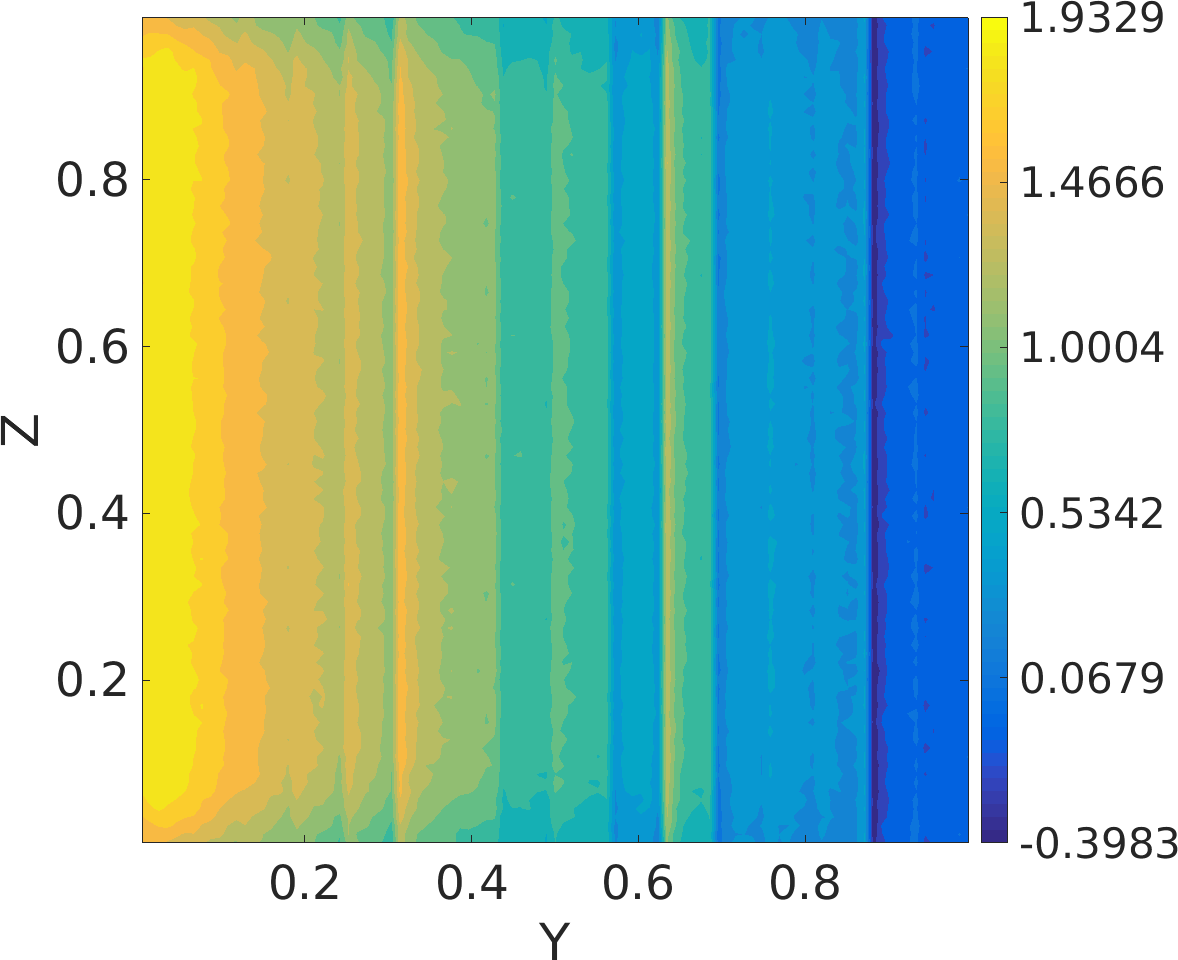}
		\caption{16 Strips}
		\label{Fig:Case B Nusselt Diff 16 Strips}
	\end{subfigure}
	\begin{subfigure}[t]{0.49\textwidth}
		\includegraphics[width=\textwidth]{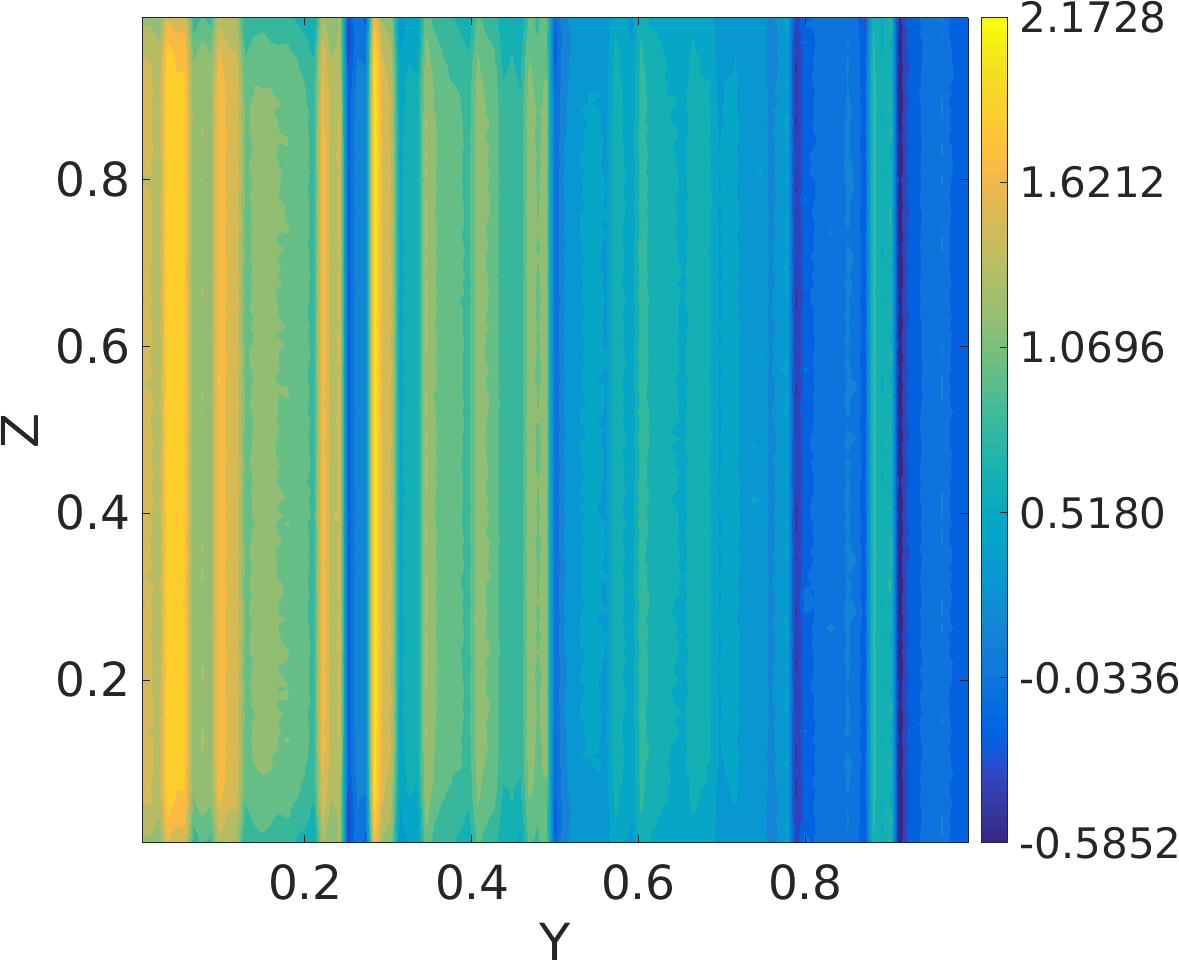}
		\caption{32 Strips}
		\label{Fig:Case B Nusselt Diff 32 Strips}
	\end{subfigure}	
	\caption{$Ra=10^6$ Hot Wall Nusselt Number: Difference between Stochastic Mean and Deterministic Value}
	\label{Fig:Case B Nusselt Diff}
\end{figure}
\begin{figure}[H]
	\centering
	\begin{subfigure}[t]{0.49\textwidth}
		\includegraphics[width=\textwidth]{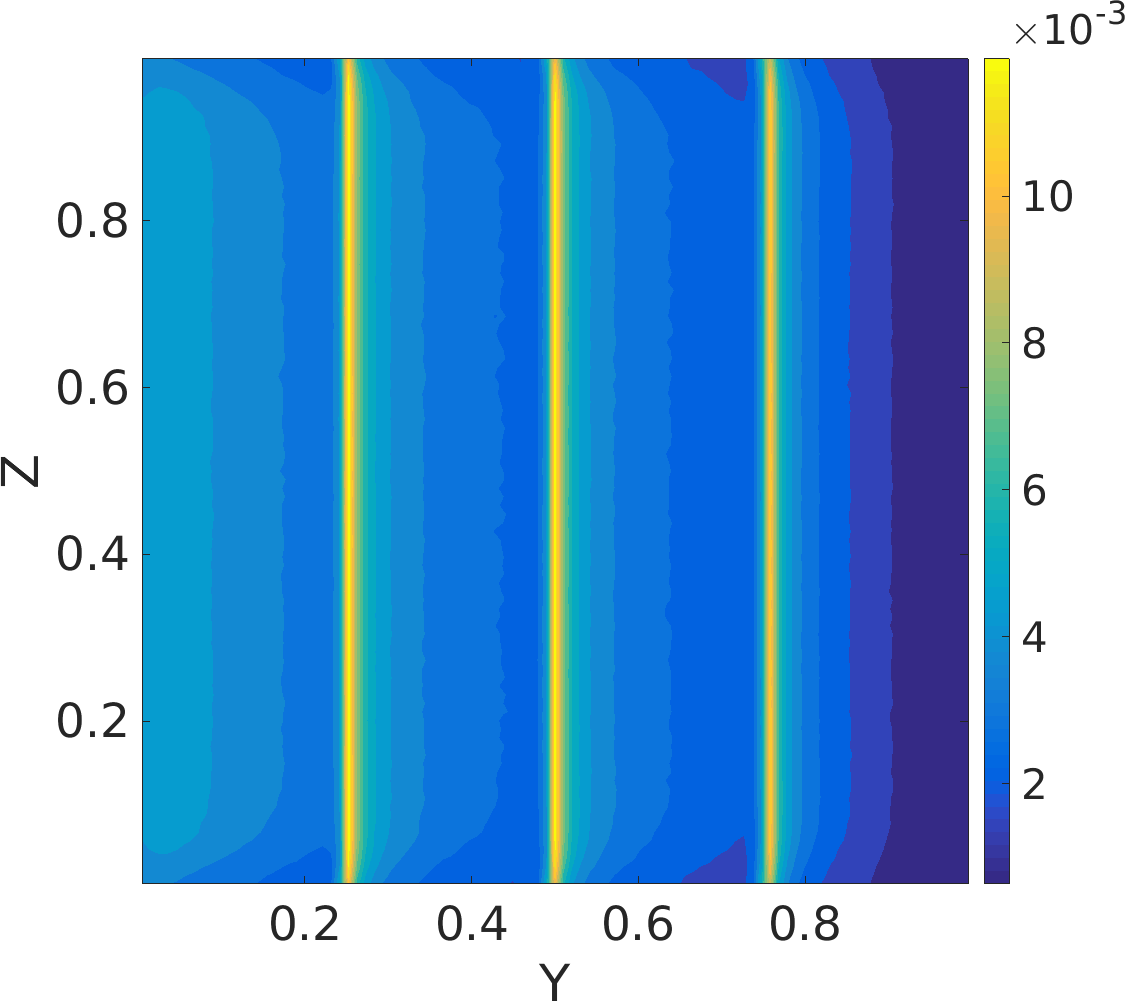}
		\caption{4 Strips}
		\label{Fig:Case B Nusselt Std 4 Strips}
	\end{subfigure}
	\begin{subfigure}[t]{0.49\textwidth}
		\includegraphics[width=\textwidth]{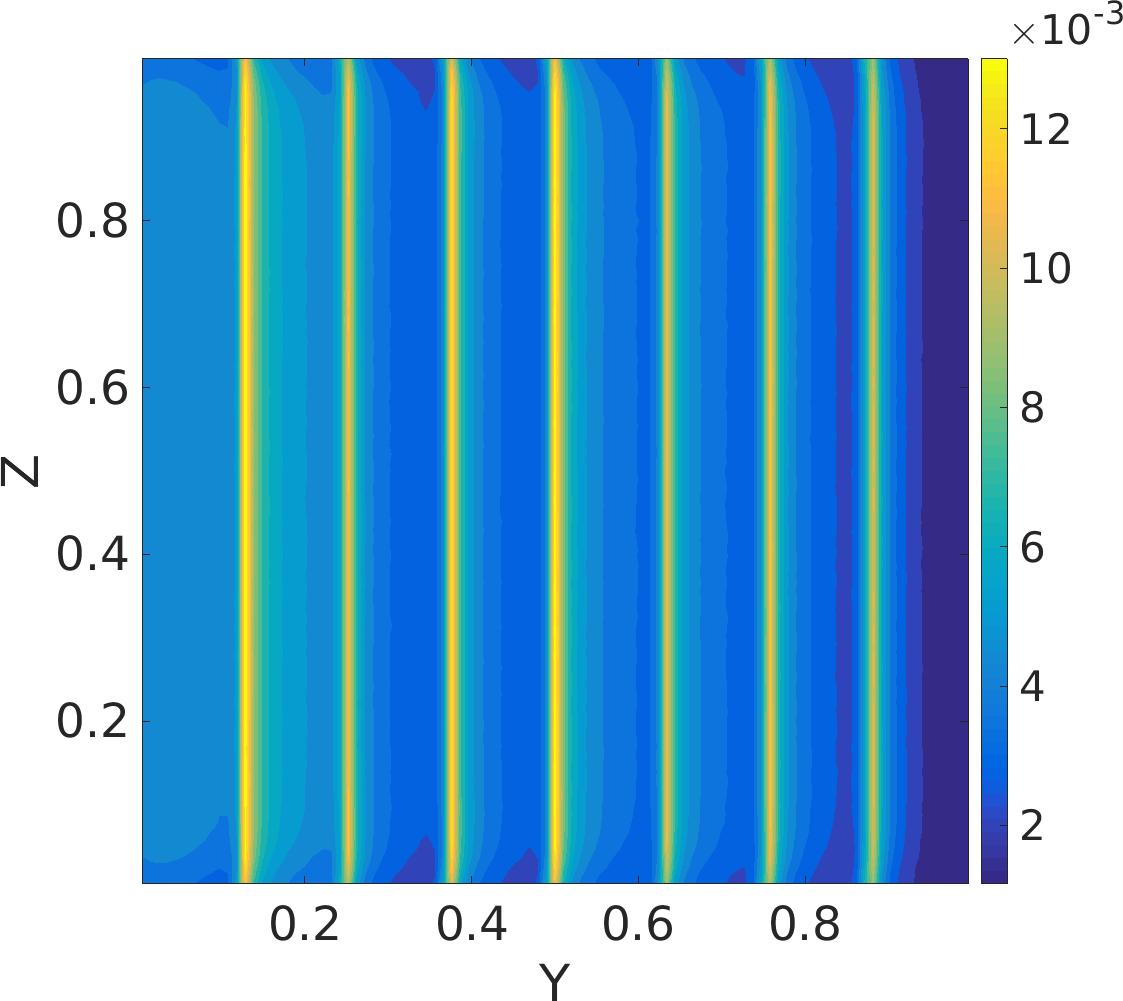}
		\caption{8 Strips}
		\label{Fig:Case B Nusselt Std 8 Strips}
	\end{subfigure}
	\begin{subfigure}[t]{0.49\textwidth}
		\includegraphics[width=\textwidth]{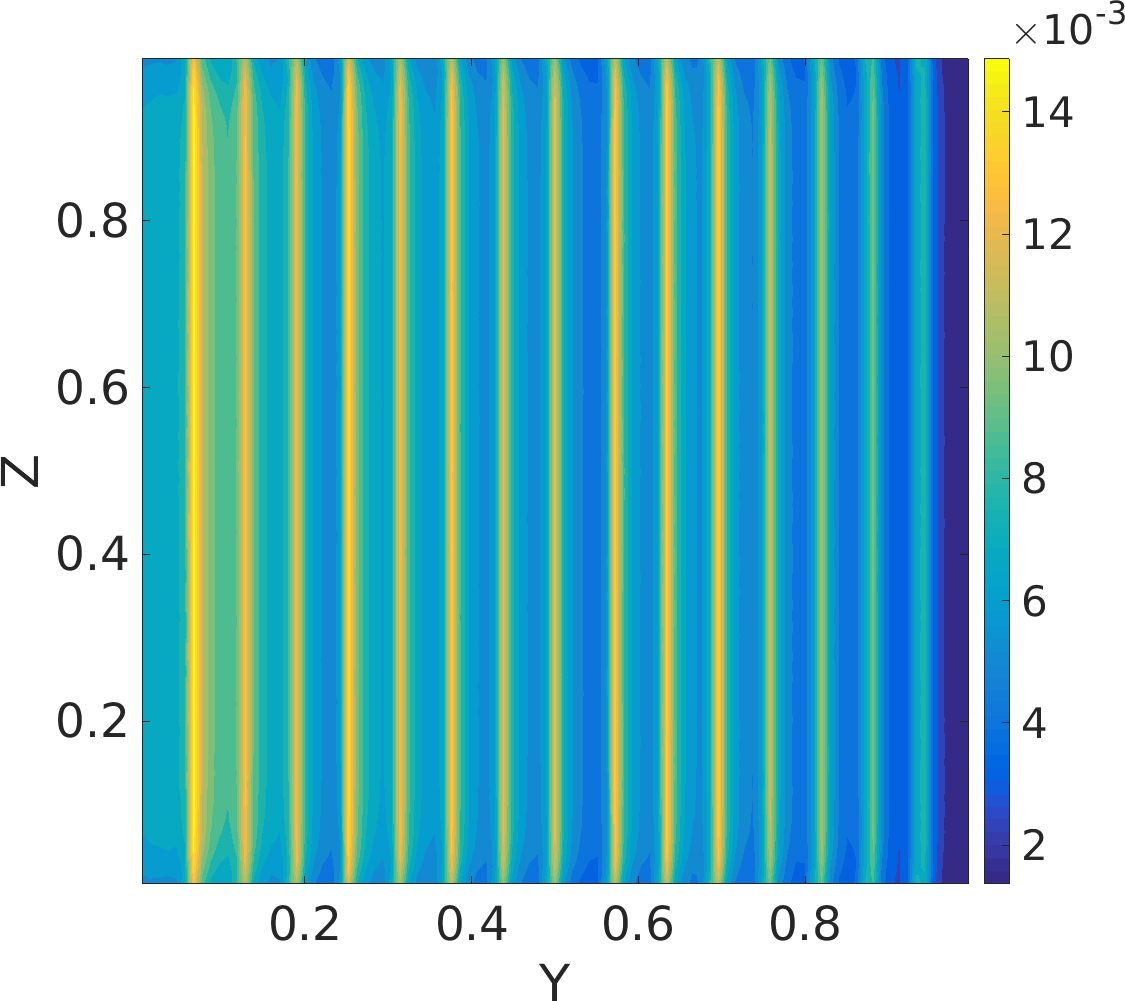}
		\caption{16 Strips}
		\label{Fig:Case B Nusselt Std 16 Strips}
	\end{subfigure}
	\begin{subfigure}[t]{0.49\textwidth}
		\includegraphics[width=\textwidth]{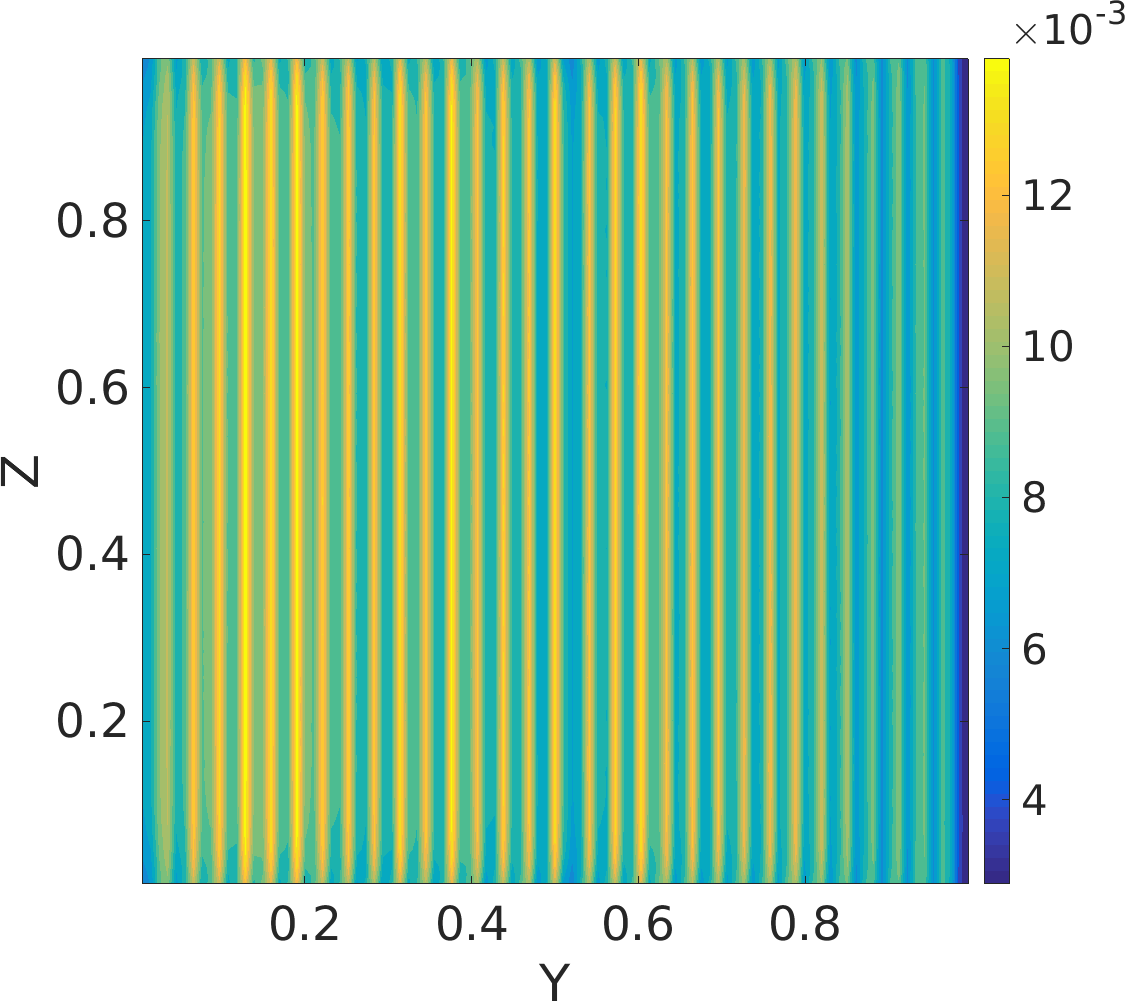}
		\caption{32 Strips}
		\label{Fig:Case B Nusselt Std 32 Strips}
	\end{subfigure}
	\caption{$Ra=10^6$ Hot Wall Nusselt Number: Stochastic Standard Deviation}
	\label{Fig:Case B Nusselt Std}
\end{figure}

\vspace{-0.125cm}
\Cref{Fig:Case B Nusselt Diff} plots the difference between the stochastic mean and the deterministic value of the Nusselt number on the hot wall. Since this wall is subjected to the stochastic boundary condition, the demarcations of the strips can be seen. For example, \cref{Fig:Case B Nusselt Diff 4 Strips,Fig:Case B Nusselt Diff 8 Strips} have four and eight strips respectively. Similar strips are also observed on the stochastic standard deviation contours (\cref{Fig:Case B Nusselt Std}). The maximum deterministic value of the Nusselt number over the hot wall is 18.71 (\cref{Fig:Ra_1E6_value_at_input_mean_Nu}). The maximum values of the difference are 1.62, 1.91, 1.93 and 2.17 for 4, 8, 16 and 32 strips respectively. Hence, this shift in the stochastic mean from the deterministic value is comparable to the input uncertainty of 10\%. On the other hand, the stochastic standard deviation is of the order of $10^{-2}$. Thus, it is seen that the shift in mean is more pronounced than its deviation similar to case A.
\subsubsection{Velocity and Temperature}
\begin{figure}[H]
	\centering
	\begin{subfigure}[t]{0.49\textwidth}
		\includegraphics[width=\textwidth]{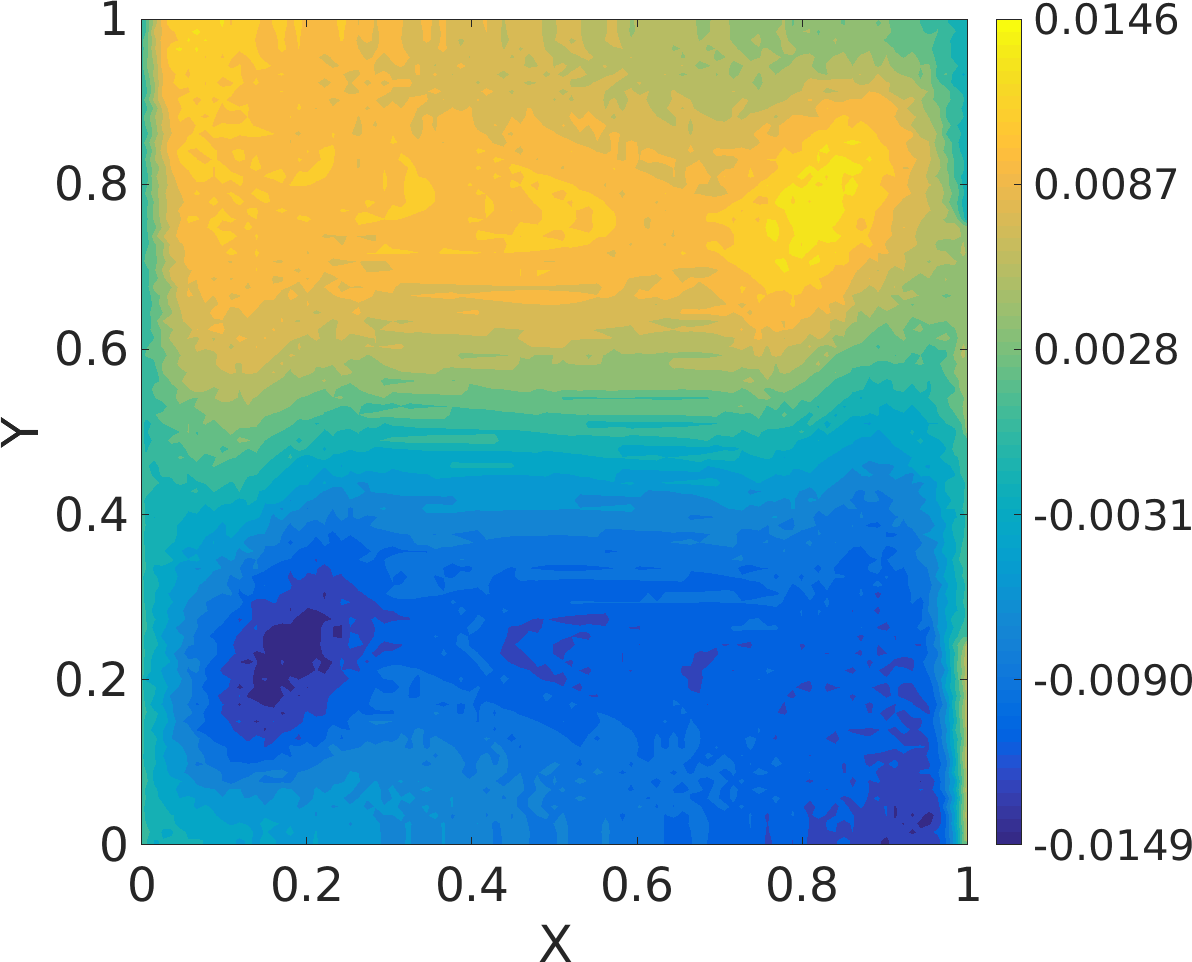}
		\caption{4 Strips}
		\label{Fig:Case B T Diff 4 Strips}
	\end{subfigure}
	\begin{subfigure}[t]{0.49\textwidth}
		\includegraphics[width=\textwidth]{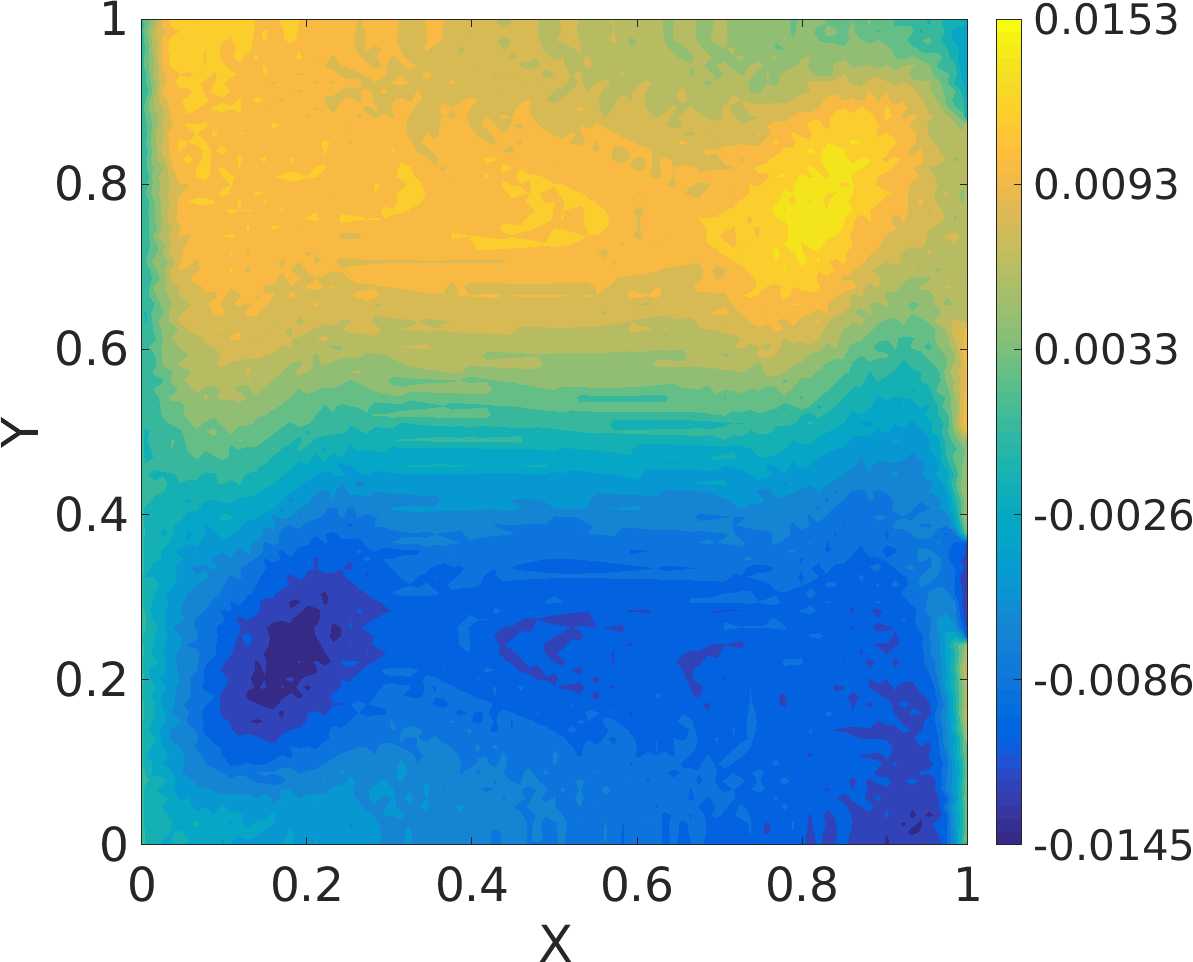}
		\caption{8 Strips}
		\label{Fig:Case B T Diff 8 Strips}
	\end{subfigure}
	\begin{subfigure}[t]{0.49\textwidth}
		\includegraphics[width=\textwidth]{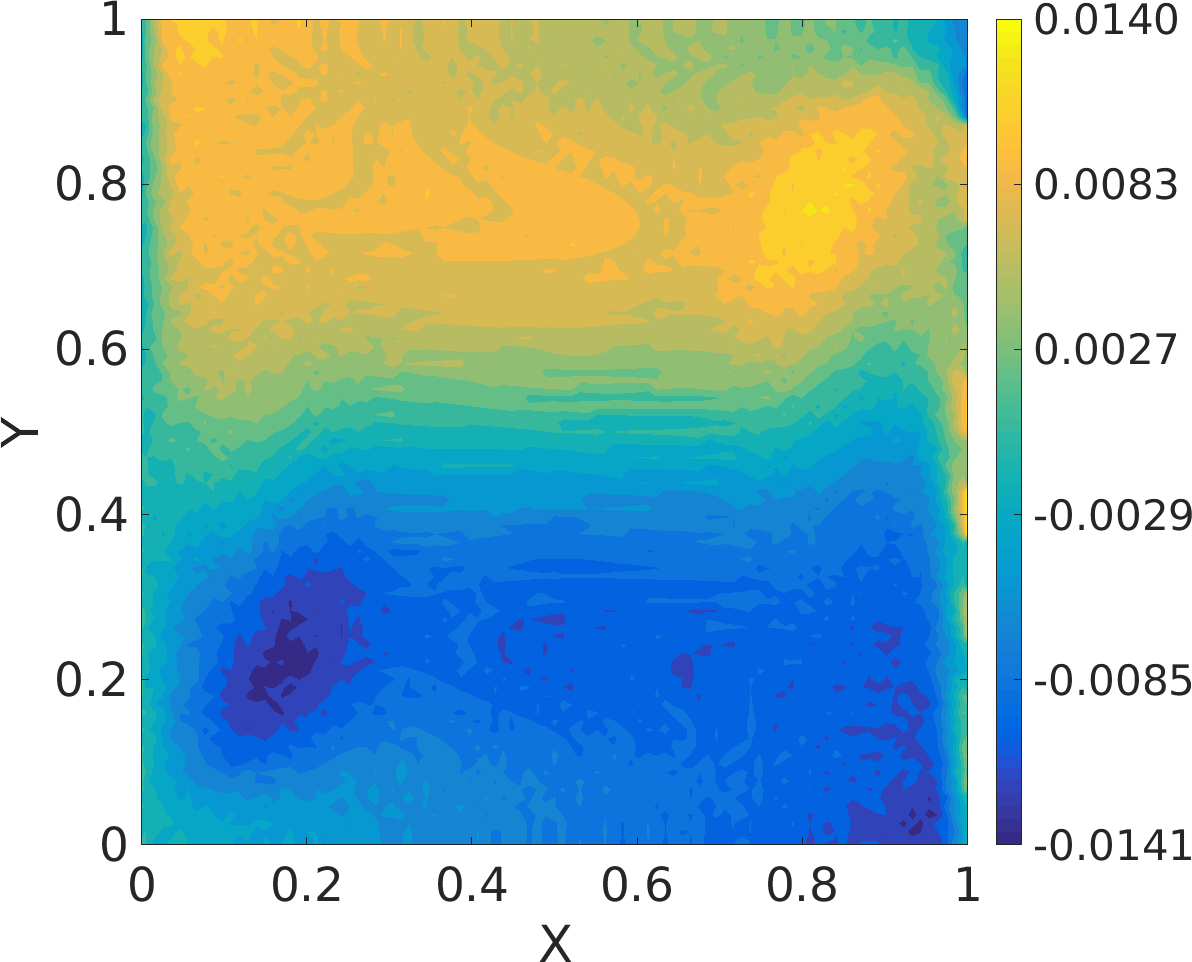}
		\caption{16 Strips}
		\label{Fig:Case B T Diff 16 Strips}
	\end{subfigure}
	\begin{subfigure}[t]{0.49\textwidth}
		\includegraphics[width=\textwidth]{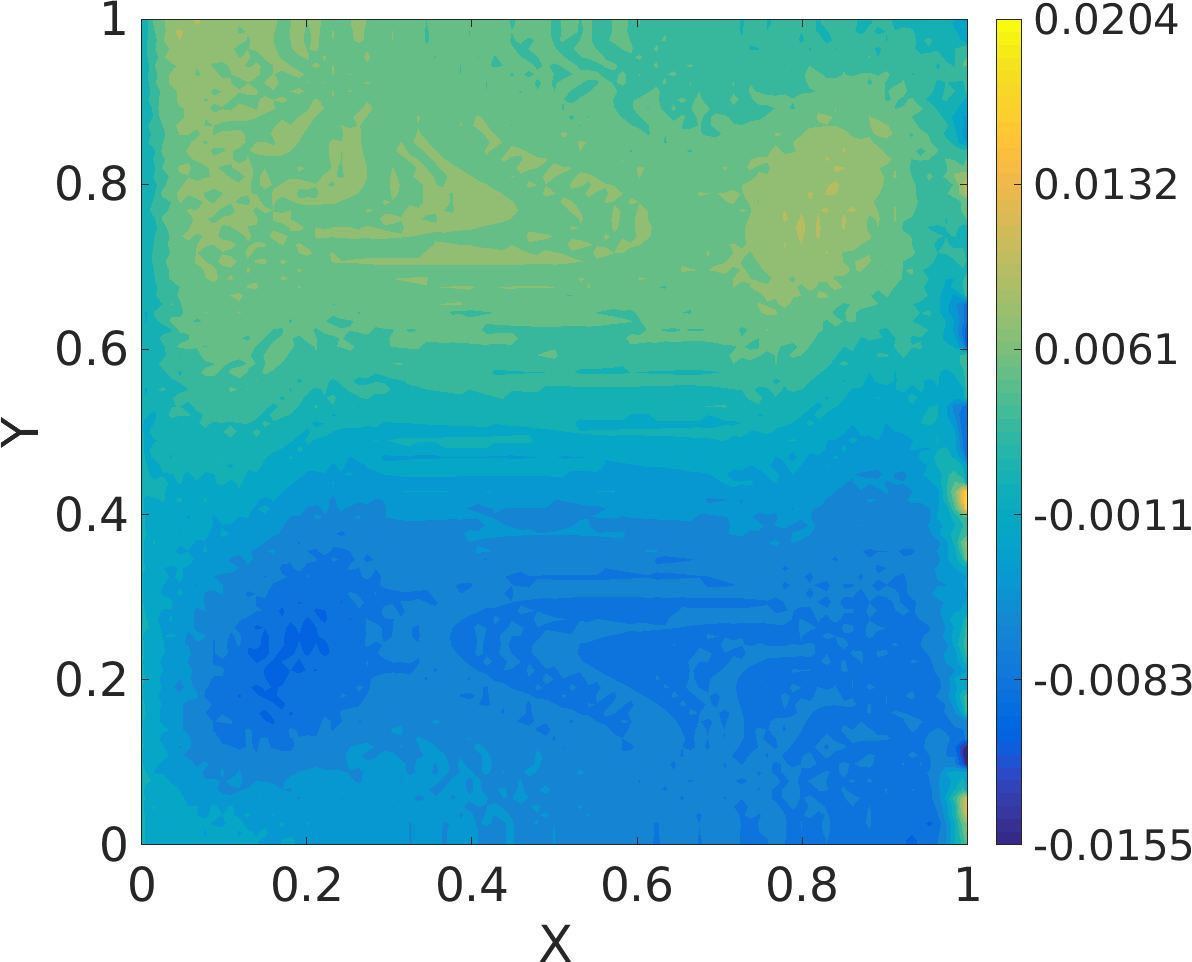}
		\caption{32 Strips}
		\label{Fig:Case B T Diff 32 Strips}
	\end{subfigure}	
	\caption{$Ra=10^6$ $Z$ Midplane Temperature: Difference between Stochastic Mean and Deterministic Value}
	\label{Fig:Case B T Diff}
\end{figure}
\Cref{Fig:Case B T Diff,Fig:Case B u Diff,Fig:Case B v Diff} plot the difference between the stochastic mean and the deterministic value of the temperature, X and Y velocities respectively along the $Z=0.5$ midplane. The effect of the number of strips along the boundary condition is clearly visible in the Nusselt number plots (\cref{Fig:Case B Nusselt Diff,Fig:Case B Nusselt Std}). On the other hand, the number of strips does not affect the temperature and velocity contours. The temperature and velocity contour plots for all the 4 cases look similar. The difference between the stochastic mean and the deterministic value of the temperature is of the order of $3-4\%$ of the deterministic value. For velocities the difference is around $4-5\%$ of the deterministic value. This is smaller than that of the Nusselt number which is $8-12\%$. Thus, the effect of boundary uncertainty is more pronounced in the Nusselt number compared to the temperature and velocities. The stochastic standard deviation is quite similar for all the four cases and thus, only the case with 32 strips is plotted here (\cref{Fig:Case B	Std}). It can be seen it is orders of magnitude smaller than the shift in the mean.
\begin{figure}[H]
	\centering
	\begin{subfigure}[t]{0.49\textwidth}
		\includegraphics[width=\textwidth]{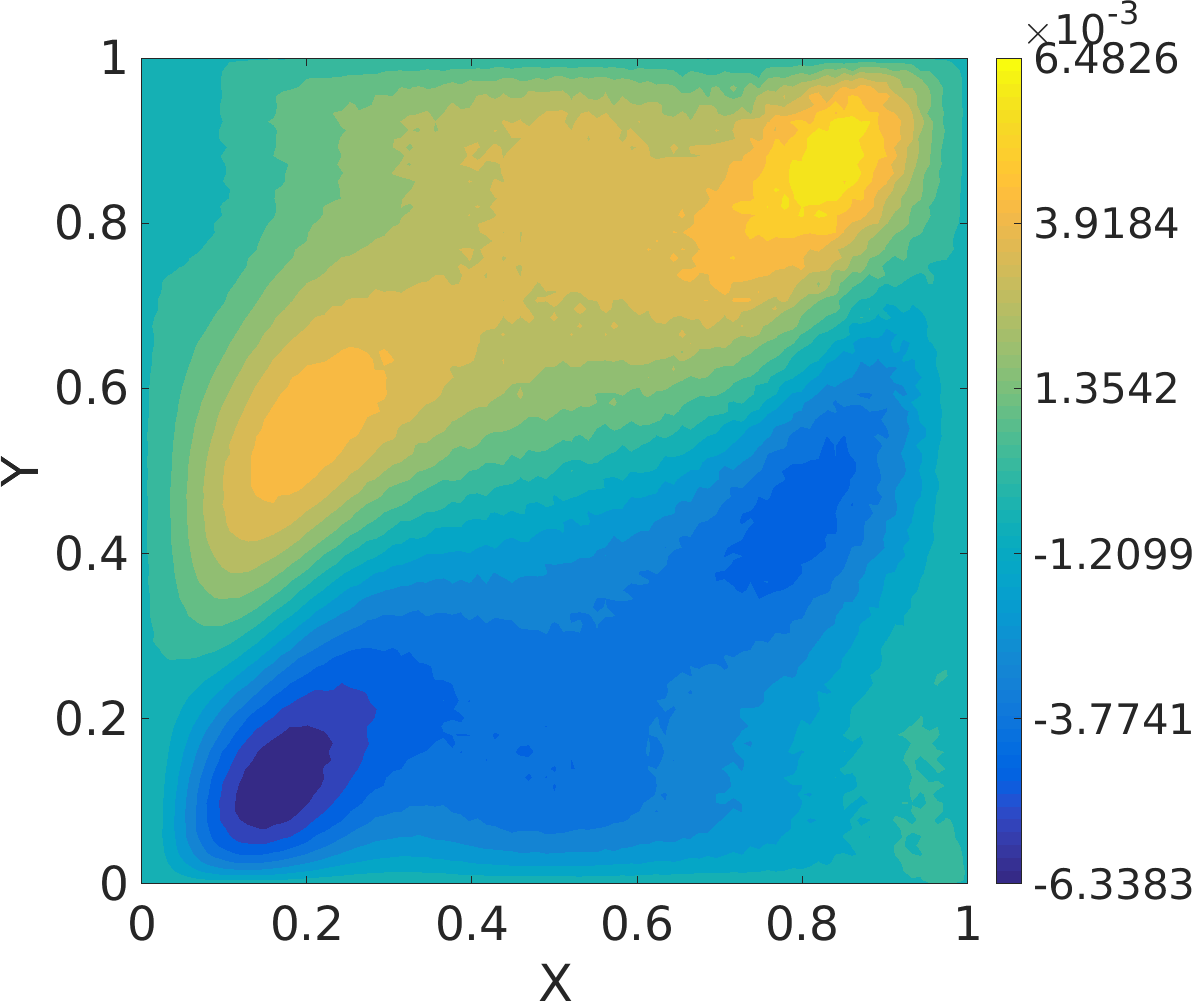}
		\caption{4 Strips}
		\label{Fig:Case B u Diff 4 Strips}
	\end{subfigure}
	\begin{subfigure}[t]{0.49\textwidth}
		\includegraphics[width=\textwidth]{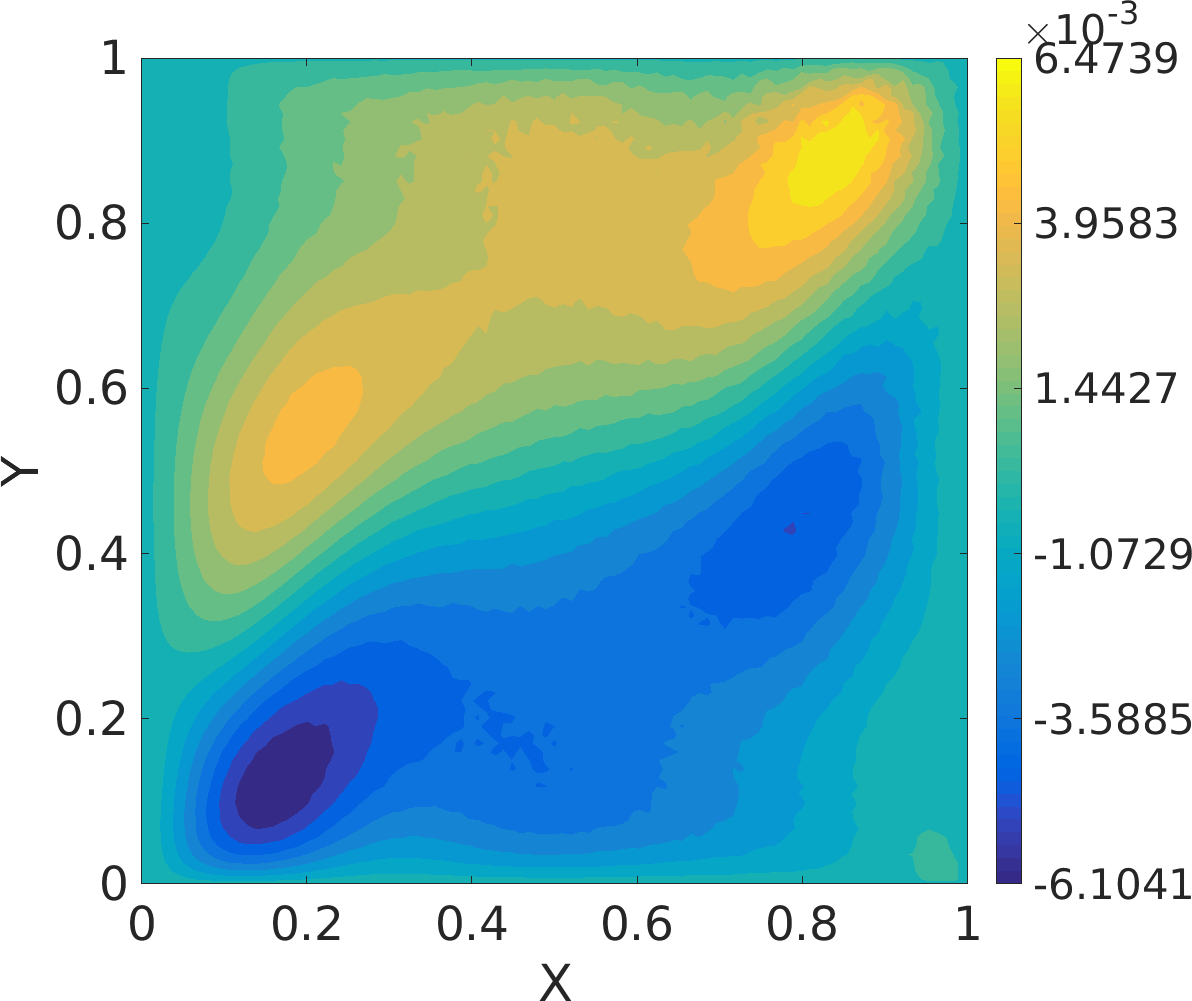}
		\caption{8 Strips}
		\label{Fig:Case B u Diff 8 Strips}
	\end{subfigure}
\end{figure}
\begin{figure}[H]\ContinuedFloat
	\begin{subfigure}[t]{0.49\textwidth}
		\includegraphics[width=\textwidth]{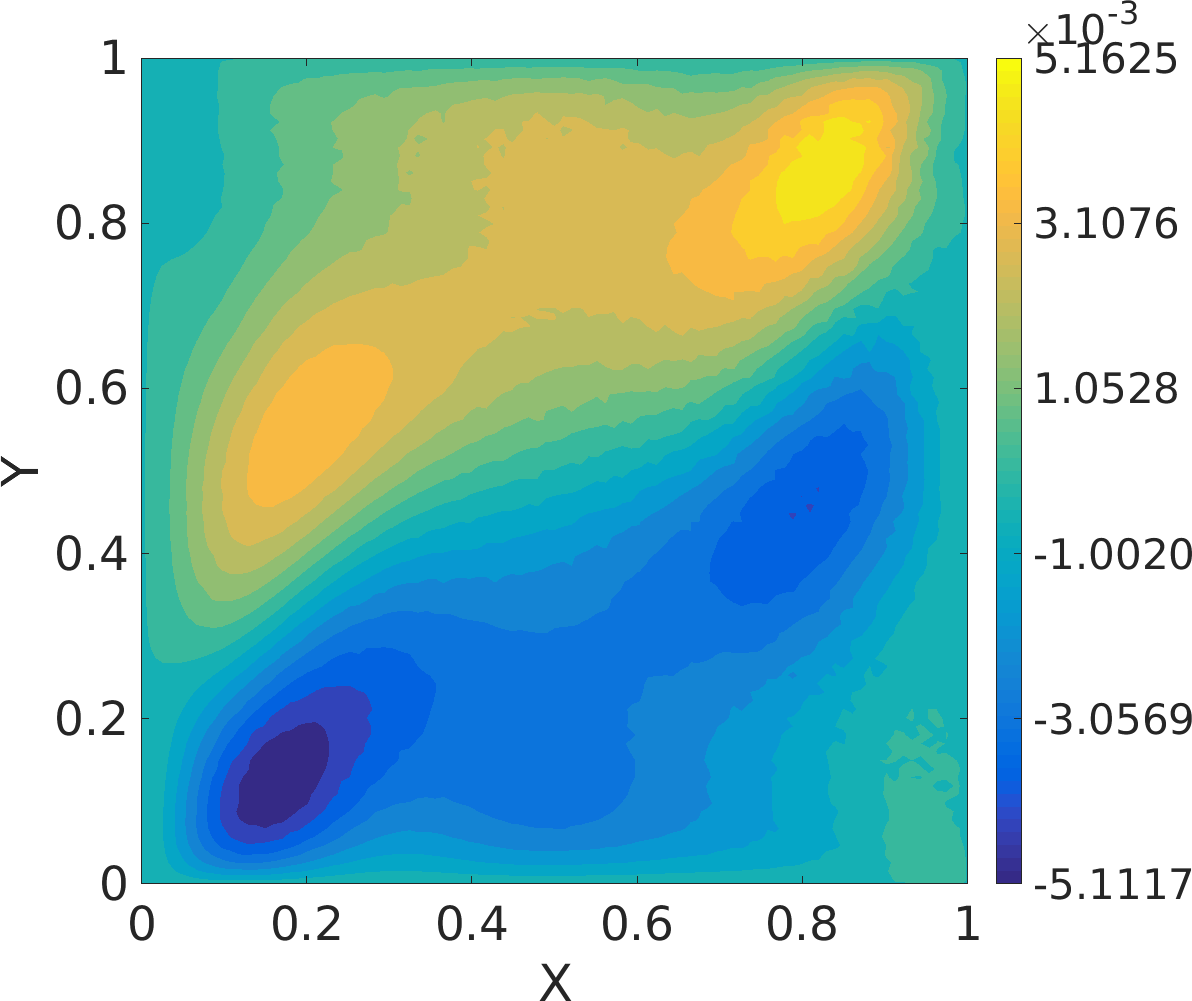}
		\caption{16 Strips}
		\label{Fig:Case B u Diff 16 Strips}
	\end{subfigure}
	\begin{subfigure}[t]{0.49\textwidth}
		\includegraphics[width=\textwidth]{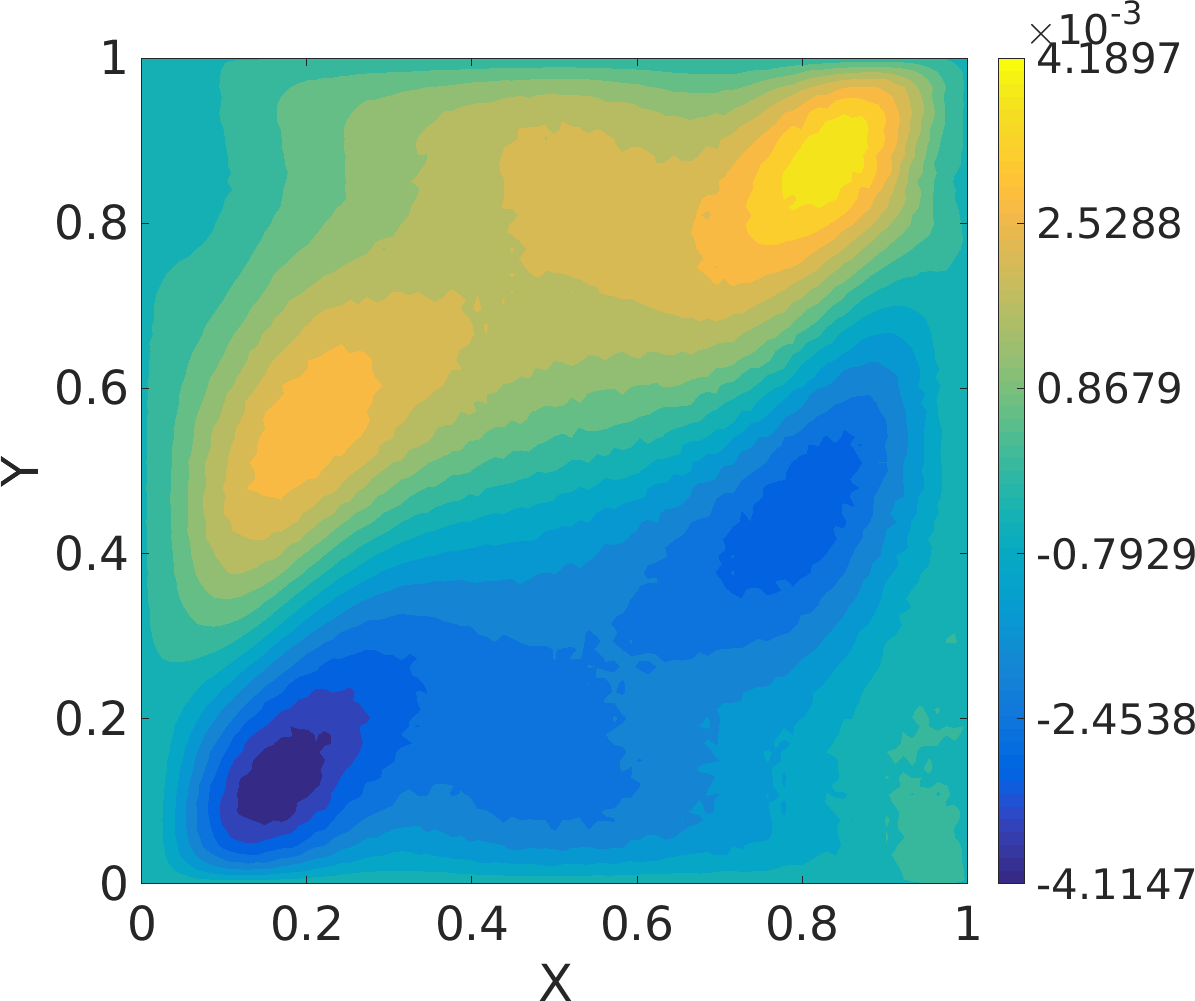}
		\caption{32 Strips}
		\label{Fig:Case B u Diff 32 Strips}
	\end{subfigure}	
	\caption{$Ra=10^6$ $Z$ Midplane X-Velocity: Difference between Stochastic Mean and Deterministic Value}
	\label{Fig:Case B u Diff}
\end{figure}
\begin{figure}[H]
	\centering
	\begin{subfigure}[t]{0.49\textwidth}
		\includegraphics[width=\textwidth]{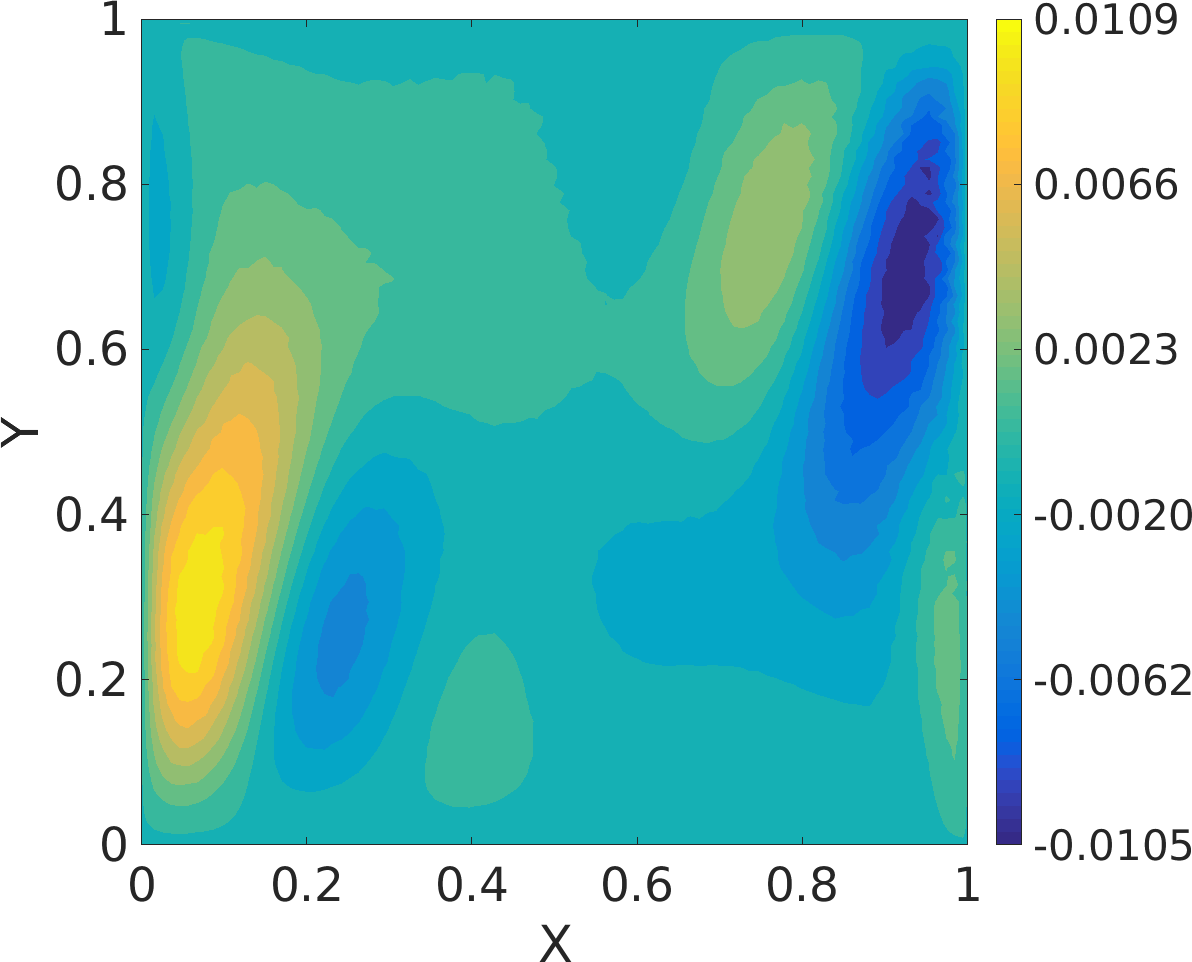}
		\caption{4 Strips}
		\label{Fig:Case B v Diff 4 Strips}
	\end{subfigure}
	\begin{subfigure}[t]{0.49\textwidth}
		\includegraphics[width=\textwidth]{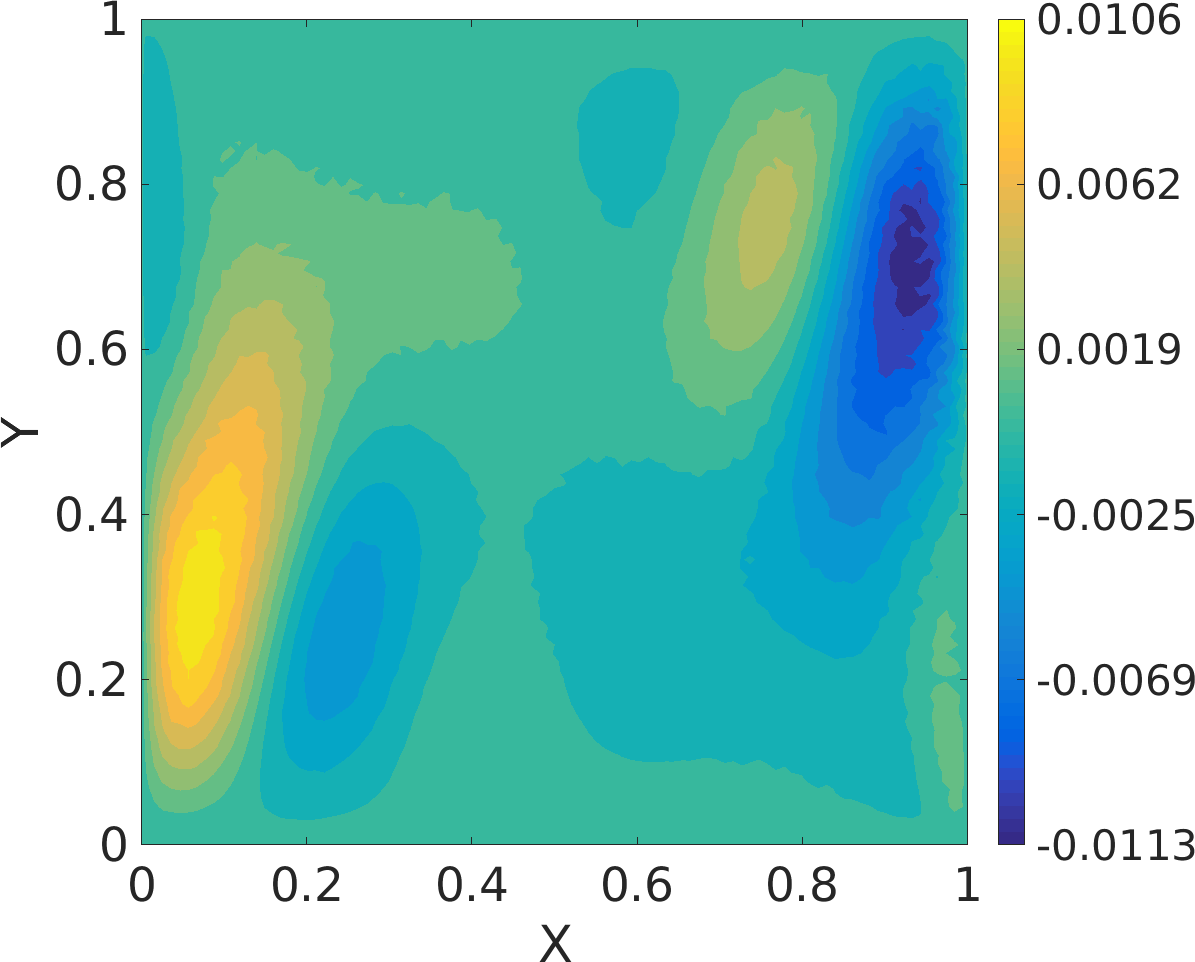}
		\caption{8 Strips}
		\label{Fig:Case B v Diff 8 Strips}
	\end{subfigure}
\end{figure}
\begin{figure}[H]\ContinuedFloat
	\begin{subfigure}[t]{0.49\textwidth}
		\includegraphics[width=\textwidth]{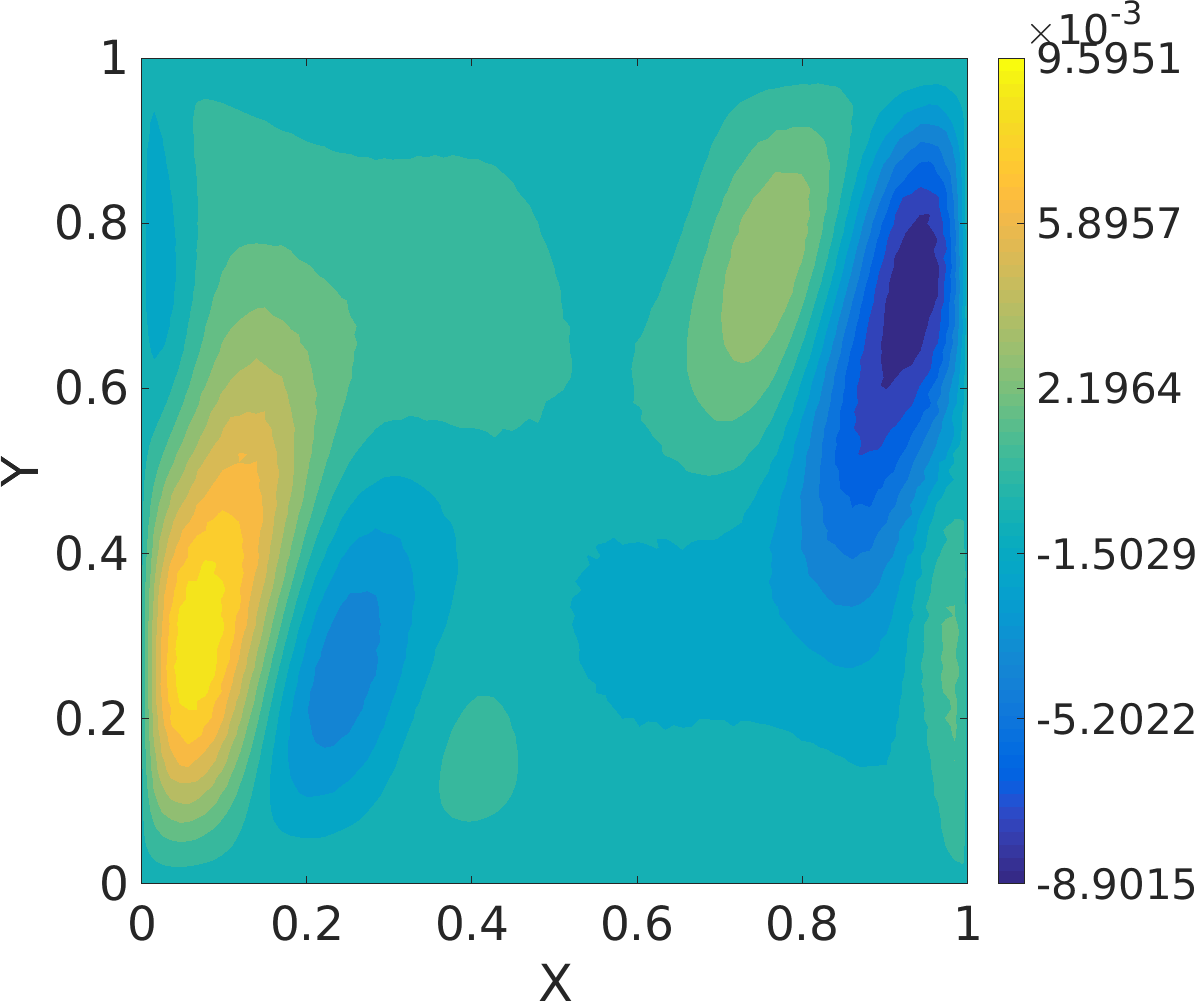}
		\caption{16 Strips}
		\label{Fig:Case B v Diff 16 Strips}
	\end{subfigure}
	\begin{subfigure}[t]{0.49\textwidth}
		\includegraphics[width=\textwidth]{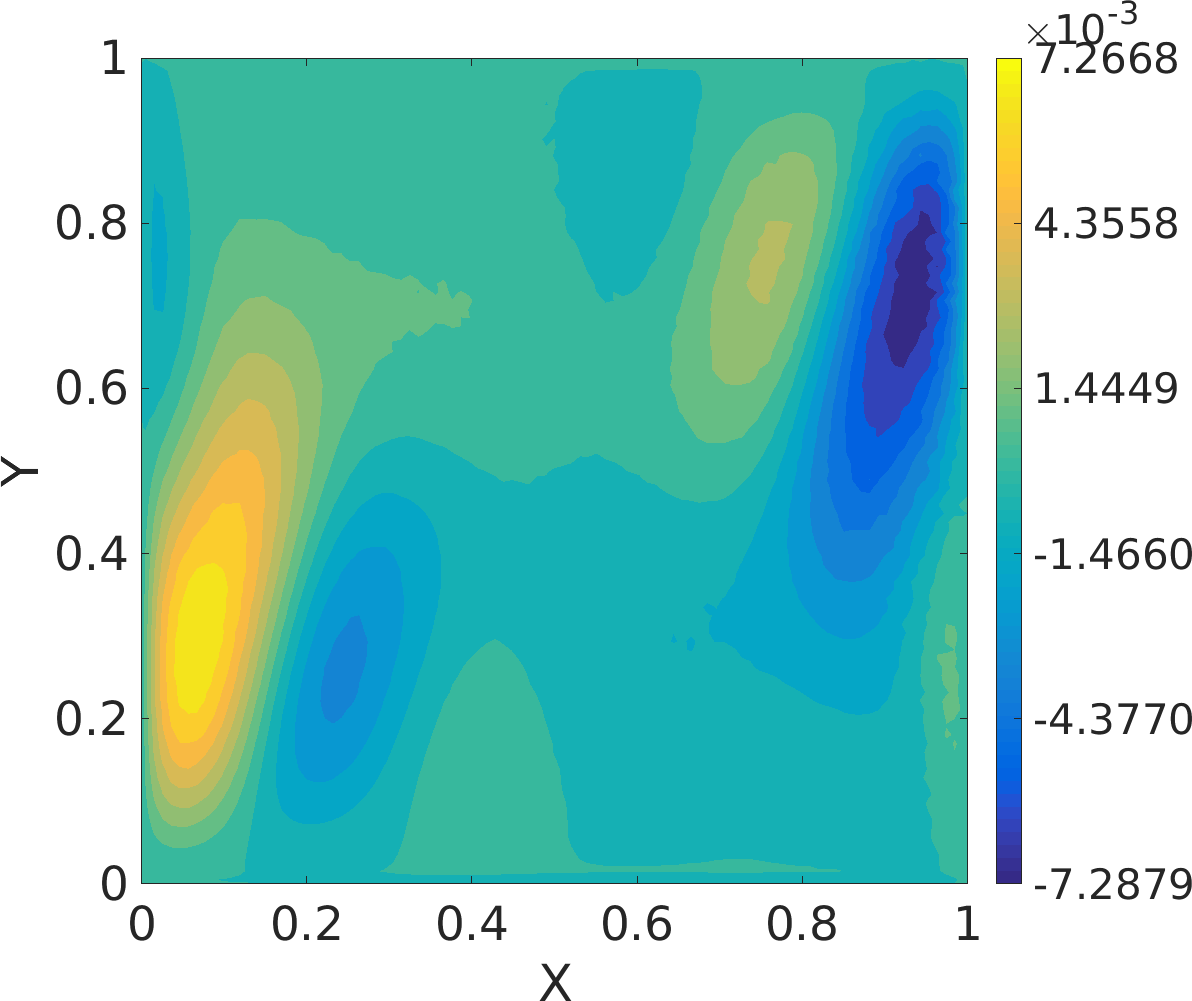}
		\caption{32 Strips}
		\label{Fig:Case B v Diff 32 Strips}
	\end{subfigure}	
	\caption{$Ra=10^6$ $Z$ Midplane Y-Velocity: Difference between Stochastic Mean and Deterministic Value}
	\label{Fig:Case B v Diff}
\end{figure}
\begin{figure}[H]
	\centering
	\begin{subfigure}[t]{0.32\textwidth}
		\includegraphics[width=\textwidth]{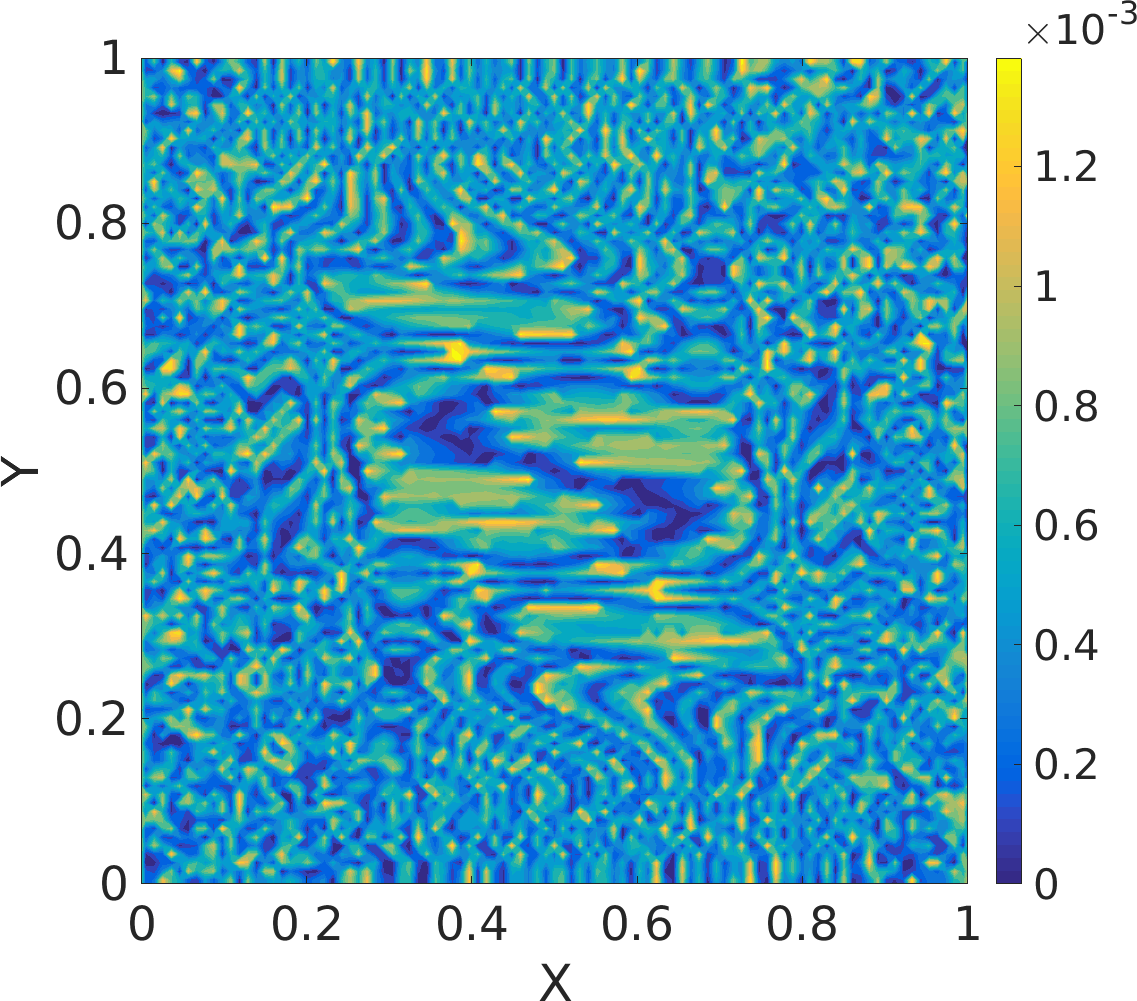}
		\caption{Temperature}
	\end{subfigure}
	\begin{subfigure}[t]{0.32\textwidth}
		\includegraphics[width=\textwidth]{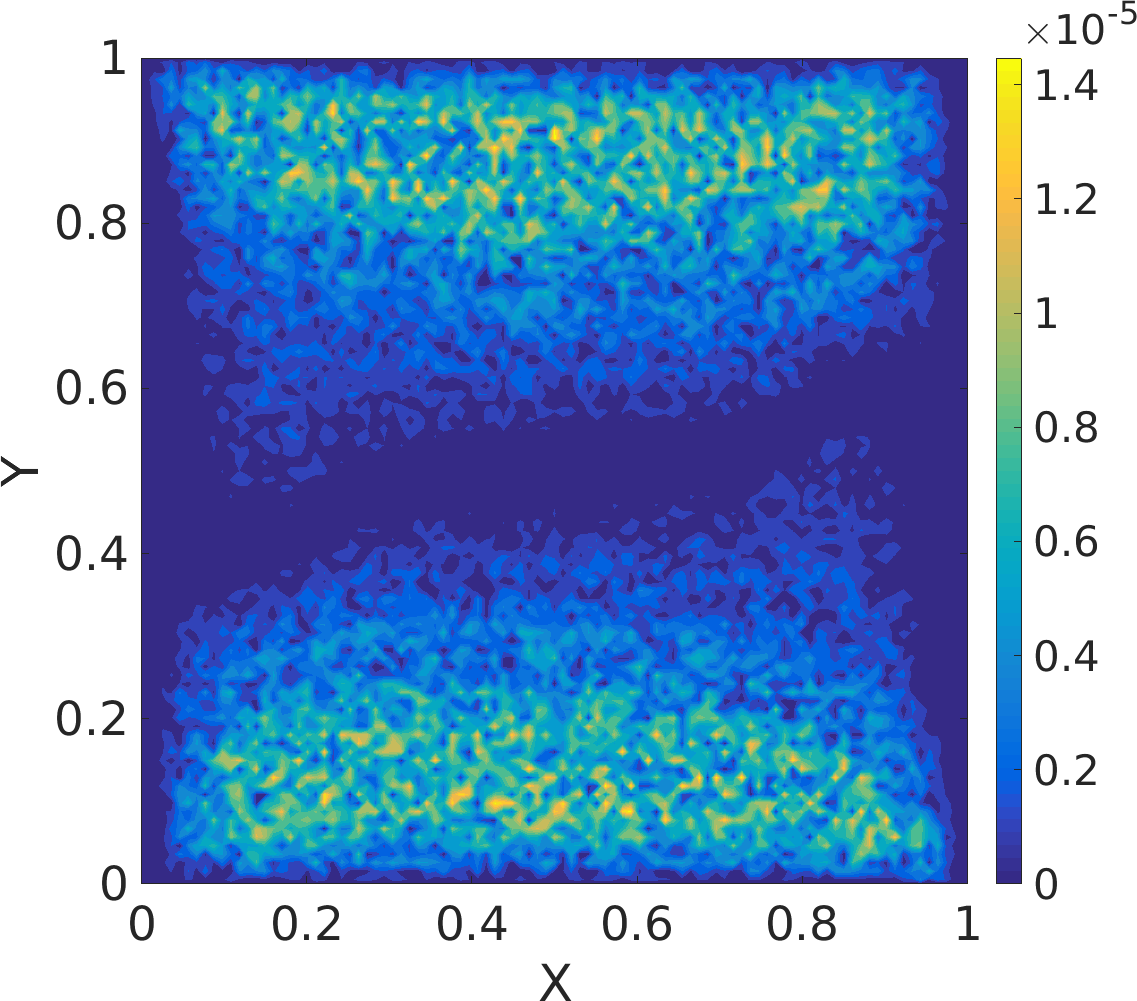}
		\caption{X-Velocity}
	\end{subfigure}
	\begin{subfigure}[t]{0.32\textwidth}
		\includegraphics[width=\textwidth]{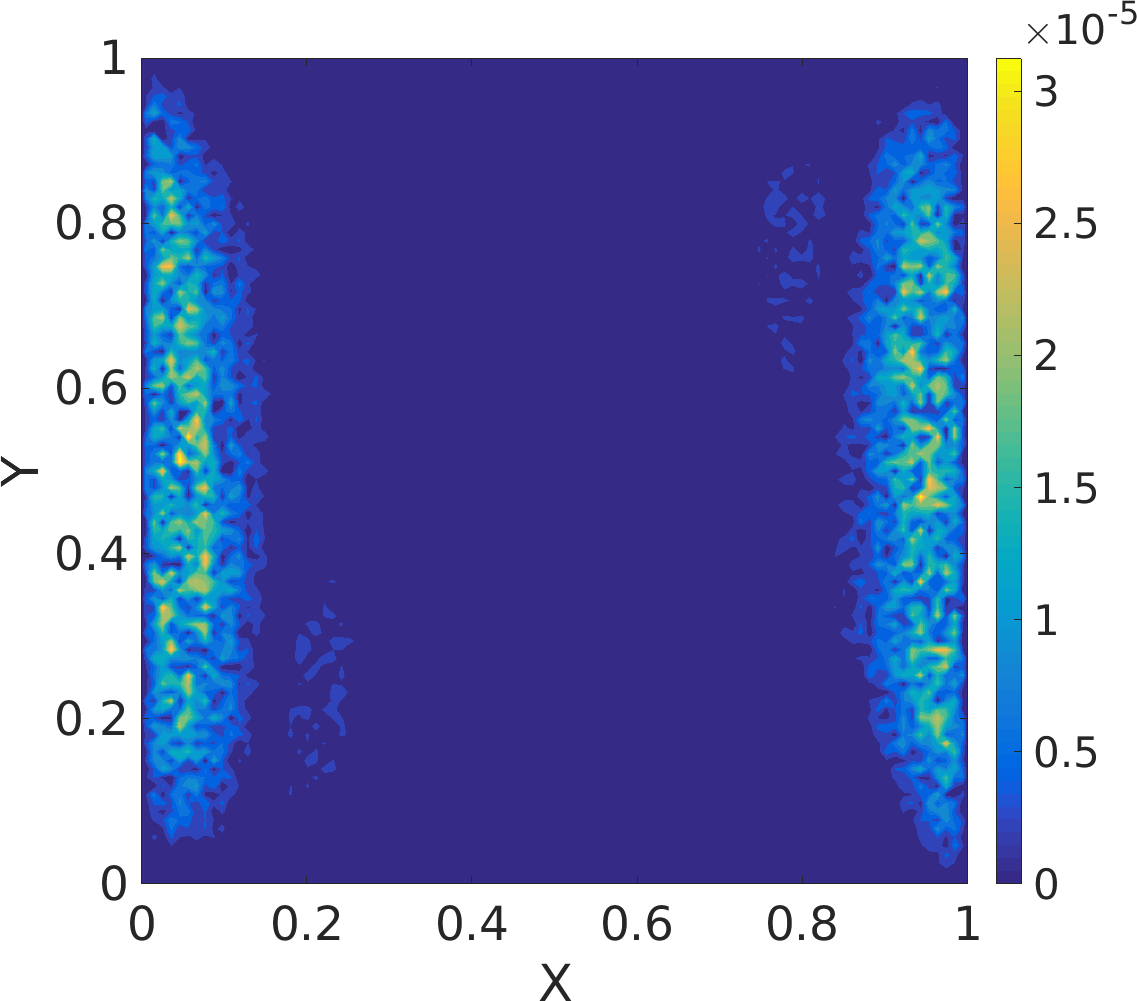}
		\caption{Y-Velocity}
	\end{subfigure}
	\caption{Stochastic Standard Deviation at $Z$ Mid-plane $Ra=10^6$ (32 Strips)}
	\label{Fig:Case B Std}
\end{figure}
%\subsubsection{Sensitivity Analysis}
\section{Conclusions}
This paper presents the input uncertainty propagation results for a three dimensional natural convection problem in a differentially heated cubical enclosure with two different cases. Case A assumes that the boundary conditions are uniform. Thus, uncertainty is introduced in the two non-dimensional parameters (Rayleigh and Prandtl numbers). For this case, the polynomial chaos expansion (PCE) method is used as a surrogate model with stochastic collocation to estimate the PCE coefficients. Case B deals with non-uniform stochastic boundary condition with deterministic material properties of the fluid. Since the temperature difference between the opposite walls drives the natural convection flow, the cold wall is held at a constant temperature and uncertainty is introduced in the hot wall temperature. A deep neural network based surrogate model is used here for estimating the output statistics with Monte Carlo method.
\par It is observed that the mean value of an output parameter averaged over the stochastic variation of the input can be different from the deterministic output value. The difference normalized by the deterministic value is defined as the relative shift of mean. \Cref{Tab:Shift of Mean Relative to Deterministic Values} summarizes the shift of mean relative to deterministic values of each of the output. Mean Rayleigh numbers of $10^5$ and $10^6$ are considered in case A whereas, case B deals only with a mean of $10^6$. The relative values of input standard deviation with respect to mean is 2\% and 3.33\% for case A and case B respectively. It can be seen that for case A, the stochastic effect is negligible in the Nusselt number and temperature whereas, it is of similar order as input in the velocities. On the other hand, for case B, the stochastic effect is much higher in Nusselt number. For both the cases, the standard deviation is much lower than the shift of mean. In general, it is concluded that the effect of uncertainties in the boundary conditions are much higher than the uncertainties in the Rayleigh and Prandtl numbers with uniform boundary conditions.
\begin{table}[H]
	\centering
	\begin{tabular}{|c|c|c|c|c|}
		\hline
		& \textbf{Input} & \textbf{Nusselt No.} & \textbf{Temperature} & \textbf{Velocities} \\ \hline
		Case A $Ra=10^5$ & $\sigma=2\% \mu$ & 0.3 \% & 0.1 \% & 3 \% \\ \hline
		Case A $Ra=10^6$ & $\sigma=2\% \mu$ & 0.3 \% & 0.3 \% & 1.5 \% \\ \hline
		Case B $Ra=10^6$ & $\sigma=3.33\% \mu$ & 12 \% & 4 \% & 5 \% \\ \hline
	\end{tabular}
	\caption{Shift of Mean Relative to Deterministic Values}
	\label{Tab:Shift of Mean Relative to Deterministic Values}
\end{table}
\par This paper also demonstrates the use of deep neural network for uncertainty quantification in natural convection problem for a high dimensional stochastic problem. The novel approach of dividing the boundary surface into domains and treating each domain value as a stochastic input is demonstrated here for the first time (to the best of our knowledge). This approach tries to mimic the experiments, when multiple feedback control systems are used to impose the boundary conditions. Note that although in this work, it is assumed that the domain values are independent and identically distributed variables, this assumption can be relaxed during sampling if additional data of the feedback control system is available. The approach can also be used for other class of problems like optimization and inverse problems.
\section*{Acknowledgments}
This work was funded by the Digital Manufacturing and Design Innovation Institute with support in part by the U.S. Department of the Army. Any opinions, findings, and conclusions or recommendations expressed in this material are those of the author(s) and do not necessarily reflect the views of the Department of the Army. %The authors also thank Steve Udvardy and Beau Glim of North American Die Casting Association (NADCA), Chicago White Metal Casting, Mercury Castings, RCM Industries, Twin City Die Castings Company and Visi-Trak Worldwide for providing die cast parts used for calibration and validation. Technical discussions with Alex Monroe of Mercury Castings and his suggestions are appreciated.
\section*{References}
\bibliography{References}

\begin{thebibliography}{50}
\providecommand{\natexlab}[1]{#1}
\providecommand{\url}[1]{\texttt{#1}}
\providecommand{\urlprefix}{URL }
\expandafter\ifx\csname urlstyle\endcsname\relax
  \providecommand{\doi}[1]{doi:\discretionary{}{}{}#1}\else
  \providecommand{\doi}[1]{doi:\discretionary{}{}{}\begingroup
  \urlstyle{rm}\url{#1}\endgroup}\fi
\providecommand{\bibinfo}[2]{#2}

\bibitem[{De~Vahl~Davis(1983)}]{de1983natural}
\bibinfo{author}{G.~De~Vahl~Davis}, \bibinfo{title}{Natural convection of air
  in a square cavity: a bench mark numerical solution},
  \bibinfo{journal}{International Journal for Numerical Methods in Fluids}
  \bibinfo{volume}{3}~(\bibinfo{number}{3}) (\bibinfo{year}{1983})
  \bibinfo{pages}{249--264}.

\bibitem[{Fusegi et~al.(1991)Fusegi, Hyun, Kuwahara, and
  Farouk}]{fusegi1991numerical}
\bibinfo{author}{T.~Fusegi}, \bibinfo{author}{J.~Hyun},
  \bibinfo{author}{K.~Kuwahara}, \bibinfo{author}{B.~Farouk}, \bibinfo{title}{A
  numerical study of three-dimensional natural convection in a differentially
  heated cubical enclosure}, \bibinfo{journal}{{International Journal of Heat
  and Mass Transfer}} \bibinfo{volume}{34}~(\bibinfo{number}{6})
  (\bibinfo{year}{1991}) \bibinfo{pages}{1543--1557}.

\bibitem[{Le~Qu{\'e}r{\'e}(1991)}]{le1991accurate}
\bibinfo{author}{P.~Le~Qu{\'e}r{\'e}}, \bibinfo{title}{Accurate solutions to
  the square thermally driven cavity at high Rayleigh number},
  \bibinfo{journal}{Computers \& Fluids}
  \bibinfo{volume}{20}~(\bibinfo{number}{1}) (\bibinfo{year}{1991})
  \bibinfo{pages}{29--41}.

\bibitem[{Hu et~al.(2017)Hu, Li, Wu, Li, and Li}]{hu2017flow}
\bibinfo{author}{Y.~Hu}, \bibinfo{author}{Y.~Li}, \bibinfo{author}{C.~Wu},
  \bibinfo{author}{S.~Li}, \bibinfo{author}{M.~Li}, \bibinfo{title}{Flow
  pattern and heat transfer in Rayleigh-B{\'e}nard convection of cold water
  near its density maximum in a rectangular cavity},
  \bibinfo{journal}{International Journal of Heat and Mass Transfer}
  \bibinfo{volume}{107} (\bibinfo{year}{2017}) \bibinfo{pages}{1065--1075}.

\bibitem[{Li et~al.(2016)Li, Zhang, Zhang, and Wu}]{li2016three}
\bibinfo{author}{Y.~Li}, \bibinfo{author}{H.~Zhang},
  \bibinfo{author}{L.~Zhang}, \bibinfo{author}{C.~Wu},
  \bibinfo{title}{Three-dimensional numerical simulation of double-diffusive
  Rayleigh--B{\'e}nard convection in a cylindrical enclosure of aspect ratio
  2}, \bibinfo{journal}{International Journal of Heat and Mass Transfer}
  \bibinfo{volume}{98} (\bibinfo{year}{2016}) \bibinfo{pages}{472--483}.

\bibitem[{Shu et~al.(2003)Shu, Ding, and Yeo}]{shu2003local}
\bibinfo{author}{C.~Shu}, \bibinfo{author}{H.~Ding}, \bibinfo{author}{K.~Yeo},
  \bibinfo{title}{Local radial basis function-based differential quadrature
  method and its application to solve two-dimensional incompressible
  Navier--Stokes equations}, \bibinfo{journal}{Computer Methods in Applied
  Mechanics and Engineering} \bibinfo{volume}{192}~(\bibinfo{number}{7-8})
  (\bibinfo{year}{2003}) \bibinfo{pages}{941--954}.

\bibitem[{Yigit et~al.(2015)Yigit, Poole, and Chakraborty}]{yigit2015effects}
\bibinfo{author}{S.~Yigit}, \bibinfo{author}{R.~Poole},
  \bibinfo{author}{N.~Chakraborty}, \bibinfo{title}{Effects of aspect ratio on
  laminar Rayleigh--B{\'e}nard convection of power-law fluids in rectangular
  enclosures: A numerical investigation}, \bibinfo{journal}{International
  Journal of Heat and Mass Transfer} \bibinfo{volume}{91}
  (\bibinfo{year}{2015}) \bibinfo{pages}{1292--1307}.

\bibitem[{Icoz and Jaluria(2004)}]{icoz2004design}
\bibinfo{author}{T.~Icoz}, \bibinfo{author}{Y.~Jaluria}, \bibinfo{title}{Design
  of cooling systems for electronic equipment using both experimental and
  numerical inputs}, \bibinfo{journal}{Journal of Electronic Packaging}
  \bibinfo{volume}{126}~(\bibinfo{number}{4}) (\bibinfo{year}{2004})
  \bibinfo{pages}{465--471}.

\bibitem[{Sharma et~al.(2007)Sharma, Velusamy, Balaji, and
  Venkateshan}]{sharma2007conjugate}
\bibinfo{author}{A.~Sharma}, \bibinfo{author}{K.~Velusamy},
  \bibinfo{author}{C.~Balaji}, \bibinfo{author}{S.~Venkateshan},
  \bibinfo{title}{Conjugate turbulent natural convection with surface radiation
  in air filled rectangular enclosures}, \bibinfo{journal}{International
  Journal of Heat and Mass Transfer}
  \bibinfo{volume}{50}~(\bibinfo{number}{3-4}) (\bibinfo{year}{2007})
  \bibinfo{pages}{625--639}.

\bibitem[{Herrlin and Belady(2006)}]{herrlin2006gravity}
\bibinfo{author}{M.~Herrlin}, \bibinfo{author}{C.~Belady},
  \bibinfo{title}{Gravity-assisted air mixing in data centers and how it
  affects the rack cooling effectiveness}, in: \bibinfo{booktitle}{Thermal and
  Thermomechanical Phenomena in Electronics Systems, 2006. ITHERM'06. The Tenth
  Intersociety Conference on}, \bibinfo{organization}{IEEE},
  \bibinfo{pages}{5--pp}, \bibinfo{year}{2006}.

\bibitem[{Patankar(1980)}]{patankar1980numerical}
\bibinfo{author}{S.~Patankar}, \bibinfo{title}{Numerical heat transfer and
  fluid flow}, \bibinfo{publisher}{CRC Press}, \bibinfo{year}{1980}.

\bibitem[{Xiu and Karniadakis(2002{\natexlab{a}})}]{xiu2002modeling}
\bibinfo{author}{D.~Xiu}, \bibinfo{author}{G.~Karniadakis},
  \bibinfo{title}{Modeling uncertainty in steady state diffusion problems via
  generalized polynomial chaos}, \bibinfo{journal}{Computer Methods in Applied
  Mechanics and Engineering} \bibinfo{volume}{191}~(\bibinfo{number}{43})
  (\bibinfo{year}{2002}{\natexlab{a}}) \bibinfo{pages}{4927--4948}.

\bibitem[{Xiu and Karniadakis(2003)}]{xiu2003modeling}
\bibinfo{author}{D.~Xiu}, \bibinfo{author}{G.~Karniadakis},
  \bibinfo{title}{Modeling uncertainty in flow simulations via generalized
  polynomial chaos}, \bibinfo{journal}{Journal of Computational Physics}
  \bibinfo{volume}{187}~(\bibinfo{number}{1}) (\bibinfo{year}{2003})
  \bibinfo{pages}{137--167}.

\bibitem[{Maitre et~al.(2001)Maitre, Knio, Najm, and
  Ghanem}]{knio2001stochastic}
\bibinfo{author}{O.~Maitre}, \bibinfo{author}{O.~Knio},
  \bibinfo{author}{H.~Najm}, \bibinfo{author}{R.~Ghanem}, \bibinfo{title}{A
  stochastic projection method for fluid flow: I. basic formulation},
  \bibinfo{journal}{Journal of Computational Physics}
  \bibinfo{volume}{173}~(\bibinfo{number}{2}) (\bibinfo{year}{2001})
  \bibinfo{pages}{481--511}.

\bibitem[{Maitre et~al.(2004)Maitre, Reagan, Debusschere, Najm, Ghanem, and
  Knio}]{le2004natural}
\bibinfo{author}{O.~Maitre}, \bibinfo{author}{M.~Reagan},
  \bibinfo{author}{B.~Debusschere}, \bibinfo{author}{H.~Najm},
  \bibinfo{author}{R.~Ghanem}, \bibinfo{author}{O.~Knio},
  \bibinfo{title}{Natural convection in a closed cavity under stochastic
  non-Boussinesq conditions}, \bibinfo{journal}{SIAM Journal on Scientific
  Computing} \bibinfo{volume}{26}~(\bibinfo{number}{2}) (\bibinfo{year}{2004})
  \bibinfo{pages}{375--394}.

\bibitem[{Venturi et~al.(2012)Venturi, Choi, and
  Karniadakis}]{venturi2012supercritical}
\bibinfo{author}{D.~Venturi}, \bibinfo{author}{M.~Choi},
  \bibinfo{author}{G.~Karniadakis}, \bibinfo{title}{Supercritical
  quasi-conduction states in stochastic Rayleigh--B{\'e}nard convection},
  \bibinfo{journal}{International Journal of Heat and Mass Transfer}
  \bibinfo{volume}{55}~(\bibinfo{number}{13-14}) (\bibinfo{year}{2012})
  \bibinfo{pages}{3732--3743}.

\bibitem[{Ganapathysubramanian and
  Zabaras(2007)}]{ganapathysubramanian2007sparse}
\bibinfo{author}{B.~Ganapathysubramanian}, \bibinfo{author}{N.~Zabaras},
  \bibinfo{title}{Sparse grid collocation schemes for stochastic natural
  convection problems}, \bibinfo{journal}{Journal of Computational Physics}
  \bibinfo{volume}{225}~(\bibinfo{number}{1}) (\bibinfo{year}{2007})
  \bibinfo{pages}{652--685}.

\bibitem[{Carnevale et~al.(2013)Carnevale, Montomoli, D’Ammaro, Salvadori,
  and Martelli}]{carnevale2013uncertainty}
\bibinfo{author}{M.~Carnevale}, \bibinfo{author}{F.~Montomoli},
  \bibinfo{author}{A.~D’Ammaro}, \bibinfo{author}{S.~Salvadori},
  \bibinfo{author}{F.~Martelli}, \bibinfo{title}{{Uncertainty quantification: A
  stochastic method for heat transfer prediction using LES}},
  \bibinfo{journal}{Journal of Turbomachinery}
  \bibinfo{volume}{135}~(\bibinfo{number}{5}) (\bibinfo{year}{2013})
  \bibinfo{pages}{051021}.

\bibitem[{Marepalli et~al.(2014)Marepalli, Murthy, Qiu, and
  Ruan}]{marepalli2014quantifying}
\bibinfo{author}{P.~Marepalli}, \bibinfo{author}{J.~Murthy},
  \bibinfo{author}{B.~Qiu}, \bibinfo{author}{X.~Ruan},
  \bibinfo{title}{Quantifying uncertainty in multiscale heat conduction
  calculations}, \bibinfo{journal}{Journal of Heat Transfer}
  \bibinfo{volume}{136}~(\bibinfo{number}{11}) (\bibinfo{year}{2014})
  \bibinfo{pages}{111301}.

\bibitem[{Fajraoui et~al.(2017)Fajraoui, Fahs, Younes, and
  Sudret}]{fajraoui2017analyzing}
\bibinfo{author}{N.~Fajraoui}, \bibinfo{author}{M.~Fahs},
  \bibinfo{author}{A.~Younes}, \bibinfo{author}{B.~Sudret},
  \bibinfo{title}{Analyzing natural convection in porous enclosure with
  polynomial chaos expansions: Effect of thermal dispersion, anisotropic
  permeability and heterogeneity}, \bibinfo{journal}{International Journal of
  Heat and Mass Transfer} \bibinfo{volume}{115} (\bibinfo{year}{2017})
  \bibinfo{pages}{205--224}.

\bibitem[{Fezi and Krane(2017)}]{fezi2017uncertainty}
\bibinfo{author}{K.~Fezi}, \bibinfo{author}{M.~Krane},
  \bibinfo{title}{Uncertainty quantification in modeling metal alloy
  solidification}, \bibinfo{journal}{Journal of Heat Transfer}
  \bibinfo{volume}{139}~(\bibinfo{number}{8}) (\bibinfo{year}{2017})
  \bibinfo{pages}{082301}.

\bibitem[{Hosder et~al.(2006)Hosder, Walters, and Perez}]{hosder2006non}
\bibinfo{author}{S.~Hosder}, \bibinfo{author}{R.~Walters},
  \bibinfo{author}{R.~Perez}, \bibinfo{title}{{A non-intrusive polynomial chaos
  method for uncertainty propagation in CFD simulations}}, in:
  \bibinfo{booktitle}{44th AIAA Aerospace Sciences Meeting and Exhibit},
  \bibinfo{pages}{891}, \bibinfo{year}{2006}.

\bibitem[{Kumar et~al.(2016)Kumar, Raisee, and Lacor}]{kumar2016efficient}
\bibinfo{author}{D.~Kumar}, \bibinfo{author}{M.~Raisee},
  \bibinfo{author}{C.~Lacor}, \bibinfo{title}{{An efficient non-intrusive
  reduced basis model for high dimensional stochastic problems in CFD}},
  \bibinfo{journal}{Computers \& Fluids} \bibinfo{volume}{138}
  (\bibinfo{year}{2016}) \bibinfo{pages}{67--82}.

\bibitem[{Smolyak(1963)}]{smolyak1963quadrature}
\bibinfo{author}{S.~Smolyak}, \bibinfo{title}{Quadrature and interpolation
  formulas for tensor products of certain classes of functions}, in:
  \bibinfo{booktitle}{{Soviet Math. Dokl.}}, vol.~\bibinfo{volume}{4},
  \bibinfo{pages}{240--243}, \bibinfo{year}{1963}.

\bibitem[{Heiss and Winschel(2008)}]{heiss2008likelihood}
\bibinfo{author}{F.~Heiss}, \bibinfo{author}{V.~Winschel},
  \bibinfo{title}{Likelihood approximation by numerical integration on sparse
  grids}, \bibinfo{journal}{{Journal of Econometrics}}
  \bibinfo{volume}{144}~(\bibinfo{number}{1}) (\bibinfo{year}{2008})
  \bibinfo{pages}{62--80}.

\bibitem[{Caflisch(1998)}]{caflisch1998monte}
\bibinfo{author}{R.~Caflisch}, \bibinfo{title}{Monte carlo and quasi-monte
  carlo methods}, \bibinfo{journal}{Acta Numerica} \bibinfo{volume}{7}
  (\bibinfo{year}{1998}) \bibinfo{pages}{1--49}.

\bibitem[{Hornik et~al.(1989)Hornik, Stinchcombe, and
  White}]{hornik1989multilayer}
\bibinfo{author}{K.~Hornik}, \bibinfo{author}{M.~Stinchcombe},
  \bibinfo{author}{H.~White}, \bibinfo{title}{Multilayer feedforward networks
  are universal approximators}, \bibinfo{journal}{Neural Networks}
  \bibinfo{volume}{2}~(\bibinfo{number}{5}) (\bibinfo{year}{1989})
  \bibinfo{pages}{359--366}.

\bibitem[{Schmidhuber(2015)}]{schmidhuber2015deep}
\bibinfo{author}{J.~Schmidhuber}, \bibinfo{title}{Deep learning in neural
  networks: An overview}, \bibinfo{journal}{Neural Networks}
  \bibinfo{volume}{61} (\bibinfo{year}{2015}) \bibinfo{pages}{85--117}.

\bibitem[{Gholami et~al.(2015)Gholami, Bonakdari, Zaji, and
  Akhtari}]{gholami2015simulation}
\bibinfo{author}{A.~Gholami}, \bibinfo{author}{H.~Bonakdari},
  \bibinfo{author}{A.~Zaji}, \bibinfo{author}{A.~Akhtari},
  \bibinfo{title}{Simulation of open channel bend characteristics using
  computational fluid dynamics and artificial neural networks},
  \bibinfo{journal}{Engineering Applications of Computational Fluid Mechanics}
  \bibinfo{volume}{9}~(\bibinfo{number}{1}) (\bibinfo{year}{2015})
  \bibinfo{pages}{355--369}.

\bibitem[{Sablani et~al.(2005)Sablani, Kacimov, Perret, Mujumdar, and
  Campo}]{sablani2005non}
\bibinfo{author}{S.~Sablani}, \bibinfo{author}{A.~Kacimov},
  \bibinfo{author}{J.~Perret}, \bibinfo{author}{A.~Mujumdar},
  \bibinfo{author}{A.~Campo}, \bibinfo{title}{Non-iterative estimation of heat
  transfer coefficients using artificial neural network models},
  \bibinfo{journal}{International Journal of Heat and Mass Transfer}
  \bibinfo{volume}{48}~(\bibinfo{number}{3-4}) (\bibinfo{year}{2005})
  \bibinfo{pages}{665--679}.

\bibitem[{Santra et~al.(2009)Santra, Chakraborty, and
  Sen}]{santra2009prediction}
\bibinfo{author}{A.~Santra}, \bibinfo{author}{N.~Chakraborty},
  \bibinfo{author}{S.~Sen}, \bibinfo{title}{Prediction of heat transfer due to
  presence of copper--water nanofluid using resilient-propagation neural
  network}, \bibinfo{journal}{International Journal of Thermal Sciences}
  \bibinfo{volume}{48}~(\bibinfo{number}{7}) (\bibinfo{year}{2009})
  \bibinfo{pages}{1311--1318}.

\bibitem[{Cz{\'e}l et~al.(2014)Cz{\'e}l, Woodbury, and
  Gr{\'o}f}]{czel2014simultaneous}
\bibinfo{author}{B.~Cz{\'e}l}, \bibinfo{author}{K.~Woodbury},
  \bibinfo{author}{G.~Gr{\'o}f}, \bibinfo{title}{Simultaneous estimation of
  temperature-dependent volumetric heat capacity and thermal conductivity
  functions via neural networks}, \bibinfo{journal}{International Journal of
  Heat and Mass Transfer} \bibinfo{volume}{68} (\bibinfo{year}{2014})
  \bibinfo{pages}{1--13}.

\bibitem[{Zhang and Duraisamy(2015)}]{zhang2015machine}
\bibinfo{author}{Z.~Zhang}, \bibinfo{author}{K.~Duraisamy},
  \bibinfo{title}{Machine learning methods for data-driven turbulence
  modeling}, in: \bibinfo{booktitle}{22nd AIAA Computational Fluid Dynamics
  Conference}, \bibinfo{pages}{2460}, \bibinfo{year}{2015}.

\bibitem[{Tripathy and Bilionis(2018)}]{tripathy2018deep}
\bibinfo{author}{R.~Tripathy}, \bibinfo{author}{I.~Bilionis},
  \bibinfo{title}{{Deep UQ: Learning deep neural network surrogate models for
  high dimensional uncertainty quantification}}, \bibinfo{journal}{arXiv
  preprint arXiv:1802.00850} .

\bibitem[{Khoo et~al.(2017)Khoo, Lu, and Ying}]{khoo2017solving}
\bibinfo{author}{Y.~Khoo}, \bibinfo{author}{J.~Lu}, \bibinfo{author}{L.~Ying},
  \bibinfo{title}{Solving parametric PDE problems with artificial neural
  networks}, \bibinfo{journal}{arXiv preprint arXiv:1707.03351} .

\bibitem[{Zhang et~al.(2018)Zhang, Lu, Guo, and
  Karniadakis}]{zhang2018quantifying}
\bibinfo{author}{D.~Zhang}, \bibinfo{author}{L.~Lu}, \bibinfo{author}{L.~Guo},
  \bibinfo{author}{G.~Karniadakis}, \bibinfo{title}{Quantifying total
  uncertainty in physics-informed neural networks for solving forward and
  inverse stochastic problems}, \bibinfo{journal}{arXiv preprint
  arXiv:1809.08327} .

\bibitem[{Shahane et~al.(2018)Shahane, Aluru, Ferreira, Kapoor, and
  Vanka}]{shahane2018finite}
\bibinfo{author}{S.~Shahane}, \bibinfo{author}{N.~Aluru},
  \bibinfo{author}{P.~Ferreira}, \bibinfo{author}{S.~Kapoor},
  \bibinfo{author}{S.~Vanka}, \bibinfo{title}{Finite Volume Simulation
  Framework for Die Casting with Uncertainty Quantification},
  \bibinfo{journal}{arXiv preprint arXiv:1810.08572} .

\bibitem[{Harlow and Welch(1965)}]{harlow1965numerical}
\bibinfo{author}{F.~Harlow}, \bibinfo{author}{J.~Welch},
  \bibinfo{title}{Numerical calculation of time-dependent viscous
  incompressible flow of fluid with free surface}, \bibinfo{journal}{The
  physics of fluids} \bibinfo{volume}{8}~(\bibinfo{number}{12})
  (\bibinfo{year}{1965}) \bibinfo{pages}{2182--2189}.

\bibitem[{Xiu and Karniadakis(2002{\natexlab{b}})}]{xiu2002wiener}
\bibinfo{author}{D.~Xiu}, \bibinfo{author}{G.~Karniadakis},
  \bibinfo{title}{{The Wiener--Askey polynomial chaos for stochastic
  differential equations}}, \bibinfo{journal}{{SIAM Journal on Scientific
  Computing}} \bibinfo{volume}{24}~(\bibinfo{number}{2})
  (\bibinfo{year}{2002}{\natexlab{b}}) \bibinfo{pages}{619--644}.

\bibitem[{Wiener(1938)}]{wiener1938homogeneous}
\bibinfo{author}{N.~Wiener}, \bibinfo{title}{The homogeneous chaos},
  \bibinfo{journal}{{American Journal of Mathematics}}
  \bibinfo{volume}{60}~(\bibinfo{number}{4}).

\bibitem[{Smith(2013)}]{smith2013uncertainty}
\bibinfo{author}{R.~Smith}, \bibinfo{title}{Uncertainty quantification: theory,
  implementation, and applications}, \bibinfo{year}{2013}.

\bibitem[{Marelli and Sudret(2014)}]{marelli2014uqlab}
\bibinfo{author}{S.~Marelli}, \bibinfo{author}{B.~Sudret},
  \bibinfo{title}{{UQLab: A framework for uncertainty quantification in
  Matlab}}, in: \bibinfo{booktitle}{Vulnerability, Uncertainty, and Risk:
  Quantification, Mitigation, and Management}, \bibinfo{pages}{2554--2563},
  \bibinfo{year}{2014}.

\bibitem[{Goodfellow et~al.(2016)Goodfellow, Bengio, and
  Courville}]{goodfellow2016deep}
\bibinfo{author}{I.~Goodfellow}, \bibinfo{author}{Y.~Bengio},
  \bibinfo{author}{A.~Courville}, \bibinfo{title}{Deep Learning},
  \bibinfo{publisher}{MIT Press},
  \bibinfo{note}{\url{http://www.deeplearningbook.org}}, \bibinfo{year}{2016}.

\bibitem[{Chauvin and Rumelhart(2013)}]{chauvin2013backpropagation}
\bibinfo{author}{Y.~Chauvin}, \bibinfo{author}{D.~Rumelhart},
  \bibinfo{title}{Backpropagation: theory, architectures, and applications},
  \bibinfo{publisher}{Psychology Press}, \bibinfo{year}{2013}.

\bibitem[{Kingma and Ba(2014)}]{kingma2014adam}
\bibinfo{author}{D.~Kingma}, \bibinfo{author}{J.~Ba}, \bibinfo{title}{Adam: A
  method for stochastic optimization}, \bibinfo{journal}{arXiv preprint
  arXiv:1412.6980} .

\bibitem[{Sobol(1993)}]{sobol1993sensitivity}
\bibinfo{author}{I.~Sobol}, \bibinfo{title}{Sensitivity estimates for nonlinear
  mathematical models}, \bibinfo{journal}{Mathematical modelling and
  computational experiments} \bibinfo{volume}{1}~(\bibinfo{number}{4})
  (\bibinfo{year}{1993}) \bibinfo{pages}{407--414}.

\bibitem[{pyd(2018)}]{pydoe}
\bibinfo{title}{{pyDOE}: The experimental design package for python},
  \urlprefix\url{https://pythonhosted.org/pyDOE/}, \bibinfo{year}{2018}.

\bibitem[{sci(2018)}]{scipy_stats}
\bibinfo{title}{Statistical functions (scipy.stats)},
  \urlprefix\url{https://docs.scipy.org/doc/scipy/reference/stats.html},
  \bibinfo{year}{2018}.

\bibitem[{Abadi et~al.(2015)Abadi, Agarwal, Barham, Brevdo, Chen, Citro,
  Corrado, Davis, Dean, Devin, Ghemawat, Goodfellow, Harp, Irving, Isard, Jia,
  Jozefowicz, Kaiser, Kudlur, Levenberg, Man\'{e}, Monga, Moore, Murray, Olah,
  Schuster, Shlens, Steiner, Sutskever, Talwar, Tucker, Vanhoucke, Vasudevan,
  Vi\'{e}gas, Vinyals, Warden, Wattenberg, Wicke, Yu, and
  Zheng}]{tensorflow2015-whitepaper}
\bibinfo{author}{M.~Abadi}, \bibinfo{author}{A.~Agarwal},
  \bibinfo{author}{P.~Barham}, \bibinfo{author}{E.~Brevdo},
  \bibinfo{author}{Z.~Chen}, \bibinfo{author}{C.~Citro},
  \bibinfo{author}{G.~Corrado}, \bibinfo{author}{A.~Davis},
  \bibinfo{author}{J.~Dean}, \bibinfo{author}{M.~Devin},
  \bibinfo{author}{S.~Ghemawat}, \bibinfo{author}{I.~Goodfellow},
  \bibinfo{author}{A.~Harp}, \bibinfo{author}{G.~Irving},
  \bibinfo{author}{M.~Isard}, \bibinfo{author}{Y.~Jia},
  \bibinfo{author}{R.~Jozefowicz}, \bibinfo{author}{L.~Kaiser},
  \bibinfo{author}{M.~Kudlur}, \bibinfo{author}{J.~Levenberg},
  \bibinfo{author}{D.~Man\'{e}}, \bibinfo{author}{R.~Monga},
  \bibinfo{author}{S.~Moore}, \bibinfo{author}{D.~Murray},
  \bibinfo{author}{C.~Olah}, \bibinfo{author}{M.~Schuster},
  \bibinfo{author}{J.~Shlens}, \bibinfo{author}{B.~Steiner},
  \bibinfo{author}{I.~Sutskever}, \bibinfo{author}{K.~Talwar},
  \bibinfo{author}{P.~Tucker}, \bibinfo{author}{V.~Vanhoucke},
  \bibinfo{author}{V.~Vasudevan}, \bibinfo{author}{F.~Vi\'{e}gas},
  \bibinfo{author}{O.~Vinyals}, \bibinfo{author}{P.~Warden},
  \bibinfo{author}{M.~Wattenberg}, \bibinfo{author}{M.~Wicke},
  \bibinfo{author}{Y.~Yu}, \bibinfo{author}{X.~Zheng},
  \bibinfo{title}{{TensorFlow}: Large-Scale Machine Learning on Heterogeneous
  Systems}, \urlprefix\url{https://www.tensorflow.org/},
  \bibinfo{note}{software available from tensorflow.org}, \bibinfo{year}{2015}.

\bibitem[{Chollet et~al.(2015)}]{chollet2015keras}
\bibinfo{author}{F.~Chollet}, et~al., \bibinfo{title}{Keras},
  \bibinfo{howpublished}{\url{https://keras.io}}, \bibinfo{year}{2015}.

\end{thebibliography}
\end{document}